%% file: uzman_arithmeticity_main.tex
\documentclass[a4paper,11pt,oneside]{amsart}

\usepackage[english]{babel}
\usepackage[T1]{fontenc}
\usepackage[utf8]{inputenc}
\usepackage{textcomp}
\usepackage{adjustbox}
\usepackage{parallel}
\usepackage{stmaryrd}
\usepackage[usenames,dvipsnames]{xcolor}
\usepackage[textwidth=1in,colorinlistoftodos]{todonotes} 
        
\usepackage{latexsym}
\usepackage{float}
\usepackage{mathtools,amscd,amsthm,amssymb,amsmath,wasysym,marvosym}
\usepackage{graphicx}    
\usepackage{xr}
\usepackage{centernot}
\usepackage{hhline}
\usepackage{imakeidx}
\usepackage{euscript}
\usepackage{faktor}
\usepackage[
    marginparwidth = 1in,
    top = 4cm,
    bottom = 4cm,
    right = 1in,
    left = 1in,
    headheight = 15pt
    ]{geometry} 
\usepackage{fancyhdr}    
\usepackage[shortlabels]{enumitem}
    \setlist[description]{leftmargin=\parindent,labelindent=\parindent}
    \setlist[enumerate]{align=left, leftmargin=0pt, listparindent=\parindent, labelwidth=0pt, itemindent=!}
    
\usepackage{comment}
   \specialcomment{AU}{\begingroup\color{MidnightBlue}\em [AU]\;}{\hfill$\Huge\lrcorner$\endgroup}

\PassOptionsToPackage{hyphens}{url}
\usepackage[colorlinks=true,linkcolor=Mahogany,citecolor=Mahogany,urlcolor=Mahogany,backref=page]{hyperref}
\usepackage[final]{showlabels}
\usepackage{tikz-cd}
\usepackage{extarrows}
\usetikzlibrary{matrix,arrows,decorations.pathmorphing}

\newcommand{\fn}{\footnote}

\def\undertilde#1{\mathord{\vtop{\ialign{##\crcr
$\hfil\displaystyle{#1}\hfil$\crcr\noalign{\kern1.5pt\nointerlineskip}
$\hfil\tilde{}\hfil$\crcr\noalign{\kern1.5pt}}}}}
\newenvironment{ncase}{\left.\begin{cases}}{\end{cases}\right\}}

\newcommand{\OL}[1]{\overline{#1}}
\newcommand{\UL}[1]{\underline{#1}}

\newcommand{\blt}{\bullet}
\newcommand{\wt}{\widetilde}

\newcommand{\done}{\hfill$\lrcorner$}

\newcommand{\pt}{\mbox{pt}}
\newcommand{\restr}[2]{{\left.{#1}\right|_{#2}}}

\newcommand{\ot}{\shortleftarrow}

\newcommand{\EE}{\mathbb{E}}

\newcommand{\II}{\mathbb{I}}

\newcommand{\LL}{\mathbb{L}}

\newcommand{\PP}{\mathbb{P}}

\newcommand{\RR}{\mathbb{R}}
    \newcommand{\RRR}{\mathbb{R}_{\geq0}}
    \newcommand{\RRP}{\mathbb{R}_{>0}}

\newcommand{\TT}{\mathbb{T}}

\newcommand{\ZZ}{\mathbb{Z}}
    \newcommand{\ZZR}{\mathbb{Z}_{\geq0}}
    \newcommand{\ZZP}{\mathbb{Z}_{\geq1}}

\newcommand{\mce}{\mathcal{E}}
\newcommand{\mcf}{\mathcal{F}}
\newcommand{\mcg}{\mathcal{G}}

\newcommand{\mcl}{\mathcal{L}}

\newcommand{\mco}{\mathcal{O}}
\newcommand{\mcp}{\mathcal{P}}
\newcommand{\mcq}{\mathcal{Q}}

\newcommand{\mcs}{\mathcal{S}}

\newcommand{\mcu}{\mathcal{U}}

\newcommand{\mcw}{\mathcal{W}}

\newcommand{\eua}{\EuScript{A}}

\newcommand{\euc}{\EuScript{C}}
\newcommand{\eud}{\EuScript{D}}
\newcommand{\eue}{\EuScript{E}}

\newcommand{\euw}{\EuScript{W}}

\newcommand{\mfe}{\mathfrak{e}}

\newcommand{\mfg}{\mathfrak{g}}

\newcommand{\mfC}{\mathfrak{C}}

\newcommand{\mfE}{\mathfrak{E}}

\newcommand{\mfH}{\mathfrak{H}}

\newcommand{\mfR}{\mathfrak{R}}

\newcommand{\mfT}{\mathfrak{T}}

\newcommand{\cat}[1]{\underline{\overline{\mbox{#1}}}}
    \newcommand{\tcat}[1]{\underline{\overline{\mbox{\tiny #1}}}}

\DeclareMathOperator{\Hom}{Hom}

\DeclareMathOperator{\Aut}{Aut}

\DeclareMathOperator{\id}{id}

\DeclareMathOperator{\im}{im}
\DeclareMathOperator{\proj}{proj}

\newcommand{\fwd}[1]{\overrightarrow{\;#1\;\;}}
\newcommand{\bck}[1]{\overleftarrow{\;\;#1\;}}

\DeclareMathOperator{\Sub}{Sub} 

\newcommand{\all}{\forall}
\newcommand{\exi}{\exists}
\newcommand{\ifr}{\Rightarrow}

\renewcommand{\iff}{\Leftrightarrow}

\DeclareMathOperator{\pwr}{\EuScript{P}}

\DeclareMathOperator{\open}{\EuScript{T}}

\DeclareMathOperator{\bor}{\EuScript{B}}

\DeclareMathOperator{\Par}{Par}

    \DeclareMathOperator{\mPar}{mPar}
    \DeclareMathOperator{\aePar}{\aem\hspace{-5pt} \Par}

    \DeclareMathOperator{\aemPar}{\aem\hspace{-5pt} \mPar}
\DeclareMathOperator{\Discrete}{\mathbf{0}}
\DeclareMathOperator{\Indiscrete}{\mathbf{1}}
\DeclareMathOperator{\Hull}{Hull} 

\DeclareMathOperator{\Bij}{Bij}

\DeclareMathOperator{\Gr}{Gr}
\DeclareMathOperator{\Spl}{Spl}

\DeclareMathOperator{\GL}{GL}
	\DeclareMathOperator{\DGL}{DGL}
\DeclareMathOperator{\Aff}{Aff}
	    \newcommand{\aff}{\text{aff}}
		\DeclareMathOperator{\DAff}{DAff}

\DeclareMathOperator{\lie}{lie}

\DeclareMathOperator{\Span}{Span}

\let\oldae\ae

    \renewcommand{\ae}{\oldae{ }\hspace{-3pt}}
    \newcommand{\aem}{\mbox{\small\oldae{ }}}
\newcommand{\dbar}{\mbox{\dh}}    
\DeclareMathOperator{\leb}{leb}
\DeclareMathOperator{\haar}{haar}

\newcommand{\loc}{\mbox{\tiny loc}}

\newcommand{\eps}{\varepsilon}

\DeclareMathOperator{\supp}{supp}

    \DeclareMathOperator{\Radon}{Radon}
	\DeclareMathOperator{\vRadon}{vRadon}
    \DeclareMathOperator{\Prob}{Prob}

\DeclareMathOperator{\Jac}{Jac}

\DeclareMathOperator{\Ad}{Ad}

\DeclareMathOperator{\Imm}{Imm} 
    \DeclareMathOperator{\Emb}{Emb} 

\DeclareMathOperator{\Fol}{Fol}
	    \DeclareMathOperator{\aeFol}{\aem\hspace{-5pt} \Fol}

\DeclareMathOperator{\rank}{rank}

\DeclareMathOperator{\Fix}{Fix}

\DeclareMathOperator{\Diff}{Diff}

\DeclareMathOperator{\Stab}{Stab}

\newcommand{\lact}{\curvearrowright}

\DeclareMathOperator{\ent}{ent}

\DeclareMathOperator{\LSpec}{LSpec}

\DeclareMathOperator{\Cham}{Cham}
\DeclareMathOperator{\Wall}{Wall}

\DeclareMathOperator{\Osel}{Osel}

\DeclareMathOperator{\Friedent}{Fried}

\newtheoremstyle{main}
    {7pt}
    {7pt}
    {}
    {}
    {\bfseries}
    {:}
    {.5em}
    {} 
    
\newtheoremstyle{sub}
    {7pt}
    {7pt}
    {} 
    {\parindent} 
    {\bfseries} 
    {:} 
    {.5em}
    {}

\theoremstyle{main}  
    \newtheorem{dfn}{\textbf{Definition}}
        
    \newtheorem*{pf}{\textbf{Proof}}

    \newtheorem{rem}{\textbf{Remark}}
        
    \newtheorem{obs}{\textbf{Observation}}

    \newtheorem{thm}{\textbf{Theorem}}
        
    \newtheorem{cor}{\textbf{Corollary}}
        
    \newtheorem{prp}{\textbf{Proposition}}
        
    \newtheorem{lem}{\textbf{Lemma}}

\theoremstyle{sub}

\begin{document}

\title{Arithmeticity for Smooth Maximal Rank Positive Entropy Actions of $\RR^k$}

\author{Alp Uzman}
\address{Department of Mathematics, University of Utah\\
155 South 1400 East, JWB 233\\
Salt Lake City, UT 84112\\
USA}
\email{uzman@math.utah.edu}


\begin{abstract}
We establish arithmeticity in the sense of \cite{MR3503686} of  higher rank abelian actions $\alpha_\blt:\RR^k\lact M$ by diffeomorphisms on an anonymous manifold $M$ of dimension $\dim(M)=2k+1$ provided that there is an ergodic invariant Borel probability measure $\mu$ on $M$ w/r/t which each nontrivial time-$t$ map $\alpha_t$ of the action has positive entropy. Arithmeticity in this context means that the action $\alpha$ is measure theoretically isomorphic to a constant time change of the suspension of an affine Cartan action of $\ZZ^k$.  This in particular solves, up to measure theoretical isomorphism, Problem 4 from \cite[p.394]{MR2811602}.
\end{abstract}

\maketitle

\tableofcontents

\input{uzman_arithmeticity_introduction.tex}

\input{uzman_arithmeticity_outline.tex}

\input{uzman_arithmeticity_preliminaries.tex}

\input{uzman_arithmeticity_mrpehypotheses.tex}

\input{uzman_arithmeticity_affinestructures.tex}

\input{uzman_arithmeticity_affineholonomies.tex}

\input{uzman_arithmeticity_affineextension.tex}

\input{uzman_arithmeticity_homoclinicgroup.tex}

\bibliographystyle{amsalpha}
\bibliography{uzman_arithmeticity_references.bib}

\end{document}

%% file: uzman_arithmeticity_introduction.tex
\section{Introduction}

\subsection{Statement of the Main Results}

We prove the following arithmeticity result. Let $M$ be a compact $C^\infty$ manifold, $k\in\ZZ_{\geq1}$, $\theta\in]0,1]$, $r=(1,\theta)$, $\mfE^r(\RR^k\lact M)$ be the \textbf{ergodic theory with no potential} of $\RR^k$ systems on $M$ of class $C^r$; by definition this is the collection of all pairs $(\mu,\alpha)$ such that $\mu$ is a Borel probability measure on $M$ and $\alpha_\blt:\RR^k\to\Diff^r(M)$ is a group homomorphism such that $\alpha:\RR^k\times M\to M$ is $C^r$ and $\mu$ is $\alpha$-invariant; we call any such pair $(\mu,\alpha)$ an \textbf{$\RR^k$ system}. We show that in the case of simplest Lyapunov geometry compatible with hyperbolicity, any system is the suspension of a hyperbolic algebraic action on some space crystal:

\begin{thm}\label{001}

Let $(\mu,\alpha)\in\mfE^r(\RR^k\lact M)$. If

\begin{itemize}

\item $k\in\ZZ_{\geq2}$ and $\dim(M)=2k+1$,

\item $(\mu,\alpha)$ is locally free and ergodic, 

\item The system $(\mu,\alpha)$ has exactly $k+1$ distinct Lyapunov hyperplanes, and

\item $\all t\in\RR^k\setminus 0: \mfe_{(\mu,\alpha)}(t)=\ent_\mu(\alpha_t)>0$,

\end{itemize}

then there is

\begin{itemize}

\item an affine Cartan action $\gamma_\blt:\ZZ^k\to \Aff(T^{k+1})$, where $T^{k+1}$ is the torus or the $\pm$-infratorus of dimension $k+1$, and 

\item a $\kappa\in\GL(k,\RR)^\circ$

\end{itemize}

such that 

$$\exi\Phi_{(\mu,\alpha)}:(\mu,\alpha)\xrightarrow{\phantom{\cong_{\tcat{Meas}}}\cong_{\tcat{Meas}}\phantom{\cong_{\tcat{Meas}}}} (\haar_{\TT^k}\otimes^\gamma \haar_{T^{k+1}},\hbar^{\gamma}_\kappa),$$

where $\hbar^{\gamma}_\kappa$ is the suspension of $\gamma$ with a constant time change $\kappa$.  Furthermore, 

\begin{itemize}

\item The restriction of the measure theoretical isomorphism $\Phi_{(\mu,\alpha)}$ to any global stable manifold of any Weyl chamber of $(\mu,\alpha)$ is $C^r$, and

\item For any $\theta'\in]0,\theta[$ there is a open subset $U_{\theta'}\subseteq M$ with $\mu(M\setminus U_{\theta'})<\theta'$ and the measure theoretical isomorphism $\Phi_{(\mu,\alpha)}$ extends to a $C^{r-\theta'}$ injective immersion on $U_{\theta'}$.

\end{itemize}

\done
\end{thm}

Here the infratorus is defined in terms of the standard action of $\Aff(\RR^{k+1})$ of affine automorphisms of $\RR^{k+1}$ on $\RR^{k+1}$: we have that $\Aff(\RR^{k+1})\cong \RR^{k+1}\rtimes \GL(k+1,\RR)$. Then the $(k+1)$-dimensional torus $\TT^{k+1}$ is the orbit space of the subgroup $\ZZ^{k+1}\rtimes \{I_{k+1}\}\leq \Aff(\RR^{k+1})$ and the $(k+1)$-dimensional \textbf{$\pm$-infratorus} is the orbit space of the subgroup $\ZZ^{k+1}\rtimes \{\pm I_{k+1}\}\leq \Aff(\RR^{k+1})$. In general an \textbf{infratorus} is the orbit space of a subgroup $L\rtimes F\leq \Aff(\RR^{k+1})$ for $L$ a lattice and $F$ a finite group. An \textbf{affine Cartan action} $\gamma$ is an action by affine automorphisms of $T^{k+1}$ such that for any $t\in\ZZ^k\setminus0$, the linear part of the time-$t$ map $\gamma_t$ is ergodic w/r/t the Haar measure\fn{One could equivalently define $\gamma$ to be an affine Cartan action iff the linear part of $\gamma_t$ is a hyperbolic (infra-)toral automorphism by \cite[p.731, Prop.4.1]{MR1949111}. Note that for $k=1$, ergodicity is strictly weaker than hyperbolicity. Also note that in the literature there are nonequivalent definitions of Cartan actions; see e.g. \cite[pp.4-5]{spatziervinhage}. Our definition is compatible with \cite{MR2643892} and is an extension of the definition in \cite[p.731]{MR1949111}.}.

\begin{rem}\label{087}

Some hypotheses in \autoref{001} can be weakened. More specifically, if $M$ is not necessarily compact, $r=(q,\omega)$ for some $q\in\ZZP$, and for some modulus of continuity $\omega$ satisfying the Dini condition, so that the $q$th differential objects have $\omega$ as a local modulus of continuity, or $r=\infty$,  $\alpha_\blt:\RR^k\to \Diff^r(M)$ is an arbitrary family of $C^r$ diffeomorphisms with $k$ weakly generating vector fields whose norms w/r/t  some $C^0$ fiberwise norm are $\log^+$-$(k,1)$-Lorentz w/r/t $\mu$ and that commute on a set of full $\mu$ measure, and if $\alpha$ is not necessarily $\mu$-essentially locally free as a group action in the category of standard probability spaces, all other hypotheses ditto (in particular ergodicity is certainly necessary), then the theorem is still true. For the sake of readability we don't prove this generalization (although see \autoref{138} for an indication as to the kinds of arguments that can be used to handle higher regularity).

\done
\end{rem}

\begin{rem}\label{101}

There is also some redundancy in the way the hypotheses of \autoref{001} are stated and the above statement corresponds to the most inefficient choice, see \autoref{066} for a related discussion.

\done
\end{rem}

\begin{dfn}\label{083}

Let us call a system $(\mu,\alpha)$ satisfying the four hypotheses of \autoref{001} an \textbf{$\RR^k$ maximal rank positive entropy system} (\textbf{MRPES} for short). Similarly, if $N$ is a compact $C^\infty$ manifold, $(\nu,\beta)\in\mfE^r(\ZZ^k\lact N)$ is an ergodic system whose suspension system  $(\haar_{\TT^k}\otimes^\beta \nu,\hbar^{\beta})$ is an MRPES, let us call $(\nu,\beta)$ a \textbf{$\ZZ^k$ maximal rank positive entropy system}.

Moreover, we call $\mfe_{(\mu,\alpha)}:\RR^k\to \RRR$ the \textbf{entropy gauge} of the system $(\mu,\alpha)$.

\done
\end{dfn}

Many corollaries follow immediately from \autoref{001}; below is a sampler:

\begin{cor}\label{004}
Any $\RR^k$ MRPES is the suspension of a $\ZZ^k$ MRPES up to a measure theoretical isomorphism and a constant linear time change. In particular, the answer to Problem 4 of \cite{MR2811602} is negative up to $\cat{Meas}$-isomorphism.

\done
\end{cor}

\begin{cor}\label{086}

The Jacobian cocycle of any $\RR^k$ MRPES along any Lyapunov \ae-foliation is cohomologous to a cocycle constant in space with measurable transfer that is $C^r$ along any Lyapunov \ae-foliation. In particular, this solves Conjecture 1\fn{This conjecture ( \cite[p.393]{MR2811602}) claims that for any $\ZZ^k$ or $\RR^k$ MRPES the derivative cocycle along any Lyapunov \ae-foliation is measurably cohomologous to one constant in space. The conjecture for $\ZZ^k$ MRPES follows from \cite[p.394]{MR2811602}.} of \cite{MR2811602} for the $\RR^k$ case.

\done
\end{cor}

\begin{cor}\label{096}

Any Lyapunov exponent of any time-$t$ map of any $\RR^k$ MRPES is the logarithm of an algebraic number.

\done
\end{cor}

\begin{cor}\label{073}

If $k\in\ZZ_{\geq2}$, $(\mu,\alpha)\in\mfE^r(\RR^k\lact M)$ is an essentially locally free and ergodic system with exactly $k+1$ distinct Lyapunov hyperplanes, $\dim(M)=2k+1$ and for some $t^\ast\in\RR^k: (\mu,\alpha_{t^\ast})\in \mfE^r(\ZZ\lact M)$ has the $K$-property, then the Lyapunov hyperplanes of $(\mu,\alpha)$ can't be in general position, and there is some $t^\dagger\in\RR^k\setminus0$ such that $\mfe_{(\mu,\alpha)}(t^\dagger)=0$.

\done
\end{cor}

\begin{cor}\label{074}

Let $k\in\ZZ_{\geq2}$ and $(\mu,\alpha)\in\mfE^r(\RR^k\lact M)$ be an essentially locally free ergodic system with exactly $k+1$ distinct Lyapunov hyperplanes. If $\all t\in\RR^k\setminus0: \mfe_{(\mu,\alpha)}(t)>0$ and for some $t^\ast\in \RR^k$,  $(\mu,\alpha_{t^\ast})\in \mfE^r(\ZZ\lact M)$ has the $K$-property, then $\dim(M)\geq 2k+1$.

\done
\end{cor}

\begin{cor}\label{144}

Let $(\mu,\alpha)\in\mfE^r(\RR^k\lact M)$ be an MRPES. Then the set 

$$\{t\in\RR^k \,|\, (\mu,\alpha_t)\in\mfE^r(\ZZ\lact M)\text{ is not ergodic}\}$$

of non-ergodic time-$t$ systems of $(\mu,\alpha)$ is the union of exactly countably many hyperplanes that constitute a mesh determined by an injective group homomorphism $\ZZ^k\hookrightarrow\RR^k$.

\done
\end{cor}

\subsection{History and Context}

The strategy we follow to prove \autoref{001} is to adapt the machinery developed in \cite{MR3503686} from $\ZZ^k$ systems to $\RR^k$ systems (along the way we will also provide significantly more details and address various gaps).  Said machinery is built on previous advances on the geometric approach to measure rigidity of actions of $\RR^{k_1}\oplus \ZZ^{k_2}$ ($k_1+k_2\in\ZZ_{\geq2}$) by automorphisms of homogeneous spaces and more generally diffeomorphisms of manifolds, with some hyperbolicity in either case. In particular this paper is a continuation of the chain of work \cite{MR1406432,MR1619571,MR1858547,MR2261075, MR2285730, MR2643892, MR2811602,MR2729332,MR3503686}\fn{Note that, as of the writing of this paper, \cite{MR3666074} is the latest addition to this chain of papers. See \cite{MR2457050} for a survey of the papers preceding \cite{MR3503686}. The survey \cite{MR2191215} contextualizes the geometric method in the broader study of the ergodic theory of $\RR^{k_1}\oplus \ZZ^{k_2}$ actions, whereas the survey \cite{rodriguezhertzf.2021} contextualizes the works in question as part of late Prof. Katok's life's work. Finally \cite[Sec.3]{MR2261070} surveys the connections between measure rigidity and other areas.}.

\subsection{Discrete v. Connected Time}

For comparison purposes it is useful to recall the $\ZZ^k$ arithmeticity theorem from \cite{MR3503686}:

\begin{thm}[Katok - Rodr\'{i}guez Hertz\fn{\cite[pp.137-138, Thm.1]{MR3503686}}]\label{022}

Let $N$ be a compact $C^\infty$ manifold. If $(\nu,\beta)\in\mfE^r(\ZZ^k\lact N)$ is an MRPES, then there is

\begin{itemize}

\item an affine Cartan action $\delta_\blt:\ZZ^k\to \Aff(T^{k+1})$, and 

\item an injective group homomorphism $j:\ZZ^k\hookrightarrow\ZZ^k$ with $F=\faktor{\ZZ^k}{\fwd{j}(\ZZ^k)}$ finite

\end{itemize}

such that 

$$\exi \Phi_{(\nu,\beta)}:(\nu,\beta) \xrightarrow{\phantom{\cong_{\tcat{Meas}}}\cong_{\tcat{Meas}}\phantom{\cong_{\tcat{Meas}}}} (\haar_{F}\otimes^\delta \haar_{T^{k+1}},\hbar^\delta).$$

Furthermore, 

\begin{itemize}

\item The restriction of the measure theoretical isomorphism $\Phi_{(\nu,\beta)}$ to any global stable manifold of any Weyl chamber of $(\nu,\beta)$ is a $C^r$ diffeomorphism, and

\item For any $\theta'\in]0,\theta[$, there is a open subset $U_{\theta'}\subseteq N$ with $\nu(N\setminus U_{\theta'})<\theta'$ and the measure theoretical isomorphism extends to a $C^{r-\theta'}$ injective immersion on $U_{\theta'}$.

\end{itemize}

\done
\end{thm}

\begin{dfn}\label{039}
We call $\Phi_{(\mu,\alpha)}$ of \autoref{001} the \textbf{arithmeticity isomorphism} of $(\mu,\alpha)$. Similarly $\Phi_{(\nu,\beta)}$ of \autoref{022} is the \textbf{arithmeticity isomorphism} of $(\nu,\beta)$. 

\done
\end{dfn}

Note that the positive entropy hypothesis of \autoref{022} is formulated in terms of the elements of the suspension system, since the positive entropy along "irrational directions" is crucial for the proof. See \autoref{023} for a characterization of the positive entropy hypothesis in terms of Lyapunov geometry. There are essentially two new complications in the proof of \autoref{001} compared to the proof of \autoref{022}; both due to the presence of the orbit directions. The first is that in the proof of $\ZZ^k$ arithmeticity by a classical theorem of Pesin\fn{\cite[p.94, Thm.7.9]{MR0466791}} the finite time crystal can be separated at the first instance (this is their "weak mixing reduction"); this classical theorem of Pesin does not work for $\RR^k$ actions (it has indeed an analog for $\RR^k$ actions, which provides only an alternative in terms of the eigenfunctions of the Koopman action\fn{\cite[p.1227, Thm.9.7]{MR0488169}; compare \cite[p.29, Thm.14]{MR0242194}.}). The second complication is that in the $\RR^k$ case the Katok - Rodr\'{i}guez Hertz machinery may a priori spill over in the orbit directions. In more concrete terms, these make the biggest difference toward the end in \autoref{135} and \autoref{069}, where we establish the suspension structure. In the $\ZZ^k$ case the homoclinic group is recognized to be a lattice in an abelian group; in our case it is recognized to be a lattice in a solvable group.

\subsection{Uniform v. Non-Uniform Normal Hyperbolicity to the Orbit Foliation}

Our main theorem \autoref{001} can also be seen as a non-uniformly normally hyperbolic (to the orbit foliation) analog of a result of Matsumoto on uniformly normally hyperbolic (to the orbit foliation):

\begin{thm}[Matsumoto\fn{ \cite[p.42,Thm.3.1]{MR1939183}}]\label{005}

Let $M$ be a closed oriented $C^\infty$ manifold, $k\in\ZZ_{\geq2}$, $\alpha_\blt:\RR^k\to\Diff^\infty_+(M)$ be a $C^\infty$ action. If

\begin{itemize}

\item $\alpha$ is locally free,

\item $\dim(M)=2k+1$,

\item There is an $\Ad^\alpha$-invariant $C^0$ splitting $TM=O\oplus \bigoplus_{i\in\OL{k+1}} E^i$, where $O$ is the subbundle tangent to the orbit foliation of $\alpha$ and each $E^i$ is a topological line bundle with the property that there is a $\xi_\blt: \OL{k+1}\to \RR^k$ such that for any $i\in\OL{k+1},$, the time-$\xi_i$ map $\alpha_{\xi_i}$ contracts $E^i$ exponentially fast and expands  $\bigoplus_{\substack{j\neq i}} E^j$ exponentially fast,

\end{itemize}

then there is

\begin{itemize}

\item an affine Cartan action $\gamma_\blt:\ZZ^k\to \Aff_+(\TT^{k+1})$, and

\item $\kappa\in \GL(k,\RR)^\circ$

\end{itemize}

such that

$$\exi\Phi_{\alpha}:\alpha\xrightarrow{\phantom{\cong}\cong_{C^\infty}\phantom{\cong}} \hbar^\gamma_\kappa,$$

where again $\hbar^\gamma_\kappa$ is the suspension of $\gamma$ with  a constant time change $\kappa$.
\done
\end{thm}

Matsumoto calls an action of $\RR^k$ satisfying the hypotheses of \autoref{005} a "split Anosov action". The $C^\infty$ setting is amenable to differential topological methods; as such despite the fact that Matsumoto's theorem looks formally very similar to our main theorem (partially also due to the author's paraphrasing) the methods involved are completely different.

\subsection{Time Crystal Interpretation\label{142}}

In \cite[p.596]{MR1858547} and \cite[p.510]{MR1898802}\fn{Also see \cite[p.583]{MR2342699}.} the $k$-torus over which the suspension of a $\ZZ^k$ action fibers is called a "'time' torus". We take this nomenclature more seriously and borrowing a name from contemporary physics\fn{ See \cite{PhysRevLett.109.160401,doi:10.1063/PT.3.4020,ynlz}.} we call the base torus of the suspension constructed in \autoref{001} the \textbf{time crystal} and the fiber torus or $\pm$-infratorus of the suspension the \textbf{space crystal} of the system $(\mu,\alpha)$. We also call the total space of the suspension the \textbf{spacetime crystal presentation} of $(\mu,\alpha)$. This evocative language is meant to simplify referring to the similar algebraic structures along different directions. Further we may interpret our result as saying that up to a measure-preserving coordinate change the time crystal of $(\mu,\alpha)$ does keep track of the broken continuous time symmetry, all the while being part of the state space $M$: a part of the state space functionally works as an internal clock which keeps track of the discrete fiber dynamics. Analogously the time crystal of \autoref{022} is the finite factor group $F$. The reader ought to be warned that this use of the phrase "time crystal" does \textit{not} match with the use of the term in the physics literature\fn{The relation between the two uses is roughly the same relation between Pugh's and Anosov's closing lemmas (see \cite{MR3059636}).}.

\subsection{Miscellaneous Notation}

Here we collect some miscellaneous notation; other notation we'll use will be introduced when appropriate in context, if its meaning is not clear from context. We repeat and reintroduce some notation in-text.

Overall $\blt$ signifies and emphasizes the first priority place for an argument, and for $p\in\ZZP$ we put $\UL{p}=\{0,1,...,p-1\}=\OL{p}-1$, $\UL{p}+1=\{1,2,...,p\}=\OL{p}$, $\OL{-p}=\{-1,-2,...,-p\}=-\OL{p}$ and $\UL{-p}=\{0,-1,...,-p+1\}=-\UL{p}$.

$\pwr(X),\bor(X),\open(X)$ stand for the collection of all subsets, all (Borel) measurable subsets, and all open subsets of a set, measurable space, and topological space $X$, respectively. 

For $(M,d)$ a metric space, $x\in M$ and $r\in\RRP$: $M[x|\leq r]$ denotes the closed ball in $M$ centered at $x$ with radius $r$; $M[x|<r]$ similarly denotes the open ball.

For $(X,\mu)$ a measure space, $\all x\in_\mu X$ means that the point $x$ is chosen from an appropriate full $\mu$-measure subset of $X$. For two measurable subsets $A,B\subseteq X$, $A\subseteq_\mu B$ means that $\mu(A\setminus B)=0$ and $A =_\mu B$ means $\mu(A\triangle B)=0$. Similarly if $f,g$ are two measurable functions with \ae-domain $X$, $f=_\mu g$ means that $\all x\in_\mu X: f(x)=g(x)$. The ligature \ae ("ash") stands for a measure whose name is suppressed; it also stands for the phrase "almost everywhere"; when in the form $\mu$-\ae the former usage is dropped. $L^0$ means the space of measurable functions \ae-defined and \ae-identified if a measure is fixed.

Let $P=\prod_{\alpha\in A}P_\alpha$ be a set of tuples of numbers or more general objects, where $A$ is an anonymous indexing set, $L:P\to \RR$ and $R:P\to \RR$ be two functions. For two subsets $Q_1,Q_2\subseteq P$ and $p\in P$ we write $L(p) \lesssim_{Q_1}^{Q_2} R(p)$ if for some function $C:P\to \RRP$ that is possibly nonconstant on $Q_1$ and certainly constant on $Q_2$, for any $p\in P$, we have $L(p)\leq C(p)R(p)$. If $Q_1=P_\alpha$, $Q_2=P_\beta$ we also write $L(p)\lesssim_{p_\alpha}^{p_\beta} R(p)$. Here neither of $Q_1,Q_2$  needs to be optimal. 

$\cat{Lie},\cat{lie},\cat{Man}^s,\cat{Man}^s_\rho,\cat{sMble},\cat{sMeas},\cat{sProb}$ stand for the categories of Lie groups, Lie algebras, $C^s$ manifolds, $C^s$ $\rho$-manifolds, standard measurable spaces, standard measure spaces, and standard probability spaces with the standard choices for arrows. The latter three "Borelesque" categories are often the associated \ae-factor categories, signified by \ae{} in diagrams.  We'll conflate a category and the category of $G$-objects (as well as $G$-systems) of that category for $G$ a group object in the category; accordingly instead of "conjugacy" we'll say "isomorphism". For $f$ a morphism in a category $\fwd{f}$ stands for the induced map under a covariant functor and $\bck{f}$ stands for the induced map under a contravariant functor. For any morphism object $f$, the notation $f^{-1}$ will only be used when it can be interpreted as an object of the same type as $f$. 

$\Imm^s$ and $\Emb^s$ stand for $C^s$ injective immersions and embeddings, respectively; all such function spaces are topologized via the uniform $C^s$ topology on compact subsets; we'll denote the associated chosen distance by $d_{C^s}$. 

For semidirect products we always denote the group $G$ that fits into a split short exact sequence $L\to G\to R$ (in any group object category) by $G\cong L\rtimes R$ or $G\cong R\ltimes L$. Thus $L$ is a normal subgroup in $G$ and $R$ acts on $L$.

\subsection{Acknowledgements}

This work is part of the author's PhD thesis at the Pennsylvania State University\fn{\url{https://etda.libraries.psu.edu/catalog/20180azu4}}. The results in this paper were first announced at the  "Pennsylvania State University Semi-annual Workshop in Dynamical Systems and Related Topics" held on October 28 - 31, 2021. The results were also presented at the poster session at 
Simons Center for Geometry and Physics workshop "Flexibility and rigidity in dynamical systems" held  on March 7-11, 2022. The author would like to thank his thesis advisor Prof. Federico Rodr\'iguez Hertz countless hours of discussion and limitless support. The author also thanks Prof. Kurt Vinhage; the discussions during the author's visit to the University of Utah in early December 2022 have been very useful. This work has been partially supported by the NSF and the Anatole Katok Center for Dynamics and Geometry. Commutative diagrams were drawn using the tikzcd editor\fn{\url{https://tikzcd.yichuanshen.de/}} of Y. Shen.

%% file: uzman_arithmeticity_outline.tex
\section{Outline of the Argument\label{065}}

In this section we give an outline of the proof of \autoref{001}. Our starting point is that by \cite[p.363, Main Thm.(2)]{MR2811602} $\mu$ is absolutely continuous w/r/t the Lebesgue measure class of $M$. Throughout we roughly follow the outline of \cite{MR3503686}.

\begin{description}

\item[Nonstationary Linearizations for Lyapunov \ae-Foliations]

For any nonorbital Lyapunov exponent $\chi$ of $(\mu,\alpha)$, denote by $\mcl^\chi$ the associated one dimensional \ae-foliation and by $L^\chi\leq TM$ the corresponding tangent \ae-subbundle. We first introduce an \ae-unique family of nonstationary linearizations

$$\all x\in_\mu M:\Lambda_x^\chi: L_x^\chi\to \mcl_x^\chi$$

that depend measurably on the basepoint and is $C^r$ if the basepoint is fixed. Within any Lusin-Pesin set the linearizations depend uniformly continuously on the basepoint in the $C^r$ topology.

\item[Affine Structures for Stable and Unstable \ae-Foliations]

For any Weyl chamber $\euc$ of $(\mu,\alpha)$ denote by   $\mcs^\euc,\mcu^\euc,\mco,\mco\mcs^\euc,\mco\mcu^\euc$ the associated stable, unstable, orbit, orbit-stable and orbit-unstable \ae-foliations and by $S^\euc, U^\euc, O, OS^\euc,OU^\euc\leq TM$ the corresponding tangent \ae-subbundles. Then the nonstationary linearizations $\Lambda^\chi_\blt$ in the previous item can be uniquely assembled into $C^r$ affine manifold structures on global manifolds

\begin{align*}
\all x\in_\mu M:
\Sigma_x^\euc &: S_x^\euc\to \mcs_x^\euc, \\
\Gamma\Sigma_x^\euc &: \RR^k\times S_x^\euc\to \mco\mcs_x^\euc,\\
\Upsilon_x^\euc &: U_x^\euc\to \mcu_x^\euc,\\
\Gamma\Upsilon_x^\euc &: \RR^k\times U_x^\euc\to\mco\mcu_x^\euc,
\end{align*}

where the fibers of stable and unstable \ae-subbundles are endowed with the natural affine manifold structure, and  $\RR^k$ is endowed with the Euclidean affine manifold structure.

Furthermore we have a diagonal property: for $x\in_\mu M$ if $(\chi,\euc)=-$ the $1$-dimensional global Lyapunov manifold $\mcl_x^\chi$ inside $\mcs_x^\euc\subseteq\mco\mcs_x^\euc$ corresponds to a coordinate axis in $\RR^s\leq \RR^k\times \RR^s$ and similarly for global unstable and orbit-unstable manifolds. Consequently the global manifolds admit unique $C^r$ $\eua$-affine manifold structures, where $\eua$ is the subgroup of the affine group of appropriate dimension with the property that each element of $\eua$ has diagonal linear part.

\item[Affine Holonomies]

Next we show that for any Weyl chamber $\euc$ of $(\mu,\alpha)$, for any $x\in_\mu M$, and for any $y\in_{\aem} \mcu_x^\euc$, there is a well-defined holonomy map 

$$\mcu_{y\ot x}^\euc: \mco\mcs_x^\euc\to  \mco\mcs_y^\euc$$

along the unstable \ae-foliation associated to $\euc$ with the property that w/r/t the affine structures $\Gamma\Sigma_x^\euc$ and $\Gamma\Sigma_y^\euc$, the holonomy $\mcu_{y\ot x}^\euc$ is affine with linear part of the form $I_k\times D_s$, where $I_k$ is the $k\times k$ identity matrix and $D_s=D_s(\euc,y\leftarrow x)$ is a diagonal $s\times s$ matrix with $s=s(\euc)$ the rank of the \ae-subbundle $S^\euc$. Similarly for any $z\in_{\aem} \mcs_x^\euc$ there is a well-defined holonomy map $\mcs_{z\ot x}^\euc$. The two families of holonomy cohere, in the sense that for any $x\in_\mu M$ and for any $(y,z)\in_{\aem} \mcs_x^\euc\times \mcu_x^\euc \subseteq \mco\mcs_x^\euc \times \mco\mcu_x^\euc$:

$$\mcu_{z\ot x}^\euc(y)=\mcs_{y\ot x}^\euc(z).$$

\item[Measurable Covering Map and Diagonal Affine Extension]

For any Weyl chamber $\euc$ and for $x\in_\mu M$, the affine structures and the holonomy maps assemble into a measurable map $\Phi_x=\Phi_x^\euc: (T_xM,0)\to (M,x)$ that we call the measurable covering map of $(\mu,\alpha)$. $\Phi_x$ has a certain diagonal property in that it is compatible with the orbits and the one dimensional Lyapunov subspaces. Further, the image of $\Phi_x$ is a full $\mu$-measure subset of $M$, $\left(\Phi_x\right)^{-1}$ carries the measure $\mu$ to a Haar measure on $T_xM$, and  locally $\Phi_x$ is a measure theoretical isomorphism.

The change of the basepoint $x$ for the measurable covering map $\Phi_x$ contributes a pullback by a diagonal affine automorphism. We may change the basepoint either vertically along the tangent space $T_xM$ or horizontally along the orbits of $\alpha$; these two types of changes are intimately related. By thinking of $TM$ as an \ae-bundle over $M$ we obtain a diagonal affine extension $\RR^k\lact TM$ with $\Phi:TM\to M$, $(x,v)\mapsto \Phi_x(v)$ the factor map.

\item[Homoclinic Group and the Construction of the Algebraic Model]

Thinking of diagonal affine isomorphisms between different tangent spaces as the symmetries of the measurable covering maps $\Phi_\blt$ gives the homoclinic groupoid of $(\mu,\alpha)$, and specifying a base point $x$ $\in_\mu$ $ M$ gives the homoclinic group $\mfH_x\leq \Aff(T_xM)$. 

We show that $\mfH_x\cong (\ZZ^{k+1}\rtimes F) \rtimes \ZZ^k$, where $F$ is a group of involutions of order at most two, and consequently factoring it out induces the isomorphism $\Phi_{(\mu,\alpha)}$ that transforms $\alpha$ to the suspension of an affine Cartan action and $\mu$ to the suspension measure induced by Haar measure on $T^{k+1}$. Finally smoothness of $\Phi_{(\mu,\alpha)}$ along the stable and unstable \ae-foliations of any Weyl chamber $\euc$ of $(\mu,\alpha)$ as well as smoothness in the sense of Whitney is by a straightforward application of a  Journ\'{e} lemma. 

\end{description}

%% file: uzman_arithmeticity_preliminaries.tex
\section{Preliminaries\label{071}}

In this section we go over the preliminaries. This section contains no proofs as all the material is fairly straightforward except possibly the notation or formalism. The reader is advised to skip this section and return to it when a clarification of definitions or notations is needed.

\subsection{Manifolds with Geometric Structures}

Here we include some basic definitions from the theory of geometric structures\fn{See \cite{goldmanw.m.2021,MR1435975,MR2249478} for versions more appropriate for geometry.}; as far as the author is aware the detailed specific definitions we need in this paper are not available in the literature.

\begin{dfn}\label{036}
Let $X$ be a  $C^\infty$ manifold, $G$ be a Lie group, $s\in\ZZ_{\geq1}\times [0,1]$, $\rho_\blt: G\to \Diff^s(X)$ be a $C^s$ action. A \textbf{local $C^s$ $\rho$-diffeomorphism}  of $X$ is a $C^s$-diffeomorphism $\psi:A\to B$, where $A$ and $B$ are open subsets of $X$, such that there is some $g\in G$ with $\psi=\restr{\rho_g}{A}$.

Let $Y$ be a $C^s$ manifold. A \textbf{$C^s$ $\rho$-structure} (or \textbf{$(G,X)$-structure}, or \textbf{$(G\stackrel{\rho}{\lact} X)$-structure}) on $Y$ is an open cover $\mcu$ of $Y$ and for each $U\in \mcu$  a $C^s$ diffeomorphism $\phi\in \Diff^s(U,\fwd{\phi}(U))$ with $\fwd{\phi}(U)\subseteq X$ open such that for any $U_1,U_2\in\mcu$, 

$$\mfR_{U_2\leftarrow U_1}=\restr{\phi_2\circ\phi_1^{-1}}{\fwd{\phi_1}(U_1\cap U_2)}: \fwd{\phi_1}(U_1\cap U_2)\to \fwd{\phi_2}(U_1\cap U_2)$$

is a local $C^s$ $\rho$-diffeomorphism of $X$. $Y$ is a \textbf{$C^s$ $\rho$-manifold} if it's endowed with a maximal $C^s$ $\rho$-structure. A pair $(U,\phi)$ as above is a \textbf{$C^s$ $\rho$-chart}. We call a $C^s$ $\rho$-chart $(U,\phi)$ \textbf{global} if $U=Y$. A \textbf{$C^s$ globally-$\rho$-manifold} is a $C^s$ $\rho$-manifold with all $C^s$ $\rho$-charts global.

If $Y$ is a $C^s$ $\rho$-manifold, then a \textbf{$C^s$ $\rho$-diffeomorphism} of $Y$ is a $C^s$ diffeomorphism of $Y$ that is a local $C^s$ $\rho$-diffeomorphism of $X$ in any $C^s$ $\rho$-chart. Let us denote by $\Diff^s_\rho(Y)$ the subgroup of $\Diff^s(Y)$ consisting of $C^s$ $\rho$-diffeomorphisms of $Y$.

\done
\end{dfn}

\begin{dfn}\label{038}
A \textbf{$C^s$ affine manifold} is a $C^s$ $\rho$-manifold with $\rho_\blt: \Aff(\RR^d)\hookrightarrow\Diff^s(\RR^d)$ the standard affine action.

If $\eua\leq \Aff(\RR^d)$, a \textbf{$C^s$ $\eua$-affine manifold} is a $C^s$ $\rho$-manifold, where $\rho$ is the restriction of the standard affine action to $\eua$; we'll denote the associated $\rho$-diffeomorphism group by $\Diff^s_\eua(M)$. A \textbf{$C^s$ globally-$\eua$-affine manifold} is likewise a $C^s$ globally-$\rho$-manifold where $\rho$ is the restriction of the standard affine action to $\eua$. Let us denote by $\DGL(\RR^d)$ the group of diagonal invertible linear automorphisms\fn{Note that this notation implies that there is a chosen but suppressed splitting of $\RR^d$ into $1$ dimensional subspaces.} of $\RR^d$ and $\DAff(\RR^d)\leq \Aff(\RR^d)$ the subgroup of $\Aff(\RR^d)$ consisting of all affine automorphisms of $\RR^d$ with diagonal linear part; so $\DAff(\RR^d)=\RR^d\rtimes \DGL(\RR^d)$.

\done
\end{dfn}

\begin{rem}\label{037}
In general $\Diff^s_\rho(Y)$ is not a closed subgroup of $\Diff^s(Y)$, nor are the groups $\Diff^s_\rho(Y)$ and $\im(\rho)\leq \Diff^s_\rho(X)$ isomorphic. If $Y$ is a globally-$\rho$-manifold then any chart provides an isomorphism $\Diff^s_\rho(Y)\cong\im(\rho)$. If $\eua\leq\Aff(\RR^d)$, then $\Diff^s_\eua(Y)$ is a closed subgroup of $\Diff^s(Y)$, the latter being endowed with the uniform $C^s$ topology. Thus for any globally-$\eua$-affine manifold $Y$, $\Diff^s_\eua(Y)\cong \eua$.

\done
\end{rem}

\subsection{Suspensions and Time Changes\label{052}}

In this section we briefly discuss the two categorical constructions that we will recover in the context of \autoref{001}, namely, suspensions and time changes. Let $G$ be a unimodular Lie group, $\Gamma\leq G$ be a lattice, $s\in\ZZ_{\geq1}\times[0,1]$,  $N$ be a $C^s$ manifold and $(\nu,\beta)\in\mfE^s(\Gamma\lact N)$. We do not disregard the case of discrete $G$; in this case $\Gamma$ is a finite index subgroup. The aim of the suspension construction is to produce a system, in an \ae-unique manner in the first coordinate, in $\mfE^s(G\lact M)$ for some $C^s$ manifold $M$. First define the diagonal $C^s$ action

$$\dbar^{\beta}_\blt: \Gamma\to \Diff^s(G\times N),\,\, t\mapsto [(g,x)\mapsto (gt^{-1}, \beta_t(x))].$$

$\dbar^\beta$ is a free and proper action, so that the orbit space $\faktor{G\times N}{\dbar^\beta}= G\otimes^\beta N$ has a unique $C^s$ manifold structure w/r/t which the canonical map $G\times N\to G\otimes^\beta N$ is a $C^s$ submersion. There is a natural $C^s$ submersion map $\pi^\beta: G\otimes^\beta N\to \faktor{G}{\Gamma}$ which gives a $C^s$ fiber bundle structure with each fiber $C^s$ diffeomorphic to $N$. Moreover the left translation action $G\to \Diff^s(G\times N), t\mapsto[(g,x)\mapsto (tg,x)]$ commutes with $\dbar^\beta$, hence induces an action $\hbar^\beta_\blt:G\to\Diff^s(G\otimes^\beta N)$ called the \textbf{suspension} of $\beta$.

Let us denote by $\haar_G$ the Haar measure on $G$ that $\pi^\beta$ pushes forward to the unique $G$ invariant Borel probability measure $\haar_{\faktor{G}{\Gamma}}$ on $\faktor{G}{\Gamma}$. Then pushing the product measure $\haar_G\otimes \nu$ on $G\times N$ forward to $G\otimes^\beta N$ gives a Borel probability measure $\haar_{\faktor{G}{\Gamma}}\otimes^\beta \nu$. It's straightforward that $\hbar^\beta$ preserves $\haar_{\faktor{G}{\Gamma}}\otimes^\beta \nu$, so that we have produced the \textbf{suspension system}

$$\left(\haar_{\faktor{G}{\Gamma}}\otimes^\beta \nu,\hbar^\beta\right)\in\mfE^s(G\lact G\otimes^\beta N).$$

Note that $\faktor{G}{\Gamma}$ also carries a left translation action $l^G$ of $G$ under which $\haar_{\faktor{G}{\Gamma}}$ is invariant, so that $\left(\haar_{\faktor{G}{\Gamma}},l^G\right)\in\mfE^{\aff}\left(G\lact \faktor{G}{\Gamma}\right)$ which is a factor system of $\left(\haar_{\faktor{G}{\Gamma}}\otimes^\beta \nu,\hbar^\beta\right)$ via $\pi^\beta$. For each $g\in G$ such that $g\Gamma=\Gamma g$, the fiber $\{g\Gamma\}\otimes^\beta N=\bck{\pi^\beta}(g\Gamma)$ is preserved by the suspension $\hbar^\beta$, and we have a $C^s$ isomorphism $\iota_{g\Gamma}^\beta: (\nu,\beta)\xrightarrow{\phantom{\cong}\cong_{C^s}\phantom{\cong}}\left(\restr{\haar_{\faktor{G}{\Gamma}}\otimes^\beta \nu}{\{g\Gamma\}\otimes^\beta N},\restr{\hbar^\beta}{\{g\Gamma\}\otimes^\beta N}\right)$ that is a $C^s$ embedding $N\hookrightarrow G\otimes^\beta N$.

\begin{obs}\label{019}
In the special case of $G=\RR^k$ and $\Gamma=\ZZ^k$, since the normalizer of $\ZZ^k$ is the whole $\RR^k$, the suspension system restricted to any fiber of the suspension manifold is a copy of the original $\ZZ^k$ system. Moreover the factor system $\left(\haar_{\TT^k},l^{\RR^k}\right)\in\mfE^{\aff}\left(\RR^k\lact \TT^k\right)$ can be replaced by the translation system $\left(\haar_{\TT^k},l^{\TT^k}\right)\in\mfE^{\aff}\left(\TT^k\lact \TT^k\right).$ Thus we may think of the the totality of the discussion above in  this section as a bundle of systems:

\begin{center}

\begin{tikzcd}
{(\nu,\beta)} \arrow[rr, "\iota^\beta", hook] &  & {\left(\haar_{\TT^k}\otimes^\beta \nu,\hbar^\beta\right)} \arrow[rr, "\pi^\beta", two heads] &  & {\left(\haar_{\TT^k},l^{\TT^k}\right)}
\end{tikzcd}

\end{center}

In \autoref{001} we recover such a bundle up to measure theoretical isomorphism, with the fiber action affine Cartan.

\done
\end{obs}

\begin{obs}\label{020}
Let us consider the case of $G=\ZZ^k$ and $\Gamma=\fwd{j}(\ZZ^k)$ the image of an embedding. $\fwd{j}(\ZZ^k)\leq \ZZ^k$ being a lattice means that the factor group $F=\faktor{\ZZ^k}{\fwd{j}(\ZZ^k)}$ is finite. Similar to \autoref{019} we have a bundle of systems:

\begin{center}
\begin{tikzcd}
{(\nu,\beta)} \arrow[rr, "\iota^\beta", hook] &  & {\left(\haar_{F}\otimes^\beta \nu,\hbar^\beta\right)} \arrow[rr, "\pi^\beta", two heads] &  & {\left(\haar_{F},l^{F}\right)}
\end{tikzcd}
\end{center}

where $\haar_F$ is the counting measure on $F$ normalized by $\dfrac{1}{\#(F)}$.

\done
\end{obs}

Let us now discuss time changes. We generalize the $C^s$ version of the theory of time changes in \cite{MR201606,MR855312} and \cite[p.74]{MR614142}\fn{Also see \cite{MR0224771} for a slightly different, infinitesimal approach.} to systems with acting group a unimodular Lie group. The aim of the time change construction is to reparameterize a system in an orbit-preserving manner. Let us endow $G$ with a suppressed left translation invariant Riemannian metric and let $(\mu,\alpha)\in\mfE^s(G\lact M)$ be a system. Let $\lambda$ be a $C^s$ \textbf{($G$-valued $1$-)cocycle} over $\alpha$, that is, $\lambda\in C^s(G\times M;G)$ is such that

$$\all t_1,t_2\in G,\all x\in M: \lambda(t_2,x) \lambda(t_1t_2,x)^{-1} \lambda(t_1,\alpha_{t_2}(x))=e_G.$$

Suppose further that for each fixed $x\in M$, the map $t\mapsto \lambda(t,x)$ is an orientation preserving diffeomorphism; that is, for any $x\in M$, $\lambda(\blt,x)\in\Diff^s_+(G,e_G)$. Then we have that $(\lambda,\id_M):G\times M\to G\times M, (t,x)\mapsto (\lambda(t,x),x)$ is a $C^s$ diffeomorphism with inverse of the form $(\kappa,\id_M)$ for some $\kappa\in C^s(G\times M;G)$ (explicitly, $\kappa: (t,x)\mapsto \lambda(\blt,x)^{-1}(t)$). Define a family of $C^s$ self-maps of $M$ by $\beta_\blt=\alpha_{\kappa\blt}: G\to C^s(M,M), t\mapsto [x\mapsto \alpha_{\kappa(t,x)}(x)]$. As $\alpha=\beta_\lambda$ also, it's straightforward to verify that $\kappa$ is a $C^s$ cocycle over $\beta$, and $\beta=\alpha_\kappa$ is a $G$ action by $C^s$ diffeomorphisms; it's called the $\kappa$ \textbf{time change} of the action $\alpha$. Note that we have a commutative diagram (in the category $\cat{Man}^s$ of $C^s$ manifolds):

\begin{center}
\begin{tikzcd}
G\times M \arrow[r, "\alpha=\beta_\lambda"] \arrow[d, "{(\lambda,\id_M)=(\kappa,\id_M)^{-1}}"'] & M \arrow[d, "\id_M"] \\
G\times M \arrow[r, "\beta=\alpha_\kappa"']                                                     & M                   
\end{tikzcd}

\end{center}

We also need to modify the $\alpha$-invariant measure $\mu$ to get a $\beta$-invariant measure. Define a new Borel probability measure $\nu$ on $M$ by 

$$\nu=\mu_\kappa:\bor(M)\to [0,1],\,\, B\mapsto \dfrac{\int_B \det\left(T^1\lambda(e_G,x)\right)\, d\mu(x)}{\int_M \det\left(T^1\lambda(e_G,x)\right)\, d\mu(x)},$$

where $T^1\lambda(e_G,x)$ is, for fixed $x\in M$, the derivative of $G\to G, t\mapsto \lambda(t,x)$ evaluated at $e_G$, and the determinant is w/r/t the Riemannian metric fixed on $G$\fn{Instead of the Riemannian Jacobian one could alternatively use Radon-Nikodym derivative of the pushforward of a Haar measure $\haar_G$ on $G$ by $\lambda(\blt,x)$ w/r/t $\haar_G$.}. It's straightforward to verify that $\nu$ is $\beta$-invariant, e.g. by way of taking the Laplace transform or derivative of $G\to \RR, t\mapsto \EE_{\nu}(f\circ \beta_t)$ for some anonymous $f\in C^0(M;\RR)$. It's also straightforward that $(\mu,\alpha)$ is ergodic iff $(\nu,\beta)$ is ergodic, provided that $G$ is a group for which a pointwise ergodic theorem holds that allows one to replace $\alpha$-time averages of the time change $\lambda$ with its $\mu$-space averages\fn{E.g. for $G$ compactly generated abelian (with possibly some further reasonable conditions on $G$ and $\lambda$); see \cite{MR1696824,MR1865397,MR384997} and \cite[p.195, Ch.5]{MR961261}. \label{104}}. Finally the entropy gauges of $(\mu,\alpha)$ and $(\nu,\beta)$ are proportional to each other by an Abramov formula:

\begin{lem}\label{105}

Let $G$ be a unimodular Lie group. Let $(\mu,\alpha)\in\mfE^s(G\lact M)$ and $\kappa$ be a $C^s$ time change of $(\mu,\alpha)$. Then we have

$$\all t\in G: \mfe_{(\mu,\alpha)}(t) =  \left(\int_M \det\left(T^1\lambda(e_G,x)\right)\, d\mu(x)\right)\,\mfe_{(\mu_\kappa,\alpha_\kappa)}(t).$$

\done
\end{lem}

Observe that we may consider any orientation preserving Lie group automorphism $\lambda\in\Aut_{\tcat{Lie},+}(G)$ as a cocyle over $\alpha$. In this case $\kappa=\lambda^{-1}\in\Aut_{\tcat{Lie},+}(G)$, $\nu=\mu_\kappa=\mu=\nu_\lambda$ and

$$\all t\in G: \mfe_{(\mu,\alpha)}(t)=\det(\lie(\lambda))\,\mfe_{(\mu_\kappa,\alpha_\kappa)}(t),$$

where $\lie(\lambda)=T_{e_G}\lambda:\lie(G)\to\lie(G)$ is the induced Lie algebra automorphism. Since such group automorphisms $\lambda$ are in bijection with cocycles over $\alpha$ that are constant in $x$, we'll call such cocycles and in particular time changes \textbf{constant in space}.

\begin{rem}\label{136}

Let us note in closing that instead of assuming $\lambda\in C^s(G\times M;G)$ we could have assumed that $\lambda\in L^0(G\times M; G)$ with $\all x\in_\mu M: \lambda(\blt,x)\in \Diff^s_+(G)$ and

$$\all t_1,t_2\in G,\all x\in_\mu M: \lambda(t_2,x) \lambda(t_1t_2,x)^{-1} \lambda(t_1,\alpha_{t_2}(x))=e_G.$$

Then the theory of time changes as discussed above still holds (with the commutative diagram now being in $\cat{sMeas}$), and in particular we still have \autoref{105}. Let us call such cocycles \textbf{$C^s$ \ae-time changes}.

\done
\end{rem}

\subsection{Oseledets Multiplicative Ergodic Theorem and Lyapunov Geometry}

Lyapunov exponents we'll be using are provided by an Oseledets Theorem for the derivative cocycle of a $C^1$ action of $\RR^k$ on a compact $C^\infty$ manifold. Below we state such an Oseledets Theorem.

\begin{thm}[Oseledets\fn{\cite{MR0240280}. Also see \cite[p.305, Prop.2.1]{MR2811602}; \cite[p.26]{MR2798364}; \cite[pp.19-20, Thm.1-2]{MR1101082} \cite[p.5, Thm.2.4]{brown2016smooth}; \cite[p.74, Thm.A]{MR1213080} (= \cite[p.13, Thm.A]{MR2690742}).}]\label{031}

Let $(\mu,\alpha)\in\mfE^1(\RR^k\lact M)$. If $(\mu,\alpha)$ is locally free and ergodic, then $k\leq \dim(M)$. If $k=\dim(M)$, then the only Lyapunov exponent of $(\mu,\alpha)$ is $0$; otherwise there is a subset $\Osel\in\bor(M,\alpha)$ with $\Osel=_\mu M$ such that

\begin{description}

\item[EXI1] there is a unique number $ l\in\OL{\dim(M)-k}$; and $\all i\in\OL{l}, \exi \delta^i\in \OL{\dim(M)}$:

$$\sum_{i\in\OL{l}}\delta^i=\dim(M)-k,$$

\item[EXI2] 
there is a unique up to reordering linear operator 

\begin{align*}
 X  
&=(X^1,X^2,...,X^{\dim(M)-k}, \underbrace{0,0,...,0}_{k \text{ many}})\\
&=(\underbrace{\chi^1,...,\chi^1}_{\delta^1 \text{ many}},\underbrace{\chi^2,...,\chi^2}_{\delta^2 \text{ many}},..., \underbrace{\chi^{l},...,\chi^{l}}_{\delta^{l} \text{ many}}, \underbrace{0,0,...,0}_{k \text{ many}})
\in\Hom(\RR^k,\RR^{\dim(M)})
\end{align*}

with $i\neq j\ifr \chi^i\neq \chi^j$ (but not necessarily $\chi^i\neq0$),

\item[EXI3] for any $x\in \Osel$  there is a unique splitting $ T_xM=O_x\oplus\left(\bigoplus_{i\in \OL{l}}L_x^i\right)$ with 

$$\all i\in\OL{l}: \dim L_x^i=\#\{j\in\OL{\dim(M)}\,|\, X^j=\chi^i\}=\delta^i$$

\end{description}

such that

\begin{description}
\item[ASYM1]  For any $C^0$ fiberwise norm on $M$,

\begin{align*}
\all x\in \Osel,\all i\in\OL{l}:&
\lim_{|t|\to\infty} \sup_{v\in L^i_x\setminus 0} \dfrac{\log |T_x\alpha_t v|- \chi^i(t)}{|t|}=0,
\end{align*}

\item[ASYM2]  For any $C^0$ density on $M$,

$$\all x\in \Osel:\lim_{|t|\to\infty}\dfrac{\log \Jac_x(\alpha_t) - \sum_{i\in\OL{l}}\delta^i\, \chi^i (t)}{|t|}=0.$$

\item[ASYM3] For any $C^0$ Riemannian metric on $M$,

\begin{align*}
\all x\in \Osel,\all I\subseteq\OL{l}:& \lim_{|t|\to\infty} \dfrac{\log \left\Vert \proj\left(\fwd{T_x\alpha_t}\left(\bigoplus_{i\in I} L_x^j\right)\right)- \proj\left(O_{\alpha_t(x)}\right)\right\Vert}{|t|}=0\,\,\text{ and }\\
& \lim_{|t|\to\infty} \dfrac{\log \left\Vert \proj\left(\fwd{T_x\alpha_t}\left(\bigoplus_{i\in I} L_x^j\right)\right)- \proj\left(\fwd{T_x\alpha_t}\left(\bigoplus_{j\in \OL{l}\setminus I} L_x^j\right)\right)\right\Vert}{|t|}=0.
\end{align*}

\end{description}

Further,

\begin{description}

\item[INV4]  $O_\blt\oplus\left(\bigoplus_{i\in \OL{l}}L_\blt^i\right)\in L^0(\Spl(TM)\to M)$ is $\Ad^\alpha$-invariant.

\end{description}

\done
\end{thm}

Here $\Spl(TM)\to M$ is the $C^\infty$ bundle of splittings of $TM$, and $\Ad^\alpha$ is the adjoint action on the measurable sections of the bundle of splittings: $\Ad^\alpha_t: E\mapsto \left[x\mapsto  \fwd{T_{\alpha_t^{-1}(x)}\alpha_t}(E_{\alpha_t^{-1}(x)})\right]$. Let us now list some of the standard definitions of Lyapunov geometry.

\begin{dfn}\label{011}

In the context of \autoref{031},

\begin{itemize}

\item $\Osel$ is the \textbf{set of Lyapunov-Perron regular points}, 

\item $X$ is the  \textbf{Lyapunov operator},

\item $\wt{X}=(\chi^1,\chi^2,...,\chi^{l})\in \Hom(\RR^k,\RR^{l})$ is the {\bf reduced Lyapunov operator},

\item $\dim(\ker(X))\in\UL{k+1}$ is the \textbf{defect},

\item $\dim(\im(X))\in \UL{\dim(M)-k+1}$ is the \textbf{codefect},

\item The set $\LSpec=\{\chi^1,\chi^2,...,\chi^{l}\}$ of linear functionals is the {\bf (nonorbital) Lyapunov spectrum},

\item Any element  $\chi\in \LSpec$ is a {\bf (nonorbital) Lyapunov exponent},

\item $\LSpec^\ast=\LSpec\setminus0$ is the \textbf{reduced Lyapunov spectrum},

\item For any $\chi\in\LSpec$ and any $x\in\Osel$,  $L^\chi_x$ is the \textbf{Lyapunov subspace} associated to $\chi$ at $x$ and $L^\chi$ is the \textbf{Lyapunov \ae-subbundle} associated to $\chi$,

\item $TM=_{\aem} \bigoplus_{i\in\OL{l}}L^i$ is the {\bf (fine) Lyapunov \ae-splitting},

\item $\delta^i$ is the  \textbf{dynamical multiplicity}  of $\chi^i$. 

\item For $t\in\RR^k$, the function $\LSpec\to \{-,0,+\}, \chi\mapsto \text{sign}(\chi(t))$ is called the \textbf{Lyapunov signature} of the time $t$,

\item For any $t\in\RR^k$ and $x\in \Osel$: the \textbf{stable, unstable, center, stable-center} and \textbf{center-unstable subspaces} of $\alpha_t$ at $x$ are defined as follows, respectively:

\begin{align*}
&S_x(\alpha_t) = \bigoplus_{\substack{\chi\in\LSpec\\\chi(t)<0}} L^\chi_x,\quad
U_x(\alpha_t) = \bigoplus_{\substack{\chi\in\LSpec\\\chi(t)>0}} L^\chi_x,\quad
C_x(\alpha_t) = \bigoplus_{\substack{\chi\in\LSpec\\\chi(t)=0}} L^\chi_x\\
&SC_x(\alpha_t) = S_x(\alpha_t)\oplus C_x(\alpha_t),\quad
CU_x(\alpha_t) = C_x(\alpha_t)\oplus U_x(\alpha_t);\\
\end{align*}

they assemble to \textbf{stable, unstable, center, stable-center} and \textbf{center-unstable \ae-subbundles} of $\alpha_t$, respectively.

\end{itemize}

\done
\end{dfn}

\begin{rem}\label{041}
In the context of \autoref{031}, we'll use the Lyapunov exponents as superscripts if we consider a Lyapunov exponent without referring to an enumeration, e.g. we may write $\bigoplus_{\chi\in \LSpec} L^\chi$ for the (fine) Lyapunov \ae-splitting. We'll refer to this notation as the \textbf{intrinsic notation} for Lyapunov geometry.
\done
\end{rem}

\begin{rem}\label{098}

Not fixing an ergodic measure at first gives a basepoint-dependent version of \autoref{031}, which provides a measurable subset $\Osel$ and $l_\blt,\delta^i_\blt, X_\blt$ become functions on $\Osel$ such that

\begin{description}

\item[INV1] $\Osel\in\bor(M;\alpha)$ and $\all \mu\in\Prob(M;\alpha): \Osel=_\mu M$,

\item[INV2] $l_\blt\in L^0(\Osel;\OL{d})$ and is $\alpha$-invariant,

\item[INV3] $X_\blt=(X^1_\blt,X^2_\blt,...,X^d_\blt) \in L^0\left(\Osel; \Hom(\RR^k,\RR^d)\right)$ and  is $\alpha$-invariant.

\end{description}

In this case we consequently have basepoint dependent Lyapunov spectrum also. The basepoint dependent version of Oseledets Theorem is useful if no measure is fixed a priori, or if the fixed measure is not ergodic w/r/t the action, or if one is considering the action of a subgroup (e.g. the iterates of a time-$t$ map). Also note that the set $\Osel$ in \autoref{031} can be considered to be attached to the action $\alpha$ instead of the system $(\mu,\alpha)$.

\done
\end{rem}

\begin{cor}\label{042}
In the context of \autoref{031}, for any $C^0$ fiberwise norm  on $X$ and for any $\eps\in\RRP$, there is a  measurable function $C_\eps:\Osel\to\RRP$ such that for any $x\in \Osel$, $i\in\OL{l}$, $v\in L^i_x$ and $t\in\RR^k$

\begin{itemize}

\item $\dfrac{1}{C_\eps(x)}\, e^{\chi^i(t)-\frac{1}{2}\eps|t|} |v| \leq |T_x\alpha_tv|\leq C_\eps(x)\, e^{\chi^i(t)+\frac{1}{2}\eps|t|}|v|$, and

\item $C_{\eps}(\alpha_t(x))\leq C_\eps(x)\,e^{\eps|t|}$.

\end{itemize}

Similarly for any $C^0$ Riemannian metric on $M$  and for any $\eps\in\RRP$, there are measurable functions $C_{\eps},K_{\eps}:\Osel\to\RRP$ such that $C_{\eps}$ is as above and for any $x\in \Osel$, $I\subseteq \OL{l}$, $v\in L^\chi_x$ and $t\in\RR^k$

\begin{itemize}

\item  $\left\Vert \proj\left( \bigoplus_{i\in I} L_x^i \right) -  \proj\left( O_x \right)  \right\Vert\geq K_\eps(x)$,

\item $\left\Vert \proj\left( \bigoplus_{i\in I} L_x^i \right) -  \proj\left( \bigoplus_{j\in \OL{l}\setminus I} L_x^j \right)  \right\Vert\geq K_\eps(x)$, and

\item $K_{\eps}(\alpha_t(x))\geq K_\eps(x)\, e^{-\eps|t|}$.

\end{itemize}

We will refer to $C_\blt,K_\blt:\RRP\to L^0(\Osel;\RRP)$ as  \textbf{comparison families} of $(\mu,\alpha)$ w/r/t the chosen fiberwise norm or Riemannian metric.

\done
\end{cor}

\begin{rem}\label{043}
\textbf{ASYM3} of \autoref{031} is traditionally formulated in terms of the angle between the subspaces $\fwd{T_x\alpha_t}\left(\bigoplus_{i\in I} L_x^i\right)$ and  $\fwd{T_x\alpha_t}\left(\bigoplus_{j\in \OL{l}\setminus I} L_x^j\right)$ of $T_{\alpha_t(x)}M$. Recall that the \textbf{angle} $\sphericalangle(W_1,W_2)$ between two nonzero subspaces $W_1,W_2\leq V$ of an inner product space $V$ is defined as

$$\sphericalangle(W_1,W_2)=\inf\left\{\left.\sphericalangle(w_1,w_2)=\arccos\left(\dfrac{\langle w_1,w_2\rangle}{|w_1||w_2|}\right)\,\right|\, w_1\in W_1\setminus0,w_2\in W_2\setminus0\right\}\in[0,\pi/2].$$

Then \textbf{ASYM3} is equivalent to

$$\all I\subseteq\OL{l_x}: \lim_{|t|\to\infty} \dfrac{\log \sin\sphericalangle\left( \fwd{T_x\alpha_t}\left(\bigoplus_{i\in I} L_x^i\right) ,\, \fwd{T_x\alpha_t}\left(\bigoplus_{j\in \OL{l}\setminus I} L_x^j\right) \right) }{|t|}=0.$$

Accordingly, the first property of $K_\eps$ in \autoref{042} is equivalent to

$$\sin\sphericalangle\left( \bigoplus_{i\in I} L_x^i\,\,,\,\,\bigoplus_{j\in \OL{l}\setminus I} L_x^j \right) \geq K_\eps(x).$$

\done
\end{rem}

\begin{dfn}\label{013}

In the context of \autoref{031},

\begin{itemize}

\item For $\chi\in \LSpec^\ast$,  $K^\chi=\ker(\chi)\leq \RR^k$ is called the \textbf{Lyapunov hyperplane}  of $\chi$.

\item Any connected component of $\RR^k\setminus \bigcup_{\chi\in\LSpec^\ast}K^\chi$ is called a \textbf{Weyl chamber}; let's denote the set of all Weyl chambers  by  $\Cham$  $=\pi_0\left(\RR^k\setminus \bigcup_{\chi\in\LSpec^\ast}K^\chi\right)$. Let us denote by $-$ the natural involution on $\Cham$.

\item Any element of $\Wall^0=\RR^k\setminus \bigcup_{\chi\in\LSpec^\ast}K^\chi$ is a \textbf{chamber time},

\item Any element of $\bigcup_{\chi\in\LSpec^\ast}K^\chi=\RR^k\setminus \Wall^0$ is a \textbf{wall time}\fn{Compared to standard terminology we prefer more descriptive nomenclature; in \cite{MR1858547} a chamber time is called a generic element, a wall time is called a singular element, and a wall time in $\Wall^1$ is called a generic singular element.}.

\end{itemize}

If $I\subseteq \LSpec^\ast$ is a subset, put 

$$\Wall^I=\left(\bigcap_{\chi\in I}K^\chi\right) \setminus \left(\bigcup_{\rho\in \LSpec^\ast\setminus I} K^\rho\right).$$ 

For brevity we also put $\Wall^\chi=\Wall^{\{\chi\}}$. Finally put for $i\in\OL{l}$:

$$\Wall^i=\bigcup_{\substack{I\subseteq \LSpec^\ast\\\#(I)=i}}\Wall^I.$$

We call an element of $\Wall^I$ an \textbf{$I$-wall time}, an element of $\Wall^{\chi}$ a \textbf{$\chi$-wall time}, and an element of $\Wall^i$ a \textbf{wall time of order $i$}.

Let $\chi\in \LSpec^\ast$ be a nonzero Lyapunov exponent and $\euc\in \Cham$ be a Weyl chamber. $\restr{\chi}{\euc}:\euc\to \RR$ is not constant, but it's either always positive or always negative. This way we have a \textbf{Lyapunov pairing} $(\blt,\blt):\LSpec^\ast \times \Cham \to \{\pm\}$. The function $(\blt,\euc):\LSpec^\ast\to \{\pm\}$ is called the \textbf{Lyapunov signature} of the chamber $\euc$. For any $\euc\in\Cham$ and $x\in\Osel$  the \textbf{stable} and  \textbf{unstable} subspaces of $\euc$ at $x$ are defined as follows, respectively:

\begin{align*}
S^\euc_x &= \bigoplus_{\substack{\chi\in\LSpec\\(\chi,\euc)=-}} L^\chi_x,\quad
U^\euc_x = \bigoplus_{\substack{\chi\in\LSpec\\(\chi,\euc)=+}} L^\chi_x.\quad
\end{align*}

\done
\end{dfn}

\begin{rem}\label{106}

Note that we have $(\chi,-\euc)=-(\chi,\euc)$. Also note that the stable and unstable subspaces of a Weyl chamber coincides with the stable and unstable subspaces of any time-$t$ map of the action if $t$ is chosen from the Weyl chamber in question:

$$\all \euc\in\Cham,\all t\in \euc: S^\euc=S(\alpha_t), U^\euc=U(\alpha_t).$$

Another remark is that we have a partition

$$\RR^k=\biguplus_{i\in \UL{\dim(M)-k+1}}\Wall^i = \Wall^0\uplus \Wall^1 \uplus\cdots\uplus \Wall^{\dim(M)-k}.$$

Finally note that in general $\Wall^i$ is not the same as $\Wall^{\chi^i}$.

\done
\end{rem}

\begin{obs}\label{029}
Let us also note the relationship between the Lyapunov exponents of $\alpha $ and the Lyapunov exponents of any of its time-$t$ map $\alpha_t$. \textbf{ASYM1} of \autoref{031} implies that

\begin{description}
\item[ASYM1$_\ast$]  For any $C^0$ fiberwise norm on $M$ and for any $t\in\RR^k$:

\begin{align*}
\all x\in \Osel,\all i\in\OL{l}:&
\lim_{|n|\to\infty} \sup_{v\in L^i_x\setminus 0} \dfrac{\log |T_x\alpha_t^n v|}{n}= \chi^i(t),
\end{align*}

\end{description}

In particular, $X(t)=(X^1(t ),X^2(t ),X^{d-k}(t),0,0,...,0)\in\RR^d$ has as its components the Lyapunov exponents with multiplicity of the diffeomorphism $\alpha_{t }$, and even though the system $(\mu,\alpha_t)$ may fail to be ergodic its Lyapunov exponents are \ae-constant. Note that the components need not be ordered, as is common in rank-$1$ dynamics; thus we have a function $\sigma_{\blt}: \RR^k\to \Bij(\OL{d-k})$ which is unique up to permutations of same entries such that

$$\all t \in\RR^k: X^{\sigma_{t}(1)}(t ) \geq X^{\sigma_{t }(2)}(t ) \geq \cdots \geq X^{\sigma_{t}(d-k)}(t ).$$

\done
\end{obs}

\subsection[\ae-Foliations, Partitions, and Conditional Measures]{\ae-Foliations, Partitions, and Conditional Measures}\label{079}

In this subsection we briefly review conditional measures of a probability measure on a manifold along the leaves of an \ae-foliation. The classical reference for this section is \cite{MR0047744}; for the adaptation of the abstract theory to the geometric setting we loosely follow \cite{MR415679} and \cite{MR2648695}. Let $X$ be a $C^\infty$ manifold of dimension $d\in\ZZP$, and let $s\in\ZZP\times[0,1]$.

\begin{dfn}\label{080}
A \textbf{ partial $L^0$ foliation with $C^s$ leaves}  of $X$ is a pair $(\supp(\mcf),\mcf)$, where $\supp(\mcf)\subseteq X$ is a measurable subset, called the \textbf{support of $\mcf$}, into connected injectively immersed $C^s$ submanifolds $\{\mcf_x\,|\, x\in \supp(\mcf)\}$ possibly of varying dimensions such that for any $x\in \supp(\mcf)$, there is a $\delta_x\in\RRP$, a measurable subset $V\subseteq \RR^{d-\dim(\mcf_x)}[0|\leq\delta_x]$ containing $0$, a measurable subset $W\subseteq \supp(\mcf)$ containing $x$, and an $\cat{Mble}$-isomorphism 

$$\phi: \left(\RR^{\dim(\mcf_x)}[0|\leq\delta_x]\times V,(0,0)\right)\to (W,x)$$

such that for any $b\in V$ \fn{For $(X,d)$ a metric space, $x\in X$ and $r\in \RRR$, $X[x|\leq r], X[x|<r], X[x|=r]$ is the closed ball, open ball and sphere, respectively, centered at $x$ of radius $r$. \label{103}}:

\begin{align*}
\phi(\blt,b)&\in \Emb^s\left(\RR^{\dim(\mcf_x)}[0|\leq\delta_x];X\right)\text{ and}\\
&\fwd{\phi}\left(\RR^{\dim(\mcf_x)}[0|\leq\delta_x]\times\{b\}\right)=W\cap \mcf_{\phi(0,b)}.
\end{align*}

If $\mu$ is a Borel probability measure on $X$ and $(\supp(\mcf),\mcf)$ is a partial $L^0$ foliation with $C^s$ leaves such that $\supp(\mcf)=_\mu X$, then $(\supp(\mcf),\mcf)$ is called an\textbf{ \ae-foliation with $C^s$ leaves}. We'll systematically suppress the support of a partial $L^0$ foliation. Denote by $\aeFol^s(X,\mu)$ the collection of all \ae-foliations of $X$ with $C^s$ leaves.

\done
\end{dfn}

\begin{rem}\label{081}
In the context of \autoref{080}, note that the "measurable partial flowboxes" $\phi$ are superfluous, and as such an \ae-foliation may fail to be a foliation of its support, even when this support inherits a manifold structure from $X$. Still we prefer to write them as \ae-foliations are indeed adaptations of foliations to the theory of nonuniform hyperbolicity.

We also remark that the \ae-foliations coming from dynamics, i.e. the Lyapunov, stable, and unstable \ae-foliations (see \autoref{051} and \autoref{064} below) have the extra property of admitting H\"{o}lder partial flowboxes on Lusin-Pesin sets (see \autoref{061} below). This extra structure will not be used in this subsection.

\done
\end{rem}

\begin{dfn}\label{120}

Let $\mu$ be a Borel probability measure on $X$ and let $\mcf$ be an \ae-foliation of $(X,\mu)$ with $C^s$ leaves. Let $L,R$ be two injectively immersed $C^s$ submanifolds of $X$ transverse to $\mcf$. Then a \textbf{local holonomy} from $L$ to $R$ along $\mcf$ is a measurable function $\mcf_{L\leftarrow R}: L'\to R$ defined on a measurable subset $L'\subseteq L$ with the property that 

$$\all l\in L': \mcf_{L\leftarrow R}(l)\in \mcf_l \cap R.$$

A \textbf{holonomy} is a local holonomy with $L'=L$. If $\mcg$ is another \ae-foliation with $C^s$ leaves that is transverse to $\mcf$, then for any $x,y\in_\mu X$ with $\mcg_x\cap \mcg_y=\emptyset$ let us make the abbreviation $\mcf_{y\leftarrow x}=\mcf_{\mcg_y\leftarrow \mcg_x}$.

\done
\end{dfn}

Our next aim is to define, given an anonymous  Borel probability measure $\mu$ on $X$ and an \ae-foliation $\mcf$ of $X$, conditional measures $\mu^\mcf_\blt$. One can't directly define this disintegration in situations where it would be important (including the context of this paper), however by systematically coarsening an \ae-foliation this can be achieved under certain circumstances (including the context of this paper). To describe this procedure we switch to a more abstract setting momentarily.

\begin{dfn}\label{089}

A \textbf{standard measurable space} $Y$ is a measurable space $\cat{Mble}$-isomorphic to a Borel measurable subset of a compact metric space. A \textbf{standard probability space} $(Y,\eta)$ is a standard measurable space $Y$ endowed with a probability measure $\eta$. An \textbf{\ae-partition} $\mcp$ of a standard probability space $(Y,\eta)$ is a collection of measurable subsets that are \ae-disjoint and that \ae-cover the whole space, that is, $\all P_1,P_2\in\mcp: P_1\cap P_2=_\eta \emptyset$, and $\bigcup_{P\in\mcp} P=_\eta Y$; syntactically we allow $\eta$-negligible subsets to be elements of partitions. Given an \ae-partition $\mcp$ of $(Y,\eta)$, denote by $\pi^\mcp: (Y,\eta)\to \left(\faktor{Y}{\mcp}, \fwd{\pi^\mcp}(\eta)\right)$ the associated quotient map in $\cat{Meas}$. A \textbf{measurable \ae-partition} is an \ae-partition $\mcp$ such that $\left(\faktor{Y}{\mcp}, \fwd{\pi^\mcp}(\eta)\right)$ is a standard probability space. Denote by $\aePar(Y,\eta)$ the filtered category of \ae-partitions with with $\mcp\to\mcq$ iff $\mcp$ \ae-refines $\mcq$ and by $\aemPar(Y,\eta)$ the full subcategory of measurable \ae-partitions. We write $\mcp \sqcap \mcq$ for the coarsest common refiner of $\mcp$ and $\mcq$ and $\mcp \sqcup \mcq$ for the finest common refinee of $\mcp$ and $\mcq$\fn{In standard notation $\mcp\to\mcq$, $\mcp \sqcap \mcq$, $\mcp \sqcup \mcq$ are denoted by $\mcp\succeq \mcq$, $\mcp \vee \mcq$, $\mcp \wedge \mcq$, (or the other way around occasionally) respectively; we instead follow categorical conventions. Accordingly under the Galois correspondence between measurable \ae-partitions and sub-$\sigma$-algebras of the measure algebra the symbols for binary operations and relations flip. \label{095}}. $\Discrete_\eta$ and $\Indiscrete_\eta$ are the measurable \ae-partition into points and the one with exactly one element, respectively. Thus we have a diagram in $\aemPar(Y,\eta)$ with the lower the term the coarser the partition:

\begin{center}
\begin{tikzcd}
                & \Discrete_\eta \arrow[d]             &                 \\
                & \mcp\sqcap\mcq \arrow[ld] \arrow[rd] &                 \\
\mcp \arrow[rd] &                                      & \mcq \arrow[ld] \\
                & \mcp\sqcup\mcq \arrow[d]             &                 \\
                & \Indiscrete_\eta                     &                
\end{tikzcd}
\end{center}

\done
\end{dfn}

\begin{lem}\label{088}

Let $(Y,\eta)$ be a standard probability space. Then the following three categories are isomorphic:

\begin{enumerate}

\item The slice category of standard probability spaces under $(Y,\eta)$, where the arrows are factor maps.

\item The category $\aemPar(Y,\eta)$ of measurable \ae-partitions of $(Y,\eta)$.

\item The slice category $\Sub_{\tcat{Boo} }(\bor(Y,\eta))$ of sub-$\sigma$-algebras of the measure algebra of $(Y,\eta)$, where the arrows are inclusions.

\end{enumerate}

The isomorphism between the first two is by way of the fiber partition, and the isomorphism between the last two is by way of saturation. Under these isomorphisms the diagram in \autoref{089} transforms as follows:

\adjustbox{scale=0.4,center}{

\begin{tikzcd}
                                                                             & {(Y,\eta)} \arrow[d, two heads]                                                                                        &                                                                              &                 & \Discrete_\eta \arrow[d]             &                 &                                      & {\bor(Y,\eta)}                                                                &                                      \\
                                                                             & {\left(\faktor{Y}{\mcp\sqcap\mcq},\fwd{\pi^{\mcp\sqcap\mcq}}(\eta)\right)} \arrow[ld, two heads] \arrow[rd, two heads] &                                                                              &                 & \mcp\sqcap\mcq \arrow[ld] \arrow[rd] &                 &                                      & {\bor(Y,\eta;\mcp)\sqcup \bor(Y,\eta;\mcq)} \arrow[u, tail]                   &                                      \\
{\left(\faktor{Y}{\mcp},\fwd{\pi^{\mcp}}(\eta)\right)} \arrow[rd, two heads] &                                                                                                                        & {\left(\faktor{Y}{\mcq},\fwd{\pi^{\mcq}}(\eta)\right)} \arrow[ld, two heads] & \mcp \arrow[rd] &                                      & \mcq \arrow[ld] & {\bor(Y,\eta;\mcp)} \arrow[ru, tail] &                                                                               & {\bor(Y,\eta;\mcq)} \arrow[lu, tail] \\
                                                                             & {\left(\faktor{Y}{\mcp\sqcup\mcq},\fwd{\pi^{\mcp\sqcup\mcq}}(\eta)\right)} \arrow[d, two heads]                        &                                                                              &                 & \mcp\sqcup\mcq \arrow[d]             &                 &                                      & {\bor(Y,\eta;\mcp)\sqcap \bor(Y,\eta;\mcq)} \arrow[lu, tail] \arrow[ru, tail] &                                      \\
                                                                             & {(\pt,\delta_{\pt})}                                                                                                   &                                                                              &                 & \Indiscrete_\eta                     &                 &                                      & {\faktor{\{\emptyset,X\}}{\eta}} \arrow[u, tail]                              &                                     
\end{tikzcd}

}

\done
\end{lem}

\begin{dfn}\label{139}

Let $(Y,\eta)$  and $(Z,\zeta)$ be standard probability spaces and let $\pi: (Y,\eta)\to (Z,\zeta)$ be a measurable measure preserving map. Denote by $\vRadon(\pi): \vRadon(Y,\eta)\to (Z,\zeta)$ the \textbf{vertical Radon measure bundle} of $(Y,\eta)$ associated to $\pi$; by definition $\vRadon(\pi)_z= \Radon(\bck{\pi}(z))$ is the space of Radon measures with support a subset of the fiber $\bck{\pi}(z)$ at $z$ for $z\in_{\zeta} Z$. A section $\sigma_\blt$ of $\vRadon(\pi)$ is \textbf{measurable} if for any bounded measurable function $\phi\in L^0_b(Y; \RR)$, the function $Z\to \RR, z\mapsto \int_Y\phi(y)\, d \sigma_z(y)$ is measurable. A \textbf{disintegration} (or \textbf{system of conditional measures}) of $\eta$ w/r/t $\pi$ is a measurable section $\sigma_\blt$ of $\vRadon(\pi)$  such that

\begin{align*}
&\all \phi\in L^0_b(Y;\RR),\all\psi\in L^0_b(Z;\RR):\\
&\phantom{\all \phi\in L^0_b(Y;\RR),}\int_Y \psi\circ \pi(y) \, \phi(y)\, d\eta(y) = \int_Z \psi(z)\,\int_{\bck{\pi}(z)} \phi(y) \, d\sigma_z(y)\, d\zeta(z).
\end{align*}

\done
\end{dfn}

\begin{lem}\label{090}

Let $Y$ be a standard measurable space. Then for any probability measure $\eta$ on $Y$ and for any measurable \ae-partition $\mcp$ of $(Y,\eta)$, there is an \ae-unique disintegration $\eta^\mcp_\blt$ of $\eta$ w/r/t $\pi^\mcp$. Explicitly this means:

\begin{align*}
\all \phi\in L^0_b(Y;\RR): 
\int_Y \, \phi(y)\, d\eta(y) 
&= \int_{\faktor{Y}{\mcp}} \,\int_P  \phi(y) \, d\eta^\mcp_P(y)\, d\fwd{\pi^\mcp}(\eta)(P)\\
&= \int_{Y} \,\int_Y  \phi(y_2) \, d\eta^\mcp_{y_1}(y_2)\, d\eta(y_1),
\end{align*}

where in the last expression we made the abbreviation $\eta^\mcp_{y_1}=\eta^\mcp_{\mcp_{y_1}}=\eta^\mcp_{\pi^\mcp(y_1)}$.

\done
\end{lem}

\begin{dfn}\label{082}

Let $\mu$ be a Borel probability measure on $X$, $\mcf\in\aeFol^s(X,\mu)$ be an \ae-foliation and let $\mcp\in\aemPar(X,\mu)$ be a measurable \ae-partition. $\mcp$ is \textbf{subordinate} to $\mcf$ if $\mcp\to\mcf$ and for any $x\in_\mu X$, there is a neighborhood $N_x\subseteq \mcf_x$ of $x$ in the intrinsic topology of $\mcf_x$ such that $x\in N_x\subseteq \mcp_x\subseteq \mcf_x$.

\done
\end{dfn}

\begin{lem}\label{091}

Let $\mu$ be a Borel probability measure on $X$ and $\mcf\in\aeFol^s(X,\mu)$ be an \ae-foliation. Let $\mcp,\mcq\in\aemPar(X,\mu)$ be two measurable \ae-partitions. If both $\mcp$ and $\mcq$ are subordinate to $\mcf$, then $\mcp\sqcap \mcq$ is also a measurable \ae-partition subordinate to $\mcf$ and there is an \ae-unique measurable $c_\blt(\mcq\leftarrow \mcp)\in L^0(M;\RRP)$ such that 

$$\all x\in_\mu M:\restr{\mu^\mcq_x}{\bor(X,\mu;\mcp\sqcap \mcq)}=c_x(\mcq\leftarrow \mcp)\,\restr{\mu^\mcp_x}{\bor(X,\mu;\mcp\sqcap \mcq)},$$

where $\bor(X,\mu;\mcp\sqcap \mcq)$ is the sub-$\sigma$-algebra of the measure algebra of $(X,\mu)$ that corresponds to $\mcp\sqcap \mcq$\fn{Under the Galois correspondence mentioned in \autoref{095}, $\bor(X,\mu;\mcp\sqcap \mcq)\cong_\mu \bor(X,\mu;\mcp) \sqcup \bor(X,\mu;\mcq)$ is the smallest sub-$\sigma$-algebra of the measure algebra of $(X,\mu)$ that contains both $\bor(X,\mu;\mcp)$ and $\bor(X,\mu;\mcq)$.}. In words, the two disintegrations differ by a multiplicative constant that only depends on the fiber \emph{and  a priori also on the basepoint of the fiber}, when restricted to the common refinement of the two \ae-partitions.

\done
\end{lem}

\begin{dfn}\label{092}

Let $\mu$ be a Borel probability measure on $X$ and $\mcf\in\aeFol^s(X,\mu)$ be an \ae-foliation. Whenever there is a measurable \ae-partition $\mcp$ of $(X,\mu)$ subordinate to $\mcf$, we define for $x\in_\mu X$, the \textbf{conditional measure} $\mu^\mcf_x$ to be $\mu^\mcp_x$. In light of \autoref{091}, this defines $\mu^\mcf_\blt$ up to a multiplicative scalar that a priori may depend on the choice of the basepoint.

\done
\end{dfn}

\begin{rem}\label{094}

Note that in \autoref{092} at times $\mu^\mcf_\blt$ is called a system of leaf-wise measures instead of conditional measures to emphasize the difference in the natures of these two objects; we don't make this terminological distinction. Formally $\mu^\mcp_\blt$ is a section of the vertical Radon measure bundle $\vRadon(\pi^\mcp)$, whereas $\mu^\mcf_\blt$ is a section of the projectivization $\PP\vRadon(\pi^\mcp)$ of the vertical Radon measure bundle associated to some subordinated measurable \ae-partition.

The fact that we define $\mu^\mcf_\blt$ with a multiplicative scalar ambiguity and a subordinating measurable \ae-partition ambiguity also means that the support of $\mu^\mcf_x$ is not a cell of any measurable \ae-partition. If a sequence $\mcp_\blt$ of subordinated measurable \ae-partitions have cells that \ae-cover the leaves of $\mcf$, then $\mu^\mcf_\blt$ can be defined as the associated sheafy object. This procedure reduces the two ambiguities mentioned above to the former one, which ambiguity is there to stay (the cone of Haar measures on the additive group $\RR$ is one dimensional). Functionally the phenomenon of \textbf{measure rigidity} is said to be observed when a local choice of a multiplicative scalar can be upgraded to a global choice of a scalar function that is \ae-constant.

\done
\end{rem}

\subsection{Pesin Theory\label{053}}

In this section we define Lusin-Pesin sets and state Pesin's Invariant Manifold Theorem, which ought to be considered as a nonlinear analog of the Oseledets Theorem. Let  $(\mu,\alpha)\in\mfE^r(\RR^k\lact M)$ be a locally free and ergodic system with $\LSpec^\ast(\mu,\alpha)\neq\emptyset$. Let us fix a $C^\infty$ Riemannian metric on $M$ with intrinsic distance function $d_M$ on $M$ and comparison families $C_\blt,K_\blt:\RRP\to L^0(\Osel(\alpha);\RRP)$ as in \autoref{043}.

\begin{dfn}\label{046}

For any $\alpha$-invariant measurable subset $\Pi\subseteq \Osel(\alpha)$, and any $(\eps,\Lambda)\in\RRP^2$ let us define the \textbf{Pesin set} $\UL{\Pi}(\eps,\Lambda)$ by

$$\UL{\Pi}(\eps,\Lambda)=\Pi\cap \bck{C_\eps}(]0,\Lambda])\cap \bck{K_\eps}([1/\Lambda,\infty[).$$

Next, let $\eta\in\left]0,\mu\left(\UL{\Pi}(\eps,\Lambda)\right)\right[$, and for each $n\in \ZZR$ let $\Phi_n:M\to Y_n$ be a measurable function with target a second countable topological space. We define as a \textbf{Lusin-Pesin set} $\Pi(\eps,\Lambda,\eta,\Phi)$ any compact subset of $\UL{\Pi}(\eps,\Lambda)$ with the properties that

\begin{itemize}

\item $\mu\left(\UL{\Pi}(\eps,\Lambda)\right)-\eta< \mu(\Pi(\eps,\Lambda,\eta,\Phi))\leq \mu(\UL{\Pi}(\eps,\Lambda))$,

\item $\all n\in\ZZR: \restr{\Phi_n}{\Pi(\eps,\Lambda,\eta,\Phi)}: \Pi(\eps,\Lambda,\eta,\Phi)\to Y_n$ is continuous.

\end{itemize}

We further assume that for any $\eta_1,\eta_2\in\left]0,\mu\left(\UL{\Pi}(\eps,\Lambda)\right)\right[$, if $\eta_1<\eta_2$ then $\Pi(\eps,\Lambda,\eta_1,\Phi) \supseteq \Pi(\eps,\Lambda,\eta_2,\Phi)$ by replacing the Lusin-Pesin set with parameter $\eta_1$ with the union of the two sets if necessary. Let us also call any $\alpha$-invariant measurable subset $\Pi\subseteq \Osel(\alpha)$ a \textbf{pre-Lusin-Pesin set}.

\done
\end{dfn}

\begin{obs}\label{108}

In the context of \autoref{046}, for any pre-Lusin-Pesin set $\Pi$ and for any $\eps\in\RRP$  we have that strict Lusin-Pesin sets $\UL{\Pi}(\eps,\Lambda)$ increase $\Lambda$ increases. Further we have:

$$\Pi=\bigcup_{\Lambda\in\RRP}\UL{\Pi}(\eps,\Lambda) =_\mu \bigcup_{\Lambda\in\RRP}\bigcup_{\eta\in \left]0,\mu\left(\UL{\Pi}(\eps,\Lambda)\right)\right[} \Pi(\eps,\Lambda,\eta,\Phi),$$

so that Lusin-Pesin sets \ae-cover the underlying pre-Lusin-Pesin set and in particular $\mu$-ae Lyapunov-Perron regular point is in some Lusin-Pesin set (with parameters depending on the point).

\done
\end{obs}

\begin{rem}\label{077}
Note that a Pesin set as defined in \autoref{046} is not necessarily compact nor $\alpha$-invariant and a Lusin-Pesin set is compact but not necessarily $\alpha$-invariant. As a general rule we suppress $\eta$ and $\Phi$ and write $\Pi(\eps,\Lambda)=\Pi(\eps,\Lambda,\eta,\Phi)$; we will also not declare if and where Pesin sets would suffice. Contextually $\eta$ is meant to be small enough and $\Phi$ is meant to contain all functions that are needed to be continuous; there are at most countably many such functions for the purposes of this paper.

An alternative approach to compact Lusin-Pesin sets is to take the closures of Pesin sets\fn{See e.g. \cite[pp.47,81,101]{MR1862379} for this approach.}; we prefer not to do this to stay inside the set $\Osel(\alpha)$ of Lyapunov-Perron regular points.

\done
\end{rem}

\begin{lem}[\fn{\cite[pp.372-373, Cor.5.3]{MR912374}, \cite[pp.91-93, App.A]{MR1858534}.}]\label{054}

Let $\Pi\subseteq \Osel(\alpha)$ be  an $\alpha$-invariant measurable subset. Then for any $(\eps,\Lambda)\in\RRP^2$ and for any $t\in\RR^k[0|=1]$\fn{In accordance with \autoref{103}, $\RR^k[0|=1]$ denotes the unit sphere in $\RR^k$ centered at $0$.}  with $\LSpec(\mu,\alpha_t)\cap \RR_{<0}\neq\emptyset$, we have that $S(\alpha_t)\in C^{(0,+)}\left(\Gr\left(T_{\Pi(\eps,\Lambda)}M\right)\to \Pi(\eps,\Lambda)\right)$, where the latter space is the space of continuous sections of the associated Grassmannian bundle that satisfy a local H\"{o}lder estimate with anonymous exponent.

\done
\end{lem}

\begin{thm}[Pesin\fn{\cite[pp.100-101, 9.3]{MR0466791}, \cite[p.1287, Thm.2.2.1]{MR0458490}, \cite[pp.516, Prop.2.2.1]{MR819556}, \cite[p.195, Thm.16]{MR730270}.}]\label{045}

Let $\Pi\subseteq\Osel(\alpha)$ be an $\alpha$-invariant measurable subset. Then for any $t\in\RR^k[0|=1]$ with $\chi\in \LSpec(\mu,\alpha_t)\cap\RR_{<0}$, $\eps\in\RRP$, and $x\in \Pi$ we have:

\begin{itemize}

\item there is a number $r_{\eps,t}(x)\in \RRP$, 

\item there is  a measurable function $D_{\eps,x,t}:\Pi\to \RRP$, and

\item an embedded $C^r$ submanifold with unique $C^r$ germ at $x$

$$\left(\RR^{\dim S_x(\alpha_{\blt t})}[0| < r_{\eps,t}(x)],0\right)\xrightarrow{\cong_{C^r}} (\mcs_{x,\loc}(\alpha_{\blt t}),x)\hookrightarrow (M,x)$$

such that

\item $T_x\mcs_{x,\loc}(\alpha_{\blt t})=S_x(\alpha_{\blt t})\geq L^\chi_x(\alpha)$,

\item $\all \tau\in\RR:D_{\eps,x,t}(\alpha_{\tau t}(x))\leq D_{\eps,x,t}(x)\, e^{10\eps|\tau|}$,

\item $\all y\in \mcs_{x,\loc}(\alpha_t),\all \tau\in\RRR: d_M(\alpha_{\tau t}(x),\alpha_{\tau t}(y))\leq D_{\eps,x,t}(x) e^{\chi(\tau t)+\eps \tau} d_M(x,y)$.

\end{itemize}

For such a $t\in\RR^k$ we also put $\mcu_{x,\loc}(\alpha_{\blt t})=\mcs_{x,\loc}\left(\alpha_{\blt (-t)}\right)$. Then $T_x\mcu_{x,\loc}(\alpha_{\blt t}) = U_x(\alpha_{\blt t})\leq \bigoplus_{\substack{\rho\in\LSpec(\alpha)\\ \rho\neq \chi}}L^\rho_x(\alpha)$.

Furthermore, 

\begin{itemize}

\item $r_{\eps,t}:\Pi\to \RRP$ is measurable,

\item There is a constant $r^\ast_{\eps,t}\in\RRP$ such that for any $x\in \Pi$ and $\tau\in\RR$ we have $r_{\eps,t}(\alpha_{\tau t}(x))\geq r^\ast_{\eps,t}\,e^{-\eps|\tau|}\,r_{\eps,t}(x)$,

\item $\mcs_{\blt,\loc}(\alpha_{\blt t}): \Pi\to \bigsqcup_{n\in\ZZR} \Emb^r(\RR^n; M)$ is measurable, where the target is endowed with the coproduct Borel $\sigma$-algebra.

\end{itemize}

\done
\end{thm}

\begin{rem}\label{050}

In the context of \autoref{045}, it is clear that if $t\in\RR^k$ is not in $K^\chi$, then one can apply the theorem to $\frac{t}{|t|}$ to construct the local stable and unstable manifolds. Alternatively one can observe that in this case $t$ and $\frac{t}{|t|}$ define the same flow up to a constant time change (see \autoref{052}), thus the local stable and unstable manifolds constructed for them are the same (possibly after one takes $r$ and $D$ smaller). With this understanding we will drop $\blt$ from the notations. We will also conflate a local stable manifold and its Euclidean parameterization.

Finally one may use a $C^2$ Riemannian metric to construct the local stable and unstable manifolds, though this will alter the relationship between the regularity of the manifolds and the regularity of the action. In fact it is reasonable to expect that with minor modifications one ought to be able to construct local stable and unstable manifolds for a $C^{(1,1)}$ Riemannian metric; for a Riemannian metric (or a more general spray) less regular than $C^{(1,1)}$ a more substantial adaptation might be required; see \cite{MR3209534,MR3919953}.

\done
\end{rem}

\begin{dfn}\label{051}
In the context of \autoref{045}, $\mcs_{x,\loc}(\alpha_{t})$ and $\mcu_{x,\loc}(\alpha_{t})$ are called the \textbf{local stable} and \textbf{unstable manifolds} of $\alpha_{\blt t}:\RR\to\Diff^r(M)$ at $x$, respectively. In either case $r_{\eps,t}(x)$ is referred to as the \textbf{size} of the local manifold of $\alpha_{t}$ at $x$.

Further, for any $x\in \Pi$ the \textbf{global stable} and \textbf{unstable manifolds} $\mcs_x(\alpha_{t})$, $\mcu_x(\alpha_{t})$ of $\alpha_{ t}$ at $x$ are defined respectively by

$$\mcs_x(\alpha_{t})= \bigcup_{\tau\in\RRR} \bck{\alpha_{\tau t}}\left(\mcs_{\alpha_{\tau t}(x),\loc}\left(\alpha_{\blt t}\right)\right),\,\, \mcu_x(\alpha_{t})= \mcs_x\left(\alpha_{-t}\right).$$

\done
\end{dfn}

\begin{rem}\label{107}

The global stable and unstable manifolds of $\alpha_t$ at $x$ defined in \autoref{051} can be characterized as global stable and unstable sets of $\alpha_t$ at $x$ of exponential rate:

$$\mcs_x(\alpha_t) = \bigcup_{\lambda\in\RRP}\left\{y\in M \,\left|\, d_M(\alpha_{nt}(y),\alpha_{nt}(x))=O_{n\to\infty}(e^{-\lambda n}) \right\}\right..$$

The global stable and unstable manifolds are injectively immersed $C^r$ submanifolds parameterized by vector spaces of appropriate dimensions, and as submanifolds are independent of the auxilliary Riemannian metric, $\eps$ and $r_{\eps,t}(x)$ and  are uniquely defined.

Also note that since $\RR^k$ is abelian, the time-$t$ map $\alpha_t$ of the action $\alpha$ for any $t\in\RR^k$ permutes the global stable manifolds of $\alpha_{t^\ast}$ for any $t^\ast\in\RR^k$ with $\chi(t^\ast)<0$ for some $\chi$. Similarly the global unstable manifolds are permuted by the action.

\done
\end{rem}

\begin{obs}\label{062}
In the context of \autoref{045},  for any $x,y\in\Pi$ we have $\mcs_{x}(\alpha_t)\cap \mcs_{y}(\alpha_t)\neq\emptyset \ifr \mcs_{x}(\alpha_t)=\mcs_{y}(\alpha_t)$, so that the global stable manifolds partition $\Pi$.

\done
\end{obs}

\begin{dfn}\label{063}
In the context of \autoref{045}, and in light of \autoref{062}, for $\Pi=\Osel(\alpha)$ if $t\in\RR^k$ and $\LSpec(\mu,\alpha_t)\cap \RR_{<0}\neq\emptyset$ we call $\mcs(\alpha_t)$ the \textbf{stable \ae-foliation} of $M$ associated to the time-$t$ map $\alpha_t$ of the action $\alpha$. Similarly if $ \LSpec(\mu,\alpha_t)\cap \RR_{>0}\neq\emptyset$ we call  $\mcu(\alpha_t)=\mcs(\alpha_{-t})$ the \textbf{unstable \ae-foliation} of $M$ associated to the time-$t$ map $\alpha_t$ of the action $\alpha$.

\done
\end{dfn}

\begin{lem}\label{061}
In the context of \autoref{045}, let also $\Lambda\in\RRP$. Then

\begin{itemize}

\item $\inf_{x\in \Pi(\eps,\Lambda)}r_{\eps,t}(x)>0$,

\item $\mcs_{\blt,\loc}(\alpha_t): \Pi(\eps,\Lambda)\to \Emb^r\left(\RR^{s(\eps,\Lambda)}; M\right)$ is in $C^{(0,+)}$, where $s(\eps,\Lambda)\in\ZZP$ is the common dimension of local stable manifolds on the Lusin-Pesin set $\Pi(\eps,\Lambda)$.

\end{itemize}

\done
\end{lem}

\begin{lem}[\fn{\cite[pp.210-211, Prop.3.1]{MR693976}}]\label{097}

In the context of \autoref{045}, for $\Pi=\Osel(\alpha)$ there is a measurable \ae-partition $\mcp\in\aemPar(M,\mu)$ subordinate to $\mcs(\alpha_t)\in\aeFol^r(M,\mu)$ such that

\begin{itemize}

\item $\mcp \to \bck{\alpha_t}(\mcp)$, 

\item $\bigsqcap_{n\in\ZZR}\bck{\alpha_{-nt}}(\mcp) = \Discrete_\mu$, 

\item $\bigsqcup_{n\in\ZZR} \bck{\alpha_{nt}}(\mcp)=\Hull_\mu(\mcs(\alpha_t))$,

\end{itemize}

\done
\end{lem}

Thus if we have an essentially locally free ergodic system $(\mu,\alpha)\in\mfE^r(\RR^k\lact M)$ on a compact $C^\infty$ manifold, and $t\in\RR^k$ is such that at each Lyapunov-Perron regular point $x\in \Osel(\alpha)$, $\alpha_t$ has at least one negative Lyapunov exponent, then $\mcs(\alpha_t)$ is an \ae-partition of $M$ into injectively immersed connected $C^r$ submanifolds, possibly of variable dimension, that depend measurably on the basepoint, with the property that, w/r/t any $C^\infty$ Riemannian metric, on each Lusin-Pesin set  the dependency of the $C^r$ germs on the basepoint is H\"{o}lder and the dimension of the submanifolds is constant. Further, $\mu$ can be disintegrated along leaves of $\mcs(\alpha_t)$ with the caveats discussed in \autoref{094}.

%% file: uzman_arithmeticity_mrpehypotheses.tex
\section{On the Maximal Rank Positive Entropy Hypotheses\label{066}}

In this section we establish some consequences of the maximal rank positive entropy hypotheses of \autoref{001} in terms of Lyapunov geometry. Recall that for $(\mu,\alpha)\in\mfE^r(\RR^k\lact M)$, $\mfe_{(\mu,\alpha)}: \RR^k\to\RRR, t\mapsto \ent_\mu(\alpha_t)$ is the entropy gauge of the system $(\mu,\alpha)$. Note that by an Abramov formula\fn{\cite[p.169, Thm.2]{MR0113985}} the entropy gauge is absolutely homogeneous, i.e., $\all \tau\in\RR, \all t\in\RR^k: \mfe_{(\mu,\alpha)}(\tau t)=|\tau|\,\mfe_{(\mu,\alpha)}(t)$, and by the thesis work of Hu\fn{\cite[p.75, Thm.B]{MR1213080} (= \cite[p.15, Thm.B]{MR2690742}); subadditivity of the entropy gauge in the case of $\mu\ll\leb_M$ was observed earlier in \cite[p.113, Prop.3]{MR677244}.} it is subadditive, i.e. $\all t_1,t_2\in\RR^k: \mfe_{(\mu,\alpha)}(t_1+t_2)\leq\mfe_{(\mu,\alpha)}(t_1)+\mfe_{(\mu,\alpha)}(t_2)$, so that $\mfe_{(\mu,\alpha)}$ is a seminorm of $\RR^k$. Thus we may think of the ergodic theory $\mfE^r(\RR^k\lact M)$ as parameterizing a family of seminorms on $\RR^k$.

\begin{dfn}\label{015}
The \textbf{Fried entropy map } on $\RR^k$  is defined as:

\begin{align*}
\Friedent: \mfE^r(\RR^k\lact M)\to \RRR,\,\,(\mu,\alpha)
\mapsto &\dfrac{\haar_{\RR^k}(\{t\in\RR^k\,|\, |t|_{\ell^1}\leq 1 \})}{\haar_{\RR^k}(\{t\in\RR^k\,|\,\mfe_{(\mu,\alpha)}(t)\leq 1\})}\\
&=  \dfrac{2^k/k!}{\haar_{\RR^k}(\{t\in\RR^k\,|\,\mfe_{(\mu,\alpha)}(t)\leq 1\})},
\end{align*}

where $\haar_{\RR^k}$ is the Haar measure on $\RR^k$ normalized so that $\haar_{\RR^k}([0,1]^k)=1$. If the denominator is infinite we take Fried entropy to be zero.

\done
\end{dfn}

We now prove several straightforward propositions regarding the maximal rank positive entropy hypotheses. In the literature one can find $\ZZ^k$ analogs of versions of these propositions\fn{See e.g. \cite[p.368]{MR2811602}, \cite[p.1209, Prop.3.1]{MR3248484}, \cite[p.136, Prop.A.]{MR3503686}.}. The proofs are mostly linear algebraic.

\begin{obs}\label{100}

Let $(\mu,\alpha)\in\mfE^r(\RR^k\lact M)$ be locally free and ergodic. If $\all t\in\RR^k\setminus 0: \mfe_{(\mu,\alpha)}(t)>0$, then $\Friedent(\mu,\alpha)>0$ by the homogeneity of entropy.

\done
\end{obs}

\begin{prp}\label{016}

Let $(\mu,\alpha)\in\mfE^r(\RR^k\lact M)$ be locally free and ergodic. If $\Friedent(\mu,\alpha)>0$, then $(\mu,\alpha)$ has at least $k+1$ distinct Lyapunov hyperplanes and they are \textbf{in general position}, i.e. 

$$\all I\subseteq \OL{k+1}: \dim\left(\bigcap_{i\in I} K^i\right)=\max\{k-\#(I),0\}.$$

Consequently we also have $2k+1\leq \dim(M)$.

\done
\end{prp}

\begin{pf}

First note that the vanishing of Fried entropy is equivalent to an infinite volume of time parameters with low entropy. Geometrically one can consider the boundary of the unit ball w/r/t $\mfe_{(\mu,\alpha)}$; zero Fried entropy would correspond to affine codimension-$1$ cones parallel to Lyapunov hyperplanes that extend without bounds as parts of the boundary in the case of positive $\mfe_{(\mu,\alpha)}$, or to codimension-$0$ cones in the case of zero $\mfe_{(\mu,\alpha)}$.

Let us first see the case $k=2$ for illustration purposes. If there are less than three hyperplanes, then either there is a quadrant or a halfplane of elements with low entropy, whence Fried entropy vanishes. Thus there are exactly three hyperplanes. If a pair of hyperplanes coincide and have the same orientation, or else the opposite orientation, then there is a half-infinite strip in a quadrant with low entropy; in either case Fried entropy vanishes. Thus the hyperplanes must be distinct, which is equivalent to them being in general position in dimension two. Further they must be oriented so as to see all signatures but $(+,+,+)$ and $(-,-,-)$ (see \autoref{044} below for more on this).

Now consider the general case. If there are less than $k+1$ Lyapunov hyperplanes, then there is a $2^{k}$-ant of elements with zero entropy, whence Fried entropy vanishes. Thus there are at least $k+1$ Lyapunov hyperplanes, say there are $\ell\in\ZZ_{\geq k+1}$ many of them. Say they are not in general position, and let $I\subseteq \OL{\ell}$ be a maximal collection of indices such that the corresponding hyperplanes $\{K^i\,|\, i\in I\}$ are not in general position. Putting $K^I=\bigcap_{i\in I} K^i$, this means that $K^I$ is at least one dimensional and $\dim(K^I)\geq k+1-\#(I)$. Let $i\in \OL{\ell}\setminus I$. Then by the maximality assumption $K^j+ K^I=\RR^{k}$, so that by the dimension formula (sum of the dimensions of the intersection and sum is the sum of the dimensions), we have that $K^j\cap K^I\leq K^I$ is a hyperplane. Thus $\{K^j\cap K^I \,|\, j\in \OL{\ell}\setminus I\}$ is a collection of $\ell-\#(I)\geq k+1-\#(I)$ many hyperplanes of $K^I$ which has at least $k+1-\#(I)$ dimensions. Then there is at the very least a $2^{k}$-ant of elements with zero entropy, whence again Fried entropy vanishes. Thus there has to be at least $k+1$ Lyapunov hyperplanes and they have to be in general position.

\done
\end{pf}

\begin{prp}\label{023}

Let $(\mu,\alpha)\in\mfE^r(\RR^k\lact M)$ be locally free and ergodic. If $(\mu,\alpha)$ has exactly $k+1$ distinct Lyapunov hyperplanes that are in general position, and there exists a $t^\ast\in\RR^k$ such that $\mfe_{(\mu,\alpha)}(t^\ast)>0$, then $\all t\in\RR^k\setminus0: \mfe_{(\mu,\alpha)}(t)>0$.

\done
\end{prp}

\begin{pf}

Suppose there are exactly $k+1$ distinct Lyapunov hyperplanes and they are in general position, and let $t^\ast\in \RR^k\setminus$ be such that $\mfe_{(\mu,\alpha)}(t^\ast)>0$. First we claim that for any chamber time $t\in \Wall^0(\mu,\alpha)$: $\mfe_{(\mu,\alpha)}(t)>0$. By the Ledrappier-Young entropy formula\fn{\cite[p.547, Thm.C]{MR819557}} (and Abramov formula), we have that for any $t\in\RR^k$:

\begin{align*}
\mfe_{(\mu,\alpha)}(t)>0
\iff \all x\in_\mu M: \mu^{\mcs(\alpha_t)}_x \text{ is nonatomic }
\iff \all x\in_\mu M: \mu^{\mcu(\alpha_t)}_x \text{ is nonatomic}.
\end{align*}

Since $\mu^{\mcs(\alpha_{t^\ast})}_\blt$ is nonatomic $\mu$-\ae, there is a $\chi\in \LSpec^\ast(\mu,\alpha)$ such that $L^\chi = S^{\euw(\chi)}\not\leq U(\alpha_{t^\ast})$. Thus $U(\alpha_{t^\ast})\leq U^{\euw(\chi)} = \bigoplus_{\rho\neq\chi} L^\rho$. Since $\mu^{\mcu(\alpha_{t^\ast})}_\blt$ is nonatomic $\mu$-\ae, so is $\mu^{\mcu^{\euw(\chi)}}_\blt$, whence $\all t\in\euw(\chi): \mfe_{(\mu,\alpha)}(t)>0$, and $\mu^{S^{\euw(\chi)}}_\blt$ is nonatomic $\mu$-\ae. Next let $\rho\in\LSpec^\ast(\mu,\alpha)$ be such that $\rho\neq \chi$. Then $U^{\euw(\rho)}=\bigoplus_{\zeta\neq\rho} L^\zeta\geq L^\chi = S^{\euw(\chi)}$, so that $\mu^{U^{\euw(\rho)}}_\blt$ is nonatomic $\mu$-\ae, hence we have that for any chamber time $t$, $\mfe_{(\mu,\alpha)}(t)>0$.

To extend the statement to wall times, first note that by the assumption on Lyapunov hyperplanes and wall time of order $i$ is an element of a $k-i$ dimensional intersection that is not an element of a $k-i+1$ dimensional intersection. In particular, $\Wall^k=\emptyset$ and $\Wall^{k+1}=\{0\}$.

Let $I\subseteq \OL{k+1}$ and put $K^I=\bigcap_{i\in I} K^i$. Let $t\in \Wall^I$ and suppose $\all j\in\OL{k+1}\setminus I: \chi^j(t)<0$, so that the signature at $t$ is $(\underbrace{0,\cdots,0}_{\#(I) \text{ many}}, \underbrace{-,\cdots,-}_{k+1-\#(I)\text{ many}})$. We claim that there is a perturbation of $t$ such that the signature becomes $( -,\cdots,-)$, hence giving a contradiction to all chamber times having positive entropy. First note that by the general position assumption for $\eps\in\RRP$ small enough any element in $\RR^k[t | <\eps]$ would have a signature of the form $(\epsilon_1,\cdots,\epsilon_{\#(I)}, \underbrace{-,\cdots,-}_{k+1-\#(I)\text{ many}})$. Then by induction one can turn each $0$ into a $-$ one by one, as by the general position assumption $\{K^j\cap K^I | j\in\OL{k+1}\setminus I\}$ is a collection of hyperplanes in $K^I$ that are in general position.

\done
\end{pf}

\begin{cor}\label{044}
Let $(\mu,\alpha)\in\mfE^r(\RR^k\lact M)$ be essentially locally free and ergodic. If $(\mu,\alpha)$ is an MRPES, then there are exactly $2^{k+1}-2$ Weyl chambers parameterized by signatures $(\epsilon_1,\epsilon_2,...,\epsilon_{k+1})$ for $\epsilon_i\in\{\pm\}$ except for $(+,+,\cdots,+)$ and $(-,-,\cdots,-)$. In particular there is a unique injective map

$$\euw: \LSpec^\ast(\mu,\alpha)\to\Cham(\mu,\alpha)$$

such that 

$$\all \chi\in\LSpec^\ast(\mu,\alpha): L^\chi=S^{\euw(\chi)}=U^{-\euw(\chi)}, \bigoplus_{\substack{\rho\in \LSpec^\ast(\mu,\alpha)\\\rho\neq\chi}} L^\rho=U^{\euw(\chi)}=S^{-\euw(\chi)}.$$

\done
\end{cor}

\begin{pf}

We use induction. For $k=2$ the statement is straightforward; the fact that we don't see $(+,+,+)$ or $(-,-,-)$ is due to the positive entropy assumption. Say the statement is true for $k=m$. If $K^1,K^2,...,K^m,K^{m+1}\leq \RR^m$ were in general position, then so would $K^1,K^2,...,K^m$. Then these first $m$ hyperplanes would separate $\RR^m$ into $2^m$ chambers and all signatures $(\epsilon_1,\cdots,\epsilon_m)$ would be available ($\dagger$). Again by the entropy assumption $K^{m+1}$ must be positioned so as to not separate the chambers with signatures $(+,\cdots,+)$ and $(-,\cdots,-)$ and bisect all other chambers. This gives $(2^m-2)2+2=2^{m+1}-2$ chambers with signatures $(\epsilon_1,\cdots,\epsilon_m,+)$, $(\epsilon_1,\cdots,\epsilon_m,-)$, $(+,\cdots,+,-)$, $(-,\cdots,-,+)$.

($\dagger$) also easily follows from induction. Indeed, the statement is clearly true for $k=2$, and if the statement is true for $k=m$, then for $k=m+1$, by the general position assumption $K^1\cap\cdots \cap K^m$ is one dimensional and $K^1\cap\cdots\cap K^m\cap K^{m+1}$ is zero dimensional. Thus $K^1,\cdots,K^m$ splits $\RR^m$ into $2^m$ chambers and $K^{m+1}$ bisects each one of these $2^m$ chambers. 

\done
\end{pf}

\begin{obs}\label{102}

As mentioned in \autoref{101}, we see that by \autoref{100}, \autoref{016} and \autoref{023} the third hypothesis of \autoref{001} is redundant; in particular among locally free ergodic $\RR^k$ systems on $2k+1$ dimensional manifolds for $k\in\ZZ_{\geq2}$, being a maximal rank positive entropy system is a $\cat{Meas}$-isomorphism invariant. Alternatively, one can reformulate the hypotheses of \autoref{001} to be more in line with the standard hypothesis scheme in measure rigidity: a locally free ergodic $\RR^k$ system is an MRPES iff there are exactly $k+1$ distinct Lyapunov hyperplanes that are in general position and at least one time-$t$ map of the system has positive entropy.

\done
\end{obs}

The next lemma addresses suspension hereditarity and time change hereditarity for the MRPES hypotheses.

\begin{lem}\label{018}

Let $N$ be a compact $C^\infty$ manifold. If $(\nu,\beta)\in\mfE^r(\ZZ^k\lact N)$ is a $\ZZ^k$ MRPES and $\kappa: \RR^k\times (\RR^k\otimes^\beta N)\to \RR^k$ is a $C^r$ time change of $\hbar^\beta$, then $\left((\haar_{\TT^k}\otimes^\beta \nu)_\kappa,\hbar^{\beta}_\kappa\right)$ is an $\RR^k$ MRPES.

\done
\end{lem}

\begin{pf}

By definition $(\nu,\beta)$ is a $\ZZ^k$ MRPES iff its suspension system $(\mu,\alpha)$ is an $\RR^k$ MRPES. Thus it suffices to show that being an MRPES is invariant under time changes. It's clear that a time change $\kappa$ does not change the dimension. After \autoref{104} it's also clear that ergodicity is preserved, and if $x\in M$ and $t\in\Stab_x(\alpha_\kappa)$, then $\kappa(t,x)\in \Stab_x(\alpha)$, and since $\kappa(\blt,x)$ is a diffeomorphism and $\Stab_x(\alpha)$ is discrete, $\Stab_x(\alpha_\kappa)$ is discrete so that local freeness is also preserved. Finally by \autoref{105}, 

$$\all t\in\RR^k\setminus0: \mfe_{(\mu,\alpha)}(t)>0 \iff \all t\in\RR^k\setminus0: \mfe_{(\mu_\kappa,\alpha_\kappa)}(t)>0,$$

and after \autoref{102} this is sufficient for the time change of the suspension $(\mu_\kappa,\alpha_\kappa)$ to be an $\RR^k$ MRPES.

\done
\end{pf}

\begin{rem}\label{137}

After \autoref{136}, if in \autoref{018} we take $\kappa$ to be a $C^s$ \ae-time change, then 

$$\left((\haar_{\TT^k}\otimes^\beta \nu)_\kappa,\hbar^{\beta}_\kappa\right)$$

is still an $\RR^k$ MRPES, although now possibly it's not a system of diffeomorphisms.

\done
\end{rem}

\begin{rem}\label{064}
In the context of \autoref{044}, and in light of the comments at the end of \autoref{053}, we have that for any $\chi\in \LSpec^\ast(\mu,\alpha), \all t\in \euw(\chi)\neq\emptyset: L^\chi(\alpha)=S(\alpha_t)$, whence we may define the \textbf{(fine) Lyapunov \ae-foliation} of $\alpha$ associated to $\chi$ by

$$\mcl^\chi(\alpha)=\mcs(\alpha_t).$$

We will also write $\mcl^\chi_{x,\loc}(\alpha)=\mcs_{x,\loc}(\alpha_t)$. Note that this notation is doubly ambiguous, as the size of any local stable manifold is suppressed, and the size of the local Lyapunov manifold is dependent on the time parameter $t$. 

Thus in the context of \autoref{001}, we have, for each nonzero Lyapunov exponent $\chi$ of the system $(\mu,\alpha)$, an $\alpha$-invariant measurable \ae-foliation $\mcl^\chi$ of $M$ each of whose leaves is a $1$-dimensional injectively immersed and connected $C^r$ submanifold, and on any Lusin-Pesin set w/r/t any $C^\infty$ Riemannian metric the $C^r$ germs of leaves of $\mcl^\chi$ depend H\"{o}lder continuously. Further, for $x\in_{\mu} M$, there is a Radon measure $\mu^\chi_x$ on $\mcl^\chi_x(\alpha)$ that is well defined up to a multiplicative scalar that a priori depends on the basepoint, even along the same leaf, such that

$$\all \phi\in L^0_b(M;\RR): \int_M \phi(x)\, d\mu(x) = \int_M \int_{\mcl^\chi_x} \phi(y)\, d\mu^\chi_x(y)\, d\mu(x).$$

Here the righthand side is ambiguous; see \autoref{094}.

\done
\end{rem}

We make one final observation regarding the orbit type of an MRPES.

\begin{lem}\label{114}

Let $(\mu,\alpha)\in \mfE^r(\RR^k\lact M)$ be ergodic. If $\all t\in\RR^k\setminus0:\mfe_{(\mu,\alpha)}(t)>0$, then $(\mu,\alpha)$ is essentially free, that is,

$$\all x\in_\mu M: \alpha_\blt(x): \RR^k\xrightarrow{\cong_{C^r}}\mco_x.$$

\done
\end{lem}

\begin{pf}

Let us consider the collection $\mfC(\RR^k)$ of closed subgroups of $\RR^k$ endowed with the Chabauty topology\fn{See \cite{delaharpepierre2008}.}. Then $\mfC(\RR^k)$ is a compact Polish space and we have

$$\Stab_\blt(\alpha):M\to \mfC(\RR^k).$$

Using the basis elements of the Chabauty topology it's straightforward to verify that $\Stab_\blt(\alpha)$ is Borel measurable, hence by ergodicity it's constant $\mu$-\ae. If $(\mu,\alpha)$ weren't essentially free then $\exi t^\ast\in\RR^k\setminus0,\all x\in_\mu M: t^\ast\in \Stab_x(\alpha)$, i.e. $\alpha_{t^\ast}=\id_M$, whence $\mfe_{(\mu,\alpha)}(t^\ast)=0$, a contradiction.

\done
\end{pf}

%% file: uzman_arithmeticity_affinestructures.tex
\section{Affine Structures for Leaves of Invariant {\ae-Foliations}\label{112}}

From now on until the end of the paper let $(\mu,\alpha)\in\mfE^r(\RR^k\lact M)$ be an MRPES as in \autoref{001}. Recall that therefore by \cite[p.363, Main Thm.(2)]{MR2811602} we have that $\mu\ll\leb_M$; by \autoref{044} and \autoref{064}, all nonorbital Lyapunov exponents are simple and non-zero, no two distinct nonorbital Lyapunov exponents are neither positively nor negatively proportional, and there is a unique map

$$\euw: \LSpec^\ast(\mu,\alpha)\to \Cham(\mu,\alpha),\,\, \chi^i\mapsto \euw(\chi^i)=\euw^i$$

with the property that $\all i\in\OL{l}, \all t\in \euw^i: L^i(\alpha)=S^{\euw^i}(\alpha)=S(\alpha_t)$. We also fix a $C^\infty$ Riemannian metric $\mfg$ on $M$ with comparison families $C_\blt, K_\blt: \RRP\to L^0(\Osel(\alpha);\RRP)$ (see \autoref{042}) and intrinsic distance function $d_M$ on $M$. For an \ae-foliation $\mcl$ of $M$ denote by $d_{\mcl_x}$ the intrinsic distance function on the leaf $\mcl_x$ induced by the induced Riemannian metric $\mfg^\mcl_x$ on $\mcl_x$.

\subsection{Nonstationary Linearizations for Lyapunov \ae-Foliations\label{067}}

In this subsection we'll define the nonstationary linearizations for Lyapunov \ae-foliations and discuss its main properties. The main strategy is to use the fact that each leaf of each Lyapunov \ae-foliation is one dimensional; whence an affinely chosen real number completely characterizes the position of any point on any leaf once a basepoint on that leaf is chosen. From a structural point of view this is also the first instance the $\pm$-ambiguity in the space crystal in \autoref{001} can be observed; for this reason the account below is written in invariant language and identifications are kept track of pedantically. We start with an adaptation of a well-known lemma to our case\fn{Compare e.g. \cite[p.151]{MR0224771}, \cite[p.34, Thm.4]{MR0399421}, \cite[pp.418-419, Prop.2]{MR721733}, \cite[p.573]{MR840722}, \cite[p.298, Cor.4.4]{MR1194019}, \cite[p.658, Lem.3.1]{MR1452186},
\cite[pp.131-134]{MR2261075}, \cite[p.572, Lem.2.2]{MR840722}, \cite[p.298, p.303]{MR1194019}, \cite[p.659]{MR1452186}, \cite[p.5]{dolgopyatd2001}.}:

\begin{prp}\label{076}

Let $\chi\in\LSpec^\ast(\mu,\alpha)$, $t\in\euw(\chi)$, $\eps\in\RRP$ be such that $\eps<-\chi(t)$. Define for any $x\in_\mu M$

$$\Delta_x^\chi = \Delta_{x,t}^\chi:(\mcl_{x,\loc}^\chi,x)\to (\RRP,1),\,\, z\mapsto \lim_{n\to \infty} \dfrac{\left\Vert T_z^\chi \alpha_{nt}\right\Vert}{\left\Vert T_x^\chi \alpha_{nt}\right\Vert}=\prod_{n\in\ZZR}\dfrac{\left\Vert T_{\alpha_{nt}(z)}^\chi \alpha_{t}\right\Vert}{\left\Vert T_{\alpha_{nt}(x)}^\chi \alpha_{t}\right\Vert},$$

where $T^\chi\alpha_s=\restr{T\alpha_s}{L^\chi}$. Then for any $x\in_\mu M$, $\Delta_x^\chi$ is well-defined. Moreover,

\begin{enumerate}

\item $\all x\in_\mu M: \Delta_x^\chi\in C^{(0,\theta)}(\mcl_{x,\loc}^\chi,x;\RRP,1)$.

\item $\all x\in_\mu M, \all z_1,z_2\in \mcl_{x,\loc}^\chi: \dfrac{\Delta_x^\chi(z_1)}{\Delta_x^\chi(z_2)}=\Delta_{z_2}^\chi(z_1)$. In particular also $\Delta_{z_2}^\chi(z_1) = \dfrac{1}{\Delta_{z_1}^\chi(z_2)}$ and $\Delta_{z_2}^\chi(x)\Delta_x^\chi(z_1)=\Delta_{z_2}^\chi(z_1)$.

\item $\all x\in_\mu M, \all z\in \mcl_{x,\loc}^\chi, \all s\in\RR^k: \Delta_{\alpha_{s}(x)}^\chi(\alpha_{s}(z))   \left\Vert T_z^\chi \alpha_{s}\right\Vert   =\left\Vert T_x^\chi \alpha_{s}\right\Vert  \Delta_x^\chi(z)  $.

\item $\Delta^\chi_\blt(\blt):\{(x,z)\in M\times M\,|\, x\in \Osel(\alpha), z\in\mcl_{x,\loc}^\chi\}\to \RRP$ is measurable and everywhere defined.

\item On any Lusin-Pesin set the limit defining $\Delta^\chi_\blt(\blt)$ is uniform. More precisely, for any pre-Lusin-Pesin set $\Pi$, any $\Lambda\in\RRP$, any $r_0\in\RRP$ that bounds from below the sizes of local stable manifolds of $\alpha_t$ on $\Pi(\eps,\Lambda)$, and $\eta\in\RRP$, there is an $N\in\ZZP$ such that for any $n\in\ZZ_{\geq N}$:

$$\sup_{x\in \Pi(\eps,\Lambda)} \sup_{z\in \mcl_{x,\loc}^\chi[x|\leq r_0]} \left|\, \prod_{i\in\UL{n}}\dfrac{\left\Vert T_{\alpha_{it}(z)}^\chi \alpha_{t}\right\Vert}{\left\Vert T_{\alpha_{it}(x)}^\chi \alpha_{t}\right\Vert} -\Delta^\chi_x(z)\,\right| <\eta. $$

\item On any Lusin-Pesin set $\Delta^\chi_\blt(\blt)$ is uniformly continuous. More precisely, for any pre-Lusin-Pesin set $\Pi$, any $\Lambda\in\RRP$, any $r_0\in\RRP$ that bounds from below the sizes of local stable manifolds of $\alpha_t$, $\Delta^\chi_\blt(\blt):\{(x,z)\in M\times M\,|\, x\in \Pi(\eps,\Lambda), z\in\mcl_{x,\loc}^\chi[x|\leq r_0]\}\to \RRP$ is uniformly continuous.

\item On any Lusin-Pesin set $\Lambda_\blt^\chi$ is uniformly continuous. More precisely, for any pre-Lusin-Pesin set $\Pi$, any $\Lambda\in\RRP$, any $r_0\in\RRP$ that bounds from below the sizes of local stable manifolds of $\alpha_t$, w/r/t the Euclidean coordinates on local manifolds we have that $\Delta_\blt^\chi: \Pi(\eps,\Lambda)\to C^{(0,\theta)}(\RR[0|\leq r_0],0;\RRP,1)$ is uniformly continuous.

\end{enumerate}

\done
\end{prp}

\begin{pf}
Since $t \in\euw(\chi)$ we have $L^\chi(\alpha)=S(\alpha_t)$. Let us first see that $\Delta_x^\chi$ is well-defined. If $z\in\mcl_{x,\loc}^\chi$ (recall \autoref{064}), then

\begin{align*}
&\left|\,\left\Vert T_{\alpha_{nt}(z)}^\chi\alpha_t\right\Vert - \left\Vert T_{\alpha_{nt}(x)}^\chi\alpha_t\right\Vert\,\right|\\
 &\phantom{T_{\alpha_{nt}(x)}\alpha_t:} =\left|\,\left\Vert T_{\alpha_{nt}(z)}\alpha_t:L_{\alpha_{nt}(z)}^\chi\to T_{\alpha_{(n+1)t}(z)}M\right\Vert\right.\\
&\phantom{T_{\alpha_{nt}(x)}\alpha_t:}\phantom{T_{\alpha_{nt}(x)}\alpha_t:} - \left.\left\Vert T_{\alpha_{nt}(x)}\alpha_t:L_{\alpha_{nt}(x)}^\chi\to T_{\alpha_{(n+1)t}(x)}M\right\Vert\,\right|\\
&\phantom{T_{\alpha_{nt}(x)}\alpha_t:} \leq\left|\,\left\Vert T_{\alpha_{nt}(z)}\alpha_t:L_{\alpha_{nt}(z)}^\chi\to T_{\alpha_{(n+1)t}(z)}M\right\Vert \right.\\
&\phantom{T_{\alpha_{nt}(x)}\alpha_t:}\phantom{T_{\alpha_{nt}(x)}\alpha_t:}- \left.\left\Vert T_{\alpha_{nt}(x)}\alpha_t:\fwd{\Pi_{\alpha_{nt}(x)\leftarrow \alpha_{nt}(z)}}(L_{\alpha_{nt}(z)}^\chi)\to T_{\alpha_{(n+1)t}(x)}M\right\Vert\,\right|\\
&\phantom{T_{\alpha_{nt}(x)}\alpha_t:}+ \left|\,\left\Vert T_{\alpha_{nt}(x)}\alpha_t:\fwd{\Pi_{\alpha_{nt}(x)\leftarrow \alpha_{nt}(z)}}(L_{\alpha_{nt}(z)}^\chi)\to T_{\alpha_{(n+1)t}(x)}M\right\Vert\right.\\
&\phantom{T_{\alpha_{nt}(x)}\alpha_t:} \phantom{T_{\alpha_{nt}(x)}\alpha_t:}- \left.\left\Vert T_{\alpha_{nt}(x)}\alpha_t:L_{\alpha_{nt}(x)}^\chi\to T_{\alpha_{(n+1)t}(x)}M\right\Vert\,\right|,
\end{align*}

where $\Pi_{q\leftarrow p}$ is the parallel transport from $T_pM$ to $T_qM$.  In the final expression the first difference has a $\theta$-H\"{o}lder modulus of continuity since $\alpha_t$ is $C^r$ for $r=(1,\theta)$ and parallel transports are $C^\infty$. Similarly the second difference has a $\theta$-H\"{o}lder modulus of continuity since $L_{\alpha_{nt}(z)}^\chi=T_{\alpha_{nt}(z)}\mcl_{\alpha_{nt}(x),\loc}^\chi$, $\mcl_{\alpha_{nt}(x),\loc}^\chi$ is a  $C^r$ manifold by \autoref{045}, $\alpha_{nt}$ is $C^r$ and parallel transports are $C^\infty$. Again by \autoref{045} we thus have

$$\left|\,\left\Vert T_{\alpha_{nt}(z)}^\chi\alpha_t\right\Vert - \left\Vert T_{\alpha_{nt}(x)}^\chi\alpha_t\right\Vert\,\right|
\lesssim_{t,x,\eps} e^{(\chi(nt)+\eps n)\theta}d_M(z,x)^\theta
.$$

Thus

$$\left|\,\dfrac{\left\Vert T_{\alpha_{nt}(z)}^\chi \alpha_{t}\right\Vert}{\left\Vert T_{\alpha_{nt}(x)}^\chi \alpha_{t}\right\Vert}-1\,\right|
\lesssim_{t,x,\eps} d_{C^1}(\id_M,\alpha_t) e^{(\chi(nt)+\eps n)\theta}d_M(z,x)^\theta,$$

and so

$$\prod_{n\in\ZZR} \dfrac{\left\Vert T_{\alpha_{nt}(z)}^\chi \alpha_{t}\right\Vert}{\left\Vert T_{\alpha_{nt}(x)}^\chi \alpha_{t}\right\Vert}
\leq \exp\left(\sum_{n\in\ZZR} \left|\,\dfrac{\left\Vert T_{\alpha_{nt}(z)}^\chi \alpha_{t}\right\Vert}{\left\Vert T_{\alpha_{nt}(x)}^\chi \alpha_{t}\right\Vert}-1\,\right|\right)<\infty,$$

so that $\Delta_x^\chi$ is well-defined and as stated in the first item has a $\theta$-H\"{o}lder modulus of continuity. The second and third items are syntactic (for item (ii) note that the above argument works for $x$ replaced with any other point in $\mcl^\chi_x$ and item (iii) is valid germinally in that one may need to reduce the size of the local manifolds, e.g. if $s$ and $t$ is in the same chamber). The fourth item follows from the fact that local Lyapunov manifolds depend measurably on the basepoint by \autoref{045}, as in our case we may choose $\Pi=_\mu M$, so that $\Delta^\chi$ is defined as the limit of the restrictions,  to a measurable subset,  of functions $M\times M\to \RRP$ that are $C^{(0,\theta)}$ in each component. Fifth and sixth items follow from the fact that on a Lusin-Pesin set with $r_0\in\RRP$ a lower bound on the local Lyapunov manifolds we have, for $z\in\mcl_{x,\loc}^\chi[x|\leq r_0]$:

$$\left|\,\dfrac{\left\Vert T_{\alpha_{nt}(z)}^\chi \alpha_{t}\right\Vert}{\left\Vert T_{\alpha_{nt}(x)}^\chi \alpha_{t}\right\Vert}-1\,\right|
\lesssim_{t,\eps,\Lambda}^{x,z} d_{C^1}(\id_M,\alpha_t) e^{(\chi(nt)+\eps n)\theta}r_0^\theta.$$

Here $\lesssim_A^B$ means that the inequality holds up to a multiplicative constant that may depend on $A$ but does not depend on $B$. The seventh item follows from the previous two items (here we conflate an embedded submanifold with its parameterization; see \autoref{050}).

\done
\end{pf}

\begin{rem}\label{078}
Note that in \autoref{076} to define $\Delta^\chi_x$ one could alternatively use the limit

$$z\mapsto \lim_{\tau\to \infty} \dfrac{\left\Vert T_z^\chi \alpha_{\tau t}\right\Vert}{\left\Vert T_x^\chi \alpha_{\tau t}\right\Vert}=\int_0^\infty\dfrac{\left\Vert T_{\alpha_{\tau t}(z)}^\chi \alpha_{t}\right\Vert}{\left\Vert T_{\alpha_{\tau t}(x)}^\chi \alpha_{t}\right\Vert}\, d\tau.$$

Also note that even though we suppressed the dependency of $\Delta^\chi_x=\Delta^\chi_{x,t}$ on $t$, it may very well depend on $t$ beyond the ambiguities mentioned in \autoref{064}.

\done
\end{rem}

\begin{rem}\label{138}

In \autoref{076} if instead of $r=(1,\theta)$, we had $r=(q,\theta)$ for $q\in\ZZ_{\geq2}$, one could argue as follows to get the more specific $C^{(q-1,\theta)}$ regularity instead of $C^{(0,\theta)}$ in the first item: For the sake of readability we switch to Euclidean coordinates, suppress different spaces of sections involved and treat all objects simply as functions, and apply a soft Fa\`{a} di Bruno argument. For $\mce$ an anonymous collection of functions, which we'll call a collection of elementary functions, and $S,P,C\in\ZZP$, denote by $\sum_S\prod_P\bigcirc_C(\mce)$ the collection of all functions that can be written as the sum of not more than $S$ summands, each of which can be written as the product\fn{Note that intrinsically "products" really are contractions of tensors.} of not more than $P$ factors, each of which can be written as the composition of not more than $C$ elements in $\mce$. For a collection $\mcf$ of functions and $s\in\ZZR$ let's also denote by $\mcf^{(s)}$ the collection of functions that are the $s$-th derivatives of functions in $\mcf$, and put $\mcf^{(\leq s)}=\bigcup_{k\in\UL{s+1}}\mcf^{(k)}$. Then a straightforward induction argument shows that

$$\left[\sum_S\prod_P\underset{C}{\bigcirc}(\mce)\right]^{(s)}
\subseteq \sum_{S\prod_{k\in\UL{s}}(P+kC)}\prod_{P+sC}\underset{C}{\bigcirc}\left(\mce^{(\leq s)}\right)
\subseteq \sum_{S(P+C)^s(s-1)!}\prod_{P+sC}\underset{C}{\bigcirc}\left(\mce^{(\leq s)}\right).$$

In our context the collection $\mce$ will correspond to $\left\{\left. T^\chi_{\alpha_{nt}(z)}\alpha_t\,\right|\, z\in \mcl_{x,\loc}^\chi,n\in\ZZR\right\}$, so that e.g. $T^\chi_z\alpha_{nt}\in\sum_1\prod_1\bigcirc_n(\mce)\cap\sum_1\prod_n\bigcirc_1(\mce)$. As such, for the $q$-th derivatives it suffices to control $\left[\sum_1\prod_n\bigcirc_{1}(\mce)\right]^{(q)}$. The strategy is to use the exponential decay of elements in $\mce$ to offset the growth due to the sums of products of functions in $\mce^{(\leq q)}$. Indeed, each summand of $\left[\sum_1\prod_n\bigcirc_{1}(\mce)\right]^{(q)}$ has at most $q$ factors with leading term a derivative of order more than $1$, and there are $n+q$ factors of each summand, whence for $\theta'\in]0,\theta]$ we have the following uniform upper bound (which may depend on $t,x,$ and $\eps$) for differential objects of same type at $x$ and $z$:

$$ d_{C^q}(\id_M,\alpha_t) (n+1)^q (q-1)! \,e^{(\chi(nt)+\eps n)\theta'}d_M(z,x)^{\theta'}.$$

\done
\end{rem}

\begin{prp}\label{017}

In the context of \autoref{076}, for any $x\in_\mu M$ there is a unique $\Lambda_x^\chi$ in  

$$\Diff^r\left(L_x^\chi,0\,;\,\mcl_x^\chi,x\right)$$

such that

\begin{enumerate}

\item $\all x\in_\mu M: T_0^\chi\Lambda^\chi_x\cong\id_{L_x^\chi}: L_x^\chi\cong T_0L_x^\chi\to T_x\mcl_x^\chi= L_x^\chi$,

\item $\all x\in_\mu M,\all v\in L_x^\chi,\all s\in\RR^k:\left|\Lambda^\chi_{\alpha_{s }(x)}\circ T_x^\chi\alpha_{s }(v)\right|=\left|\alpha_{s }\circ \Lambda^\chi_x(v)\right|$,

\item W/r/t Euclidean coordinates on global manifolds $\Lambda^\chi_\blt: M\to \Imm^r(\RR;M)$ is measurable,

\item On any Lusin-Pesin set $\Lambda^\chi_\blt$ is uniformly continuous. More precisely, for any pre-Lusin-Pesin set $\Pi$ and for any $\Lambda\in\RRP$, w/r/t Euclidean coordinates on global manifolds we have that $\Lambda^\chi_\blt:\Pi(\eps,\Lambda)\to \Imm^r(\RR;M)$ is  uniformly continuous, where the target is endowed with the $C^r$ compact-open topology,

\item $\all x\in_\mu M, \all y\in\mcl^\chi_x:\left(\Lambda^\chi_y\right)^{-1}\circ\Lambda^\chi_x:L_x^\chi\to L_y^\chi$ is affine. 

\end{enumerate}

\done
\end{prp}

\begin{pf}

Fix a $t \in\euw(\chi)$ so that $L^\chi(\alpha)=S(\alpha_t)$, let $\eps\in\RRP$ be such that $\eps<-\chi(t)$. Let us also pick a section $\xi\in L^0(L^\chi[0|=1]\to M)$\fn{In accordance with \autoref{103}, this means that $\xi$ is an \ae-vector field along $L^\chi$ with unit norm.} to act as an orientation in such a way that the identification $L^\chi\cong_\mu M\times \RR$, $c\xi_x\mapsto c$ preserves the orientation of each fiber. Note that $\xi$ also naturally orients each leaf of $\mcl^\chi$. We'll first define a map $\Phi_{x,\loc}: \mcl^\chi_{x,\loc}\to L^\chi_x$ and then extend it to a map $\Phi_x:\mcl^\chi_x\to L^\chi_x $;   $\Lambda^\chi_x$ will then be defined as the inverse of $\Phi_x$. Define $\Phi_{x,\loc}=\Phi_{x,\loc,\xi}: (\mcl^\chi_{x,\loc},x)\to (L_x^\chi,0)$ by

$$y\mapsto \begin{ncase} \left(\int_x^y \Delta_x^\chi(z)\,dz\right)\xi_x &\text{, if the segment of $\mcl^\chi_x$ from $x$ to $y$ is positively oriented}\\
-\left(\int_y^x \Delta_x^\chi(z)\,dz\right) \xi_x  &\text{, if the segment of $\mcl^\chi_x$ from $x$ to $y$ is negatively oriented}\end{ncase},$$

where the integrant is the function established in \autoref{076} and the integral is w/r/t the Riemannian metric on $\mcl^\chi_x$ induced by $\mfg$, so that under the identification $L^\chi_x\cong \RR$, $\Delta_x^\chi(y)$ is the derivative of $\Phi_{x,\loc}$ along $\mcl^\chi_{x,\loc}$ evaluated at $y$; in particular $\Phi_{x,\loc}(x)=0 \in T_xM$, $T_x^\chi\Phi_{x,\loc}\cong\id_{L^\chi_x}$ ($\dagger_1$) and $\Phi_{x,\loc}$ is a $C^r$ diffeomorphism onto its image:

\begin{center}
\begin{tikzcd}
{(\mcl_{x,\loc},x)} \arrow[r, "{\Phi_{x,\loc}}"] \arrow[rd, "{\int_x^\blt \Delta_x(z)\, dz}"'] & {(L_x,0)} \arrow[d, "\cong"] \\
                                                                                               & {(\RR,0)}                   
\end{tikzcd}

\end{center}

By item (iii) of \autoref{076} using the Euclidean parameterizations for local Lyapunov manifolds we have for any $s\in\RR^k$

\begin{align*}
\dfrac{d}{dy}\left[\int_{\alpha_s(x)}^{\alpha_s(y)} \Delta_{\alpha_s(x)}(z) \,dz\right]
&=\Delta_{\alpha_s(x)}(\alpha_s(y)) \left\Vert T_y^\chi\alpha_s\right\Vert\\
&= \left\Vert T_x^\chi\alpha_s\right\Vert \Delta_x^\chi(y)
= \dfrac{d}{dy}\left[\left\Vert T_x^\chi\alpha_s\right\Vert\int_x^y\Delta_x(z)\, dz\right]
\end{align*}

Taking the antiderivative of this equation w/r/t $y$ and evaluating at $x$ we obtain

$$|\Phi_{\alpha_s(x),\loc}\circ(\alpha_s(y))|=|T_x^\chi\alpha_s \circ(\Phi_{x,\loc}(y))|\quad\quad (\dagger_2),$$

that is, we have a diagram that commutes up to $\pm$:

\begin{center}
\begin{tikzcd}
{\mcl_{x,\loc}^\chi} \arrow[rr, "\alpha_s"] \arrow[dd, "{\Phi_{x,\loc}}"'] &     & {\mcl_{\alpha_s(x),\loc}^\chi} \arrow[dd, "{\Phi_{\alpha_s(x),\loc}}"] \\
                                                                           & \pm &                                                                        \\
L_x^\chi \arrow[rr, "T_x^\chi\alpha_s"']                                   &     & L_{\alpha_s(x)}^\chi                                                  
\end{tikzcd}

\end{center}

Next we extend $\Phi_{x,\loc}$ to a map $\Phi_x:(\mcl_x^\chi,x)\to (L_x^\chi,0)$. If $y\in \mcl_x^\chi$, then put $n_y=\inf\left\{n\in\ZZR\,\left|\, \alpha_{n t}(y)\in \mcl_{\alpha_{n t}(x),\loc}^\chi\right.\right\}\in\ZZR$. Note that as $\mcl^\chi=\mcs(\alpha_t)$, by \autoref{045} even if the sizes of local manifolds may decrease $n_\blt:\mcl_x^\chi\to\ZZR$ is a well-defined locally constant function. Thus we may define $\Phi_x(y)= \left(T_x^\chi \alpha_{n_y t}\right)^{-1}\circ \Phi_{\alpha_{n_y t(x),\loc}}\circ \alpha_{n_y t}(y)$. Further, by ($\dagger_2$) we have:

$$\all y\in \mcl_x^\chi, \all n\in\ZZ_{\geq n_y}:
\left(T_x^\chi \alpha_{n t}\right)^{-1}\circ \Phi_{\alpha_{n t(x),\loc}}\circ \alpha_{n t}(y)
=\left(T_x^\chi \alpha_{n_y t}\right)^{-1}\circ \Phi_{\alpha_{n_y t(x),\loc}}\circ \alpha_{n_y t}(y),$$

so that the definition of $\Phi_x$ on an open ball in $\mcl_x$ restricted to a smaller ball in $\mcl_x$ coincides with the definition of $\Phi_x$ on said smaller ball; consequently $\Phi_x:(\mcl_x^\chi,x)\to (L_x^\chi,0)$ is well-defined. In Euclidean coordinates $\Phi_x$ is strictly increasing, whence $\Phi_x$ is a diffeomorphism onto its image which is by definition $L_x^\chi$ (see \autoref{051}). Thus we may put $\Lambda^\chi_x=(\Phi_x)^{-1}: (L_x^\chi,0)\to (\mcl_x^\chi,x)$.

Let us now verify the declared properties of $\Lambda^\chi_\blt$. Let $x\in_\mu M$. By item (i) of \autoref{076} since $n_\blt$ is locally constant $\Phi_x=\left(\Lambda_x^\chi\right)^{-1}$ is a $C^r$ diffeomorphism, whence together with ($\dagger_1$) above we have item (i) and  $\Lambda^\chi_x\in \Diff^r(L_x^\chi,0;\mcl_x^\chi,x)$. Item (ii) follows from ($\dagger_2$) above. Item (iii) follows from Currying item (iv) of \autoref{076} and using Euclidean coordinates for Lyapunov manifolds. Similarly item (iv) follows from item (vii) of \autoref{076}.

Let us next verify item (v). First let $y\in\mcl_{x,\loc}^\chi$ and $v\in L_x^\chi$ be such that $\Lambda_x^\chi(v)\in\mcl_{x,\loc}$. Then in Euclidean coordinates we have

\begin{align*}
\dfrac{d}{dv} \left(\Lambda_y^\chi\right)^{-1}\circ \Lambda_x^\chi(v)
=\dfrac{\Delta_y^\chi(\Lambda_x^\chi(v))}{\Delta_x^\chi(\Lambda_x^\chi(v))}
=\Delta_y^\chi(\Lambda_x^\chi(v)) \Delta_{\Lambda_x^\chi(v)}^\chi(x)
=\Delta_y^\chi(x).
\end{align*} 

Here in the first equality we used the fact that $y=\int_x^{\Lambda_x^\chi(y)} \Delta_x^\chi(z)\,dz$ and in the last two equalities we used item (ii) of \autoref{076}. Consequently $\left(\Lambda_y^\chi\right)^{-1}\circ \Lambda_x^\chi$ has locally constant derivative. Next let $y\in\mcl_x^\chi$ and $v\in L_x^\chi$ be arbitrary. For any appropriately chosen $s\in\RR^k$ we have that

\begin{align*}
\left|\left(\Lambda_{\alpha_s(y)}^\chi\right)^{-1} \circ \Lambda^\chi_{\alpha_{s }(x)} \circ T_x^\chi\alpha_{s }(v)\right|
&=\left|\left(\Lambda_{\alpha_s(y)}^\chi\right)^{-1} \circ \alpha_{s }\circ \Lambda^\chi_x(v)\right|\\
&=\left|T_y^\chi\alpha_s \circ \left(\Lambda_{y}^\chi\right)^{-1} \circ \Lambda^\chi_{\alpha_{s }(x)}(v)\right|.
\end{align*}

Here in the first equality we use item (ii) and in the second equality we use ($\dagger_2$). As for the appropriate choice of $s$ we need both $y$ and $\Lambda_x^\chi(v)$ to be near $x$ simultaneously, so one can take $s$ to be any time of the form $nt$ where $n\in\ZZ$ is such that

$$n \geq \max\left\{n_y,n_{\Lambda_x^\chi(v)}\right\}.$$

Finally let us verify uniqueness. Say for $x\in_\mu M$, we have a diffeomorphism  

$$K_x\in\Diff^r(L_x^\chi,0;\mcl_x^\chi,x)$$

satisfying the properties (i) to (v). We claim that $\all x\in_\mu M: K_x = \Lambda_x^\chi$. Put $R_x=\left(K_x\right)^{-1}\circ \Lambda_x^\chi:L_x^\chi\to L_x^\chi$. As in the proof of item (v) once $R_x(v)=v$ is verified locally the equation can be extended to the whole Lyapunov subspace. Let $v\in L_x^\chi$ be such that $\Lambda_x^\chi(v)\in\mcl_{x,\loc}^\chi$. Then by item (ii) we have for any $s\in\RR^k$:

$$\left| R_{\alpha_s(x)}\circ T_x^\chi\alpha_s(v) \right| = \left| T_x^\chi\alpha_s\circ R_x(v)\right|.$$

Then by the one-dimensionality of $L_x^\chi$ we have:

\begin{align*}
\left| R_x(v) \right|  
= \dfrac{\left| R_{\alpha_s(x)}\left( T_x^\chi\alpha_s(v) \right) \right|}{\left\Vert T_x^\chi \alpha_s\right\Vert }
= \dfrac{\left| R_{\alpha_s(x)}\left( T_x^\chi\alpha_s(v) \right) \right|}{\left| R_{x}\left( T_x^\chi\alpha_s(v) \right) \right|} \,\, \dfrac{\left| R_{x}\left( T_x^\chi\alpha_s(v) \right) \right|}{\left| T_x^\chi \alpha_s  (v)\right|}\,\, |v|. \quad\quad (\dagger_3)
\end{align*}

We will chose $s$ appropriately and take a limit. By \autoref{108}, we have that for $x\in_\mu M$, there is a $\Lambda=\Lambda_x\in\RRP$ large enough such that $x\in \Pi(\eps,\Lambda)$ and that this Lusin-Pesin set has positive $\mu$-measure, where $\eps$ is the same $\eps$ we fixed in the beginning of the proof. Further we may assume that $x$ is a density point of $\Pi(\eps,\Lambda)$ without hurting the \ae-quantifiers. Then by Poincar\'{e} Recurrence Theorem there is a sequence $n_k\subseteq \ZZP$ such that $\alpha_{n_kt}(x)\in \Pi(\eps,\Lambda)$ and $\lim_{k\to\infty} \alpha_{n_kt}(x)=x$. Plugging in $s=n_kt$ in ($\dagger_3$) and taking the limit $k\to \infty$ we obtain:

\begin{align*}
\left| R_x(v) \right|  
= \dfrac{\left| R_{\alpha_{n_kt}(x)}\left( T_x^\chi\alpha_{n_kt}(v) \right) \right|}{\left| R_{x}\left( T_x^\chi\alpha_{n_kt}(v) \right) \right|} \,\, \dfrac{\left| R_{x}\left( T_x^\chi\alpha_{n_kt}(v) \right) \right|}{\left| T_x^\chi \alpha_{n_kt}  (v)\right|}\,\, |v|\xrightarrow{k\to\infty} |v|.
\end{align*}

Here the first factor goes to $1$ since $L^\chi=S(\alpha_t)$ (so that the arguments converge to $0$ in Euclidean coordinates) and by item (iv) we have that $R_\blt$ is uniformly continuous on the Lusin-Pesin set $\Pi(\eps,\Lambda)$, and similarly the second factor goes to $|T_0R_x|$ which is equal to $1$ by item (i). Thus we have $|R_x(v)|=|v|$. Again item (i) implies that $R_x$ is orientation preserving, whence $R_x(v)=v$, that is, $K_x(v)=\Lambda_x^\chi(v)$ for $|v|$ small. The extension to an arbitrary $v\in L_x^\chi$ is similar to the extension argument used in the proof of item (v); we omit writing it.

\done
\end{pf}

\begin{rem}\label{075}

In the context of \autoref{017} for the uniqueness part one does not need all the properties already proved in the proof. More specifically, for uniqueness we only need that the map whose uniqueness is under scrutiny is a $C^r$ diffeomorphism $(L_x^\chi,0)\to (\mcl_x^\chi,x)$ that satisfies item (i), item (ii) only with $s=t$, item (iii) and item (iv).

Another important point is that after \autoref{078}, compared to  the derivatives defined in \autoref{076}, $\Lambda_\blt^\chi$ does not depend on $t$, nor does it depend on the section $\xi$. Thus the family $\Lambda_\blt^\chi$ is truly a global object attached to the system $(\mu,\alpha)$.

Let us also note that the affine maps established in item (v) are orientation preserving only for $y$ near $x$; for an arbitrary $y\in\mcl_x^\chi$ the compositions may reverse orientation.

\done
\end{rem}

\begin{dfn}\label{027}

In the context of \autoref{017}, $\Lambda_\blt^\chi$ is called the \textbf{nonstationary linearization} of the Lyapunov \ae-foliation $\mcl^\chi(\alpha)$ associated to the Lyapunov exponent $\chi\in\LSpec^\ast(\mu,\alpha)$.

\done
\end{dfn}

Consequently $\left\{\left.\left(\Lambda^\chi_y\right)^{-1}\,\right|\, y\in \mcl^\chi_x\right\}$ is a $C^r$ $\Aff(\RR)$-manifold structure on $\mcl^\chi_x$ w/r/t which $\mcl^\chi_x$ is a globally-affine manifold. In particular each $\mcl^\chi_x$ is $\cat{Man}^r_{\Aff(\RR)}$-isomorphic to $\RR$ with its standard affine manifold structure. By \autoref{037} we thus have $\all x\in_\mu M:\Diff^r_{\Aff(\RR)}(\mcl^\chi_x)\cong \Aff(\RR)$. Let us finally note that as for $ x\in_\mu M$ the conditional measure $\mu^\chi_x$ is absolutely continuous w/r/t the Lebesgue class $\leb_{\mcl^\chi_x}$, we have that the affine manifold structures push $\mu^\chi_\blt$ to a family of Haar measures:

\begin{lem}[\fn{\cite[p.138, Lem.3.10]{MR2261075}, \cite[p.392, Lem.7.4]{MR2811602}}]\label{125}

For any $x\in_\mu M$: $\fwd{\left(\Lambda^\chi_x\right)^{-1}}(\mu^\chi_x)$ is a Haar measure on $L^\chi_x$.

\done
\end{lem}

\subsection{Affine Structures for  Stable, Unstable, Orbit-Stable and Orbit-Unstable \ae-Foliations}

In this subsection we assemble the nonstationary linearizations for Lyapunov \ae-foliations to obtain affine structures for the stable, unstable, orbit-stable and orbit-unstable \ae-foliations of any Weyl chamber and we'll discuss the main properties of such affine structures. Note that as a consequence of \autoref{044}, we may write the stable \ae-foliation of any Weyl chamber as the intersection of codimension-$(k+1)$ stable \ae-foliations. Indeed, let $\euc\in\Cham(\mu,\alpha)$. Then

\begin{align*}
\bigcap_{(\chi,\euc)=+}S^{-\euw(\chi)}
&=\bigcap_{(\chi,\euc)=+} \bigoplus_{\rho\neq\chi} L^\rho\\
&= \bigcap_{(\chi,\euc)=+}\left[\left(\bigoplus_{\substack{\rho\neq\chi\\ (\rho,\euc)=+}} L^\rho\right) \oplus \left(\bigoplus_{\substack{\rho\neq\chi\\ (\rho,\euc)=-}} L^\rho\right) \right]\\
&=\bigoplus_{(\chi,\euc)=-} L^\chi
=S^\euc,
\end{align*}

and in particular we have

$$\all \chi\in\LSpec^\ast(\mu,\alpha): L^\chi = S^{\euw(\chi)}=\bigcap_{(\rho,\euw(\chi))=+} S^{-\euw(\rho)} = \bigcap_{\rho\neq\chi} S^{-\euw(\rho)}.$$

Now we establish affine structures for stable \ae-foliations. We only give an outline of the proof as the normal form theory is already well-developed.

\begin{cor}\label{109}

Let $\euc\in\Cham(\mu,\alpha)$ and $s(\euc)=\rank(S^\euc)\in \OL{k}$. Then for any $x\in_\mu M$ there is a unique a diffeomorphism $\Sigma^\euc_x\in$  $\Diff^r\left(S_x^\euc,0\,;\,\mcs_x^\euc,x\right)$ such that

\begin{enumerate}

\item $\all x\in_\mu M, \all \chi\in\LSpec^\ast(\mu,\alpha): (\chi,\euc)=- \implies \restr{\Sigma^\euc_x}{L^\chi_x}\in \Diff^r(L^\chi_x,0;\mcl^\chi_x,x)$,

\item $\phi^{\Sigma^\euc}:\RR^k\times M\to \DGL\left(\RR^{s(\euc)}\right)$, $(t,x)\mapsto \left(\Sigma_{\alpha_t(x)}^\euc\right)^{-1}\circ \alpha_t\circ \Sigma_x^\euc$ defines a measurable cocycle over $\alpha$ taking values in the group of diagonal linear automorphisms of $\RR^{s(\euc)}$ w/r/t the Euclidean coordinates on global stable manifolds,

\item W/r/t Euclidean coordinates on global manifolds $\Sigma^\euc_\blt: M\to \Imm^r(\RR^{s(\euc)};M)$ is measurable,

\item On any Lusin-Pesin set $\Sigma^\euc_\blt$ is uniformly continuous. More precisely, for any pre-Lusin-Pesin set $\Pi$ and for any $\Lambda\in\RRP$, w/r/t Euclidean coordinates on global manifolds we have that $\Sigma^\euc_\blt:\Pi(\eps,\Lambda)\to \Imm^r(\RR^{s(\euc)};M)$ is uniformly continuous, where the target is endowed with the $C^r$ compact-open topology,

\item $\all x\in_\mu M, \all y\in\mcs_x^\euc :\left(\Sigma_y^\euc\right)^{-1}\circ\Sigma_x^\euc:S_x^\euc\to S_y^\euc$ is affine with diagonal linear part.

\end{enumerate}

Moreover analogous statements establish a family $\Upsilon^\euc_\blt$ for the unstable \ae-foliation of $\euc$.

\done
\end{cor}

\begin{pf}

Let $\chi^1,...,\chi^s$ be the only Lyapunov exponents that are negative on $\euc$. Note that by the rank-nullity theorem $\dim(\ker(\chi^1-\chi^2))=k-1$ and by induction $\dim(\{t\in \RR^k\,|\, \chi^1(t)=\chi^2(t)=\cdots=\chi^s(t)\})=k-s+1>0$ is a positive dimensional linear space. Consequently there is a $t\in\euc$ such that $\chi^1(t)=\chi^2(t)=\cdots = \chi^s(t)$. Then the normal form theory presented in \cite[App. A]{MR3503686} applies to the diffeomorphism $\alpha_t$ and guarantees all statements but uniqueness.  Alternatively one can apply \cite[pp. 344-345, Thm.2.3]{MR3642250} to $\alpha_t$ which guarantees uniqueness as well. Let us also note that initially the arguments in \cite{MR3503686} gives item (v) for $y\in_{\aem}\mcs^\euc_x$; here we let the first $x\in_\mu M$ absorb the cases item (v) does not hold for every $y\in\mcs^\euc_x$; this is valid since the conditionals $\mu^\chi_\blt$ are absolutely continuous w/r/t the Lebesgue classes $\leb_{\mcl^\chi_\blt}$.

\done
\end{pf}

\begin{rem}\label{115}

Le us also note that by comparing \autoref{109} with \autoref{017}, for any $\euc\in\Cham(\mu,\alpha)$ and for any $\chi\in\LSpec^\ast(\mu,\alpha)$ with $(\chi,\euc)=-$, we have $\restr{S^\euc_x}{L^\chi_x} = \Lambda^\chi_x$. Further, we have:

\begin{align*}
\all x\in_\mu M, \all t\in\RR^k, \all v\in S^\euc_x: &\phi^{\Sigma^\euc}(t,x)(v)=\pm \restr{T_x\alpha_t}{S^\euc_x}(v),\\
\all x\in_\mu M, \all t\in\RR^k, \all v\in U^\euc_x: &\phi^{\Upsilon^\euc}(t,x)(v)=\pm \restr{T_x\alpha_t}{U^\euc_x}(v).
\end{align*}

\done
\end{rem}

From \autoref{109} we also have the following assembly properties; they will come in very handy in establishing the independence of many maps of the choice of the chamber that is used to define them:

\begin{cor}\label{099}

\begin{enumerate}

\item $\all \euc,\eud\in\Cham(\mu,\alpha): \mcs^\euc\subseteq_{\mu}\mcs^\eud \implies \restr{\Sigma^\eud}{S^\euc}=_\mu\Sigma^\euc$.

\item $\all \euc,\eud,\eue\in\Cham(\mu,\alpha): S^\euc =_{\mu} S^\eud\oplus S^\eue \implies \Sigma^\euc=_\mu\Sigma^\eud \times \Sigma^\eue$.

\item $\all \euc,\eud,\eue\in\Cham(\mu,\alpha): S^\euc =_{\mu} S^\eud\cap S^\eue \implies \all x,y\in_\mu M: z\in \mcs^\eud_x\cap \mcs^\eue_y \implies \mcs^\euc_z\subseteq \mcs^\eud_x\cap \mcs^\eue_y $ and $ \Sigma^\euc_z = \restr{\Sigma^\eud_x}{S^\euc_z} = \restr{\Sigma^\eue_y}{S^\euc_z}$.

\end{enumerate}

 Moreover there are analogous assembly properties of unstable \ae-foliations $\mcu^\euc$ w/r/t the family $\Upsilon^\euc_\blt$.

\done
\end{cor}

As in the case of Lyapunov manifolds, we have that for any chamber $\euc\in\Cham(\mu,\alpha)$, and for $x\in_\mu M$,  $\left\{\left.\left(\Sigma^\euc_y\right)^{-1}\,\right|\, y\in \mcs^\euc_x\right\}$ is a $C^r$ $\DAff(\RR^{s(\euc)})$-manifold structure on $\mcs^\euc_x$ w/r/t which $\mcs^\euc_x$ is a globally-$\DAff(\RR^{s(\euc)})$-affine manifold.  The assembly properties \autoref{099} guarantee that these affine manifold structures on different stable manifolds are compatible. There are  similar statements for unstable manifolds of chambers.

Finally we extend the affine manifold structures defined on the stable and unstable manifolds to orbit-stable and orbit-unstable manifolds. Recall that by \autoref{114} $(\mu,\alpha)$ is essentially free and hence each orbit carries a canonical affine manifold structure isomorphic to the Euclidean structure of $\RR^k$.

\begin{prp}\label{111}

Let $\euc\in\Cham(\mu,\alpha)$ and $s(\euc)=\rank(S^\euc)\in \OL{k}$. Then for any $x\in_\mu M$  there is a unique a diffeomorphism $\Gamma\Sigma^\euc_x\in$  $\Diff^r\left(\RR^k\times S_x^\euc,0\,;\,\mco\mcs_x^\euc,x\right)$ such that

\begin{enumerate}

\item $\all x\in_\mu M: \restr{\Gamma\Sigma^\euc_x}{\RR^k}=\alpha_\blt(x)\in \Diff^r\left(\RR^k;\mco_x\right)$,

\item $\all x\in_\mu M: \restr{\Gamma\Sigma^\euc_x}{ S_x^\euc}=\Sigma^\euc_x\in \Diff^r(S^\euc_x;\mcs^\euc_x)$,

\item $\phi^{\Gamma\Sigma^\euc}:\RR^k\times M\to (\RR^k\times \RR^{s(\euc)})\rtimes \left(\{I_k\}\times \DGL\left(\RR^{s(\euc)}\right)\right)\leq \DAff(\RR^k\times \RR^{s(\euc)})$, $(t,x)\mapsto \left(\Gamma\Sigma_{\alpha_t(x)}^\euc\right)^{-1}\circ \alpha_t\circ \Gamma\Sigma_x^\euc$ defines a measurable cocycle over $\alpha$,

\item W/r/t Euclidean coordinates on orbits and on global invariant manifolds $\Gamma\Sigma^\euc_\blt: M\to \Imm^r(\RR^k\times\RR^{s(\euc)};M)$ is measurable,

\item On any Lusin-Pesin set $\Gamma\Sigma^\euc_\blt$ is uniformly continuous. More precisely, for any pre-Lusin-Pesin set $\Pi$ and for any $\Lambda\in\RRP$, w/r/t Euclidean coordinates on orbits and global manifolds we have that $\Gamma\Sigma^\euc_\blt:\Pi(\eps,\Lambda)\to \Imm^r(\RR^k\times\RR^{s(\euc)};M)$ is uniformly continuous, where the target is endowed with the $C^r$ compact-open topology,

\item $\all x\in_\mu M, \all y\in\mcs_x^\euc :\left(\Gamma\Sigma_y^\euc\right)^{-1}\circ\Gamma\Sigma_x^\euc:\RR^k\times S_x^\euc\to \RR^k\times S_y^\euc$ is affine with diagonal linear part and the linear part of the first component $I_k$.

\end{enumerate}

Moreover analogous statements establish a family $\Gamma\Upsilon^\euc_\blt$ for the orbit-unstable \ae-foliation of $\euc$.

\done
\end{prp}

\begin{pf}

By definition we have $\mco\mcs_x^\euc = \{\alpha_t(y)\,|\, t\in \RR^k, y\in\mcs_x^\euc\}$, and since we have already established that $\mcs_x^\euc$ is affinely isomorphic to $\RR^{s(\euc)}$ this is a bijective description. Thus

$$\Gamma\Sigma_x^\euc:(t,v)\mapsto \alpha_t\circ \Sigma^\euc_x(v)$$

is a diffeomorphism from $\RR^k\times \mcs_x^\euc$ onto $\mco\mcs_x^\euc$. The first two items are syntactic and the remaining items as well as uniqueness follow from \autoref{017} and \autoref{075}.

\done
\end{pf}

\autoref{125} together with the assembly properties \autoref{099} gives:

\begin{lem}\label{129}

Let $\euc\in\Cham(\mu,\alpha)$. Then $\all x\in_\mu M$: $\fwd{\left(\Sigma^\euc_x\right)^{-1}}\left(\mu^{S^\euc}_x\right)$, $\fwd{\left(\Upsilon^\euc_x\right)^{-1}}\left(\mu^{U^\euc}_x\right)$ are Haar measures on $S^\euc_x$ and $U^\euc_x$.

\done
\end{lem}

Let us observe in closing that by the second item in \autoref{111} assembly properties analogous to those in \autoref{099} are available after orbit saturation.

%% file: uzman_arithmeticity_affineholonomies.tex
\section[Affine Holonomies along Invariant \ae-Foliations]{Affine Holonomies along Invariant \ae-Foliations\label{068}}

We had defined in \autoref{120} holonomies along an \ae-foliation. In this section we establish the existence of global holonomies along the stable and unstable \ae-foliations of any chamber using the affine structures established in \autoref{112}. Here the higher and maximal rank assumptions play a key role; note that in the rank-one case holonomies are available locally and on positive measure subsets of transversals. We start with an observation.

\begin{obs}\label{113}

Let $Y$ be a $C^s$ complete globally affine manifold, $\mcf$ be a $C^s$ foliation of $Y$ with each leaf a $C^s$ injectively immersed affine submanifold of $Y$, and let $L,R$ be $C^s$ injectively immersed complete affine submanifolds of $Y$ transverse to $\mcf$. Then there are $C^s$ global affine charts on $Y$ such that w/r/t these charts both $L$ and $R$ are horizontal and each leaf of $\mcf$ is vertical. Furthermore substituting the first (that is, horizontal) component w/r/t these charts defines the holonomy $\mcf_{R\leftarrow L}: L\to R$ from $L$ to  $R$ along the leaves of $\mcf$ that is a $C^s$ affine diffeomorphism.

\done
\end{obs}

We construct our global holonomies piece by piece using the above observation with appropriately chosen invariant manifolds. We start with the simpler case of a chamber of the form $\euc=\euw(\chi)$ (see \autoref{044}).

\begin{prp}\label{116}

Let $\chi\in\LSpec^\ast(\mu,\alpha)$ and put $\euc=\euw(\chi)$. Then $\all x\in_\mu M, \all y\in_{\aem} \mcu^\euc_x$, there is a holonomy $\mcu_{y\leftarrow x}^\euc: \mco\mcs_x^\euc\to \mco\mcs_y^\euc$ along the unstable \ae-foliation $\mcu^\euc$ of $\euc$ that is a $C^r$ affine diffeomorphism with the property that

$$\all t\in\RR^k,\all x\in_\mu M, \all y\in_{\aem} \mcu^\euc_x:\restr{\mcu^\euc_{\alpha_t(y)\ot \alpha_t(x)}}{\mcs^\euc_{\alpha_t(x)}}\circ \alpha_t = \alpha_t\circ \restr{\mcu^\euc_{y\ot x}}{\mcs^\euc_x}.$$

Consequently, w/r/t affine coordinates the linear part of the first component of $\mcu^\euc_{y\leftarrow x}$ is $I_k$.

\done
\end{prp}

\begin{pf}

First we define holonomies along the unstable \ae-foliation between global stable manifolds and then we extend the definition to be between global orbit-stable manifolds. Note that as $\euc=\euw(\chi)$ and there are exactly $k+1$ distinct Lyapunov exponents, there are at least two Lyapunov exponents of $(\mu,\alpha)$ distinct from $\chi$, so that there are two chambers $\eud,\eue\in\Cham(\mu,\alpha)$ such that $U^\euc=_\mu U^\eud\oplus U^\eue =_\mu S^{-\euc}$; thus we also have

\begin{align*}
S^\euc &=_\mu S^\eud \cap S^\eue \\
S^\eue &=_\mu S^\euc \oplus U^\eud\quad\quad (\dagger)\\
S^\eud &=_\mu S^\euc \oplus U^\eue
\end{align*}

By Pesin theory, there is an $M_0=_\mu M$ such that $\all x\in M_0$, $\all y\in M_0\cap \mcu_x^\euc$, there are $a,b\in \Osel(\alpha)\cap\mcu^\euc_x$ such that

$$a\in\mcu^\eud_x, b\in \mcu^\eue_a, y\in \mcu^\eud_b.$$

By the absolute continuity of $\mu$ we may rearrange the quantifiers for $x$ and $y$ to match the statement of the proposition. Note that indeed only one pivot point from $x$ to $y$ in $\mcu^\euc_x$ may fail to be sufficient, but two pivot points are sufficient. We shall define the desired holonomy as the composition of three holonomies with the required properties. Since by the second equation in $(\dagger)$ we have $S^\euc_x\oplus U^\eud_x = U^{-\eue}_x$, by \autoref{109} and \autoref{099}, following \autoref{113} we may define a $C^r$ affine diffeomorphism $\mcu^{\eud}_{a\leftarrow x}: \mcs^{\euc}_x\to\mcs^{\euc}_a$ by

$$\mcu^{\eud}_{a\leftarrow x} = \Upsilon^{-\eue}_x\circ \begin{pmatrix}\text{component} \\\text{substitution}\end{pmatrix}\circ \left(\restr{\Upsilon^{-\eue}_x}{S^{\euc}_x}\right)^{-1}.$$

In particular for any $z\in \mcs^{\euc}_x$ we have $\mcu^{\eud}_{a\leftarrow x}(z) $ $\in \mcu^{\eud}_z\cap \mcs^{\euc}_a\subseteq \mcu^{\eud}_z\cap \mcs^{\eud}_a$. Note that Pesin theory provides this only locally and for a positive $\leb_{\mcs^\euc_x}$-measure set; the holonomy above is everywhere defined on the manifold $\mcs^\euc_x$ by virtue of it being defined in terms of affine coordinates; further whenever a holonomy is provided by Pesin theory its value coincides with the value of the holonomy above. Similarly, by the third and second equations in $(\dagger)$, respectively, we may define $C^r$ affine diffeomorphisms $\mcu^{\eue}_{b\leftarrow a}: \mcs^{\euc}_a\to\mcs^{\euc}_b$ and $\mcu^{\eud}_{y\leftarrow b}: \mcs^{\euc}_b\to\mcs^{\euc}_y$ by

\begin{align*}
\mcu^{\eue}_{b\leftarrow a} &= \Upsilon^{-\eud}_a\circ \begin{pmatrix}\text{component} \\\text{substitution}\end{pmatrix}\circ \left(\restr{\Upsilon^{-\eud}_a}{S^{\euc}_a}\right)^{-1}, \text{ and }\\
\mcu^{\eud}_{y\leftarrow b} &= \Upsilon^{-\eue}_b\circ \begin{pmatrix}\text{component} \\\text{substitution}\end{pmatrix}\circ \left(\restr{\Upsilon^{-\eue}_b}{S^{\euc}_b}\right)^{-1},
\end{align*}

respectively. Thus defining

$$\mcu_{y\leftarrow x}^\euc= \mcu^{\eud}_{y\leftarrow b}\circ \mcu^{\eue}_{b\leftarrow a} \circ \mcu^{\eud}_{a\leftarrow x}$$

gives us a $C^r$ affine diffeomorphism from $\mcs^\euc_x$ onto $\mcs^\euc_y$. Note that by construction for any $z\in\mcs^\euc_x$ we have $\mcu^\euc_{y\leftarrow x}(z)\in \mcu_z^\euc\cap \mcs_y^\euc$. Next let us verify the $\alpha$-invariance. Our earlier choices of pivot points now translate into

$$\alpha_t(a)\in \mcu^\eud_{\alpha_{t}(x)}, \alpha_t(b)\in \mcu^\eue_{\alpha_t(a)}, \alpha_t(y)\in \mcu^\eud_{\alpha_t(b)}.$$

The main property that facilitates $\alpha$-invariance is the fact that the cocycles $\phi^{\blt}$ attached to the stable and unstable affine parameters in item (ii) of \autoref{109} have no translational parts. Indeed we have 

$$\Upsilon^{-\eue}_{\alpha_t(x)}=\alpha_t\circ \Upsilon^{-\eue}_{x}\circ \left(\phi^{\Upsilon^{-\eue}}(t,x)\right)^{-1},$$

and consequently

\begin{align*}
\mcu^{\eud}_{\alpha_t(a)\ot \alpha_t(x)}
&= \Upsilon^{-\eue}_{\alpha_t(x)}\circ  \begin{pmatrix}\text{component} \\\text{substitution}\end{pmatrix}\circ \left(\restr{\Upsilon^{-\eue}_{\alpha_t(x)}}{S^{\euc}_{\alpha_t(x)}}\right)^{-1}\\
&=\alpha_t\circ \Upsilon^{-\eue}_x\\
&\phantom{=\alpha_t\circ}\circ \left[\left(\phi^{\Upsilon^{-\eue}}(t,x)\right)^{-1}\circ  \begin{pmatrix}\text{component} \\\text{substitution}\end{pmatrix}\circ \phi^{\Upsilon^{-\eue}}(t,x)\right]\\
&\phantom{=\alpha_t\circ=\alpha_t\circ}\circ \left(\restr{\Upsilon^{-\eue}_x}{S^{\euc}_x}\right)^{-1} \circ\restr{\alpha_{-t}}{\mcs^\euc_{\alpha_t(x)}}\\
&=\alpha_t\circ\mcu^\eud_{a\ot x}\circ\restr{\alpha_{-t}}{\mcs^\euc_{\alpha_t(x)}}.
\end{align*}

Applying this to the other components that make up $\mcu^{\eud}_{\alpha_t(y)\ot \alpha_t(x)}$ we obtain the invariance property. Finally we may extend

$$\mcu^{\euc}_{y\leftarrow x} = \Gamma\Sigma^\euc_y\circ \left(I_k\times \left(\left(\Sigma^\euc_y\right)^{-1}\circ \mcu^\euc_{y\leftarrow x}\circ \Sigma^\euc_x \right)\right)\circ \left(\Gamma\Sigma^\euc_x\right)^{-1}: \mco\mcs^\euc_x\to\mco\mcs^\euc_y.$$

By construction this is a $C^r$ affine diffeomorphism such that w/r/t affine coordinates 	the linear part of the first component is $I_k$. Let us also verify that this is a holonomy between the orbit-stable manifolds of $\euc$ along the unstable \ae-foliation of $\euc$. After the proof of \autoref{111}, let $z\in\mco\mcs^\euc_x$ and let $(t,v)\in\RR^k\times S^\euc_x$ be the unique affine coordinates so that $\Gamma\Sigma^\euc_x(t,v)=z$. Then 

$$\mcu^\euc_{y\leftarrow x}(z)
 =  \Gamma\Sigma^\euc_y \left(t,\left(\Sigma^\euc_y\right)^{-1}\circ \mcu^\euc_{y\leftarrow x}\circ \alpha_{-t}(z)\right)
=\alpha_t\circ \underbrace{\mcu^\euc_{y\leftarrow x}\circ \alpha_{-t}(z)}_{\in \,\mcu^\euc_{\alpha_{-t}(z)}\cap \mcs^\euc_y}\in \mcu_z^\euc\cap \mco\mcs^\euc_y. $$

\done
\end{pf}

\begin{rem}\label{119}

The affine holonomies we construct in \autoref{116} are independent of the splitting chambers $\eud$ and $\eue$ by the assembly properties \autoref{099}. They are also independent of the pivot points $a,b$, however we will establish this independence in \autoref{121} only after establishing the coherence of affine holonomies along stable and unstable \ae-foliations in \autoref{118} below. Until then we will need to be careful about the choices of these pivot points.

\done
\end{rem}

Next we establish global holonomies that are $C^r$ affine diffeomorphisms for anonymous chambers. For this we will use induction. Since at this point the independence of the affine holonomy on the pivot points is not established we will update these pivot points in such a way that all affine holonomies that need to be predetermined before a step are determined appropriately.

\begin{prp}\label{117}

Let $\euc\in\Cham(\mu,\alpha)$. Then $\all x\in_\mu M, \all y\in_{\aem} \mcu^\euc_x$, there is a holonomy $\mcu_{y\leftarrow x}^\euc: \mco\mcs_x^\euc\to \mco\mcs_y^\euc$ along the unstable \ae-foliation $\mcu^\euc$ of $\euc$ that is a $C^r$ $\DAff(\RR^k\times \RR^{s(\euc)})$-affine diffeomorphism with the property that

$$\all t\in\RR^k,\all x\in_\mu M, \all y\in_{\aem} \mcu^\euc_x:\restr{\mcu^\euc_{\alpha_t(y)\ot \alpha_t(x)}}{\mcs^\euc_{\alpha_t(x)}}\circ \alpha_t = \alpha_t\circ \restr{\mcu^\euc_{y\ot x}}{\mcs^\euc_x}.$$

Consequently, w/r/t affine coordinates the linear part of the first component of $\mcu^\euc_{y\leftarrow x}$ is $I_k$. Similarly there is a family $\mcs_{\blt\leftarrow\blt}^\euc$ of holonomies along the stable \ae-foliation $\mcs^\euc$ of $\euc$ that are $C^r$ $\DAff(\RR^k\times \RR^{u(\euc)})$-affine diffeomorphisms.

\done
\end{prp}

\begin{pf}

It suffices to establish only the affine holonomies along the unstable \ae-foliations. We shall use induction on the rank $s(\euc)$ of the stable \ae-folation of the chamber $\euc$. \autoref{116} covers the case $s(\euc)=1$. Suppose $s(\euc)\geq 2$, and that for any lower rank the statement is true. Then there are two chambers $\eud,\eue\in\Cham(\mu,\alpha)$ such that $S^\euc=_\mu S^\eud\oplus S^\eue=U^{-\euc}$; consequently we also have

\begin{align*}
U^\euc &=_\mu U^\eud \cap U^\eue \\
U^\eue &=_\mu U^\euc \oplus S^\eud \quad\quad (\dagger)\\
U^\eud &=_\mu U^\euc\oplus S^\eue  
\end{align*}

Let $x\in_\mu M$ and $y\in_{\aem} \mcu^\euc_x$. On the one hand, by the second equation in $(\dagger)$, $\restr{\mcu^\euc_{y\leftarrow x}}{\mcs^\eud_x}$ is an affine holonomy given by component substitution w/r/t the $\Upsilon^\eue_x$ coordinates with image $\mcs^\eud_y$. On the other hand, by induction hypothesis, since $s(\eud)<s(\euc)$, there is an affine holonomy $\mcu^\eud_{y\leftarrow x}$ that a priori depends on the choice of pivot points. By the first equation in $(\dagger)$, we may choose the pivot points involved in $\mcu^\eud_{y\leftarrow x}$ so that

$$\restr{\mcu^\euc_{y\leftarrow x}}{\mcs^\eud_x} = \mcu^\eud_{y\leftarrow x}:\mcs^\eud_x\to\mcs^\eud_y.$$

Let $z\in_{\aem} \mcs_x^\euc$. Then since $\mcs^\eud$ and $\mcs^\eue$ transversely sub-\ae-foliate $\mcs^\euc_x$, there is a unique $a=a_z\in \mcs_x^\eud\cap \mcs_z^\eue$. Put $b=b_z=\restr{\mcu^\euc_{y\leftarrow x}}{\mcs^\eud_x}(a)$. On the one hand $\restr{\mcu^\euc_{b\leftarrow a}}{\mcs^\eue_a}$ is an affine holonomy given by component substitution w/r/t the $\Upsilon^\eud_a$ coordinates with image $\mcs_b^\eue$ by the third equation in $(\dagger)$. On the other hand, by induction hypothesis there is an affine holonomy $\mcu^\eue_{b\leftarrow a}$. Then by the first equation in $(\dagger)$, we may choose the pivot points involved in the affine holonomy provided by the induction hypothesis so that

$$\restr{\mcu^\euc_{b\leftarrow a}}{\mcs^\eue_a} = \mcu^\eue_{b\leftarrow a}:\mcs^\eue_a\to\mcs^\eue_b. $$

Applying this holonomy to $z\in\mcs^\eue_a$, we have that $\restr{\mcu^\euc_{b\leftarrow a}}{\mcs^\eue_a}(z) = \mcu^\eue_{b\leftarrow a}(z)$ $\in \mcu^\euc_z\cap \mcs^\eue_b \subseteq  \mcu^\eue_z\cap \mcs^\eue_b$; we may thus define $\mcu^\eue_{y\leftarrow x}(z)$ to be this point. Finally  we may extend $\mcu^{\euc}_{y\leftarrow x}$ to an affine holonomy between orbit-stables using the formula we used for the same purpose in the proof of \autoref{117}:

$$\mcu^{\euc}_{y\leftarrow x} = \Gamma\Sigma^\euc_y\circ \left(I_k\times \left(\left(\Sigma^\euc_y\right)^{-1}\circ \mcu^\euc_{y\leftarrow x}\circ \Sigma^\euc_x \right)\right)\circ \left(\Gamma\Sigma^\euc_x\right)^{-1}: \mco\mcs^\euc_x\to\mco\mcs^\euc_y.$$

The verification that this is indeed a holonomy along $\mcu^\euc$ as well as the $\alpha$-invariance property is straightforward and is as in the previous proof.

\done
\end{pf}

\begin{rem}\label{143}

Note that we stated the $\alpha$-invariance properties of the affine holonomies in \autoref{116} and \autoref{117} for the restrictions of these holonomies to stable or unstable manifolds. That there is invariance in exactly this sense allows us to define affine holonomies between orbit-stable manifolds and orbit-unstable manifolds. Moreover when affine holonomies are not restricted to stable or unstable manifolds, the points $x$ and $\alpha_t(x)$ determine the same orbit-stable and orbit-unstable manifolds, so we may indeed write $\mcu^\euc_{\alpha_t(y)\ot\alpha_t(x)}=\mcu^\euc_{y\ot x}$.

\done
\end{rem}

\begin{prp}\label{146}

Let $\euc\in\Cham(\mu,\alpha)$. Then the affine holonomies along the unstable manifolds of $\euc$ preserve the conditional measures of $\mu$ along the orbit-stable \ae-foliation of $\euc$. More precisely, there is a measurable $c_{\blt\to\blt}:\{(x,y)\in M\times M| x\in_\mu M, y\in_{\aem}\mcu^\euc_x\}\to \RRP$ such that

$$\all x\in_\mu M,\all y\in_{\aem}\mcu^\euc_x,\fwd{\mcu^\euc_{y\ot x}}(\mu^{\mco\mcs^\euc}_x)= c_{x\to y}\, \mu^{\mco\mcs^\euc}_y.$$

There is a similar statement for the affine holonomies $\mcs^\euc_{\blt\ot\blt}$.

\done
\end{prp}

\begin{pf}

By \autoref{129} the affine structures on stable and unstable manifolds transform the appropriate conditionals of $\mu$ to Haar measures, and affine holonomies are defined compositions of affine structures and component substitution. The proportionality function $c_{\blt\to\blt}$ is due to the fact that there is no natural way to normalize the stable and unstable manifolds for a general chamber (as opposed to a chamber of the form $\euc=\mcw(\chi)$; see \autoref{017}); see also \autoref{091} and \autoref{092}.

\done
\end{pf}

Let us next qualify the dependence of the families $\mcs^\euc_{\blt\leftarrow\blt}$ and $\mcu^\euc_{\blt\leftarrow\blt}$ on the base points:

\begin{prp}\label{127}

Let $\euc\in\Cham(\mu,\alpha)$. Then

\begin{enumerate}

\item $\mcu^\euc_{\blt\leftarrow\blt}$ depends measurably on the basepoints. More precisely, w/r/t the Euclidean coordinates on orbit-stable manifolds of $\euc$, $\mcu^\euc_{\blt\leftarrow\blt}: \{(x,y)\in M\times M\,|\, x\in_\mu M, y\in_{\aem} \mcu^\euc_x\}\to \DAff(\RR^k\times \RR^{s(\euc)})$ is measurable.

\item On any Lusin-Pesin set $\mcu^\euc_{\blt\leftarrow\blt}$ is uniformly continuous. More precisely, w/r/t the Euclidean coordinates on orbit-stable manifolds of $\euc$, for any pre-Lusin-Pesin set $\Pi$, any $\Lambda\in\RRP$, any $r_0\in\RRP$ that bounds from below the sizes of local unstable manifolds of $\euc$, $\mcu^\euc_{\blt\leftarrow\blt}:\{(x,y)\in M\times M\,|\, x\in \Pi(\eps,\Lambda), y\in\mcu^\euc_{x,\loc}[x|\leq r_0]\}\to \DAff(\RR^k\times \RR^{s(\euc)})$ is uniformly continuous.

\end{enumerate}

There are similar statements for the family $\mcs^\euc_{\blt\leftarrow\blt}$.

\done
\end{prp}

\begin{pf}

Note that any function involved in the definition of $\mcu^\euc_{\blt\leftarrow\blt}$ is either a component substitution, which is $C^r$ since the affine manifold in question is $C^r$, or an affine parameter. Thus the statement follows from the third and fourth items in \autoref{109}.

\done
\end{pf}

Finally we establish the coherence of the families $\mcs^\euc_{\blt\leftarrow\blt}$ and $\mcu^\euc_{\blt\leftarrow\blt}$ of affine holonomies along the stable and unstable \ae-foliations for an anonymous chamber $\euc$. If the orbit saturations are not involved there are no obstructions to coherence, whereas in the presence of orbit saturations the only obstruction to coherence is matching of the temporal components:

\begin{prp}\label{118}

Let $\euc\in\Cham(\mu,\alpha)$. Then 

\begin{enumerate}

\item $\all x\in_\mu M,$  $\all (y,z)\in_{\aem} \mcs_x^\euc\times \mcu_x^\euc:$ $\mcu_{z\ot x}^\euc(y)=\mcs_{y\ot x}^\euc(z)$.

\item More generally, $\all x\in_\mu M,$ $\all (y,z)\in_{\aem} \mco\mcs_x^\euc\times \mco\mcu_x^\euc:$  

\begin{align*}
&\proj_{\RR^k}\circ \left(\Gamma\Sigma^\euc_x\right)^{-1}(y)
= \proj_{\RR^k}\circ \left(\Gamma\Upsilon^\euc_x\right)^{-1}(z)\\
&\iff \mcu_{z\ot x}^\euc(y)=\mcs_{y\ot x}^\euc(z).
\end{align*}

\end{enumerate}

\done
\end{prp}

\begin{pf}

Let us first prove the first item. Let us write $y'=\mcu^\euc_{z\leftarrow x}(y)$ for the sake of brevity; our aim is to realize $y'$ as the value of the affine holonomy $\mcs^\euc_{y\leftarrow x}$ (with appropriately chosen pivot points) at the point $z\in \mcu^\euc_x$. Since $s(\euc)+u(\euc)=k+1\geq3$, by considering $-\euc$ if necessary, we may assume that $s(\euc)\geq2$. Thus there are two chambers $\eud,\eue\in\Cham(\mu,\alpha)$ such that $S^\euc =_\mu S^\eud\oplus S^\eue$, and consequently also

\begin{align*}
U^\euc &=_\mu U^\eud \cap U^\eue \\
U^\eue &=_\mu U^\euc \oplus S^\eud \quad\quad (\dagger)\\
U^\eud &=_\mu U^\euc\oplus S^\eue  
\end{align*}

By Pesin theory, there are $a,b\in \Osel(\alpha)\cap \mcs^\euc_x$ such that

\begin{align*}
&\phantom{a' =}a\in\mcs^\eud_x, b\in\mcs^\eue_a, y\in\mcs^\eud_b,\\
&a' = \mcu^\euc_{z\leftarrow x}(a)\in \Osel(\alpha)\cap \mcs^\euc_z,\\
&b' = \mcu^\euc_{a'\leftarrow a}(b)\in \Osel(\alpha)\cap \mcs^\euc_{z}.
\end{align*}

By the second equation in $(\dagger)$, w/r/t $\Upsilon^\eue_x$ coordinates, after \autoref{113}, $\restr{\mcu^\euc_{z\leftarrow x}}{\mcs^\eud_x}$ and $\restr{\mcs^\eud_{a\leftarrow x}}{\mcu^\euc_x}$ are given by substitution of complementary components (that is, if the former is given by substitution of the horizontal component then the latter is given by substitution of the vertical component). Consequently the pivot points for $\mcu^\eud_{z\leftarrow x}$ and $\mcs^\euc_{a\leftarrow x}$ can be chosen so that

\begin{align*}
a' 
&= \restr{\mcu^\euc_{z\leftarrow x}}{\mcs^\eud_x}(a) = \mcu^\eud_{z\leftarrow x}(a)\\
&= \restr{\mcs^\eud_{a\leftarrow x}}{\mcu^\euc_x}(z) = \mcs^\euc_{a\leftarrow x}(z)
\in\mcs^\eud_z.
\end{align*}

Next, by the third equation in $(\dagger)$, w/r/t $\Upsilon^\eud_a$ coordinates, $\restr{\mcu^\euc_{a'\leftarrow a}}{\mcu^\euc_a}$ and $\restr{\mcs^\eue_{b\leftarrow a}}{\mcu^\euc_a}$ are given by substitution of complementary components. Consequently the pivot points for $\mcu^\eue_{a'\leftarrow a}$ and $\mcs^\euc_{b\leftarrow a}$ can be chosen so that

\begin{align*}
b' 
&= \restr{\mcu^\euc_{a'\leftarrow a}}{\mcs^\eue_a}(b) = \mcu^\eue_{a'\leftarrow a}(b)\\
&= \restr{\mcs^\eue_{b\leftarrow a}}{\mcu^\euc_a}(a') = \mcs^\euc_{b\leftarrow a}(a')
\in\mcs^\eue_{a'}.
\end{align*}

Finally again by the second equation in $(\dagger)$, w/r/t $\Upsilon^\eue_b$ coordinates, $\restr{\mcu^\euc_{b'\leftarrow b}}{\mcs^\eud_b}$ and $\restr{\mcs^\eud_{y\leftarrow b}}{\mcu^\euc_b}$ are given by substitution of complementary components and consequently the pivot points for $\mcu^\eud_{b'\leftarrow b}$ and $\mcs^\euc_{y\leftarrow b}$ can be so chosen that

\begin{align*}
 &\restr{\mcu^\euc_{b'\leftarrow b}}{\mcs^\eud_b}(y) = \mcu^\eud_{b'\leftarrow b}(y)\\
= &\restr{\mcs^\eud_{y\leftarrow b}}{\mcu^\euc_b}(b') = \mcs^\euc_{y\leftarrow b}(b')
\in\mcs^\eud_{b'}.
\end{align*}

Since $y\in\mcs^\eud_b$, this last point also coincides with $y'=\mcu^\euc_{z\leftarrow x}(y)$ due to the first description of the point. Therefore

\begin{align*}
y'
&=\mcu^\euc_{z\leftarrow x}(y)
= \mcs^\euc_{y\leftarrow b}(b')
=\mcs^\euc_{y\leftarrow b}\circ \mcs^\euc_{b\leftarrow a}(a')\\
&=\mcs^\euc_{y\leftarrow b}\circ \mcs^\euc_{b\leftarrow a}\circ \mcs^\euc_{a\leftarrow x}(z)
= \mcs^\euc_{y\leftarrow x}(z),
\end{align*}

where in the last equality we have used the fact that all the pivot points for the affine holonomies that are being composed are compatible.

For the second item, for sufficiency we may simply apply the argument for the first item with $\alpha_t(x)$ instead of $x$, where $t=\proj_{\RR^k}\circ \left(\Gamma\Sigma^\euc_x\right)^{-1}(y) = \proj_{\RR^k}\circ \left(\Gamma\Upsilon^\euc_x\right)^{-1}(z)$. Necessity of this condition is straightforward.

\done
\end{pf}

\begin{cor}\label{121}

In \autoref{116}, and consequently in \autoref{117} and \autoref{118}, the $C^r$ affine holonomies are independent of the pivot points.

\done
\end{cor}

\begin{pf}

Let $A,B$ be points alternative to $a,b$ in the proof of \autoref{116}, respectively. Independence of affine holonomies on the pivot points means that

$$ \mcu^{\eud}_{y\leftarrow b}\circ \mcu^{\eue}_{b\leftarrow a} \circ \mcu^{\eud}_{a\leftarrow x}  = \mcu^{\eud}_{y\leftarrow B}\circ \mcu^{\eue}_{B\leftarrow A} \circ \mcu^{\eud}_{A\leftarrow x}, \quad (\star)$$

where each one of the composed holonomies is given by component substitution w/r/t some affine coordinate. By factoring holonomies

$$\mcu^\eud_{A\leftarrow x} = \mcu^\eud_{A\leftarrow a}\circ \mcu^\eud_{a\leftarrow x}, \quad\quad \mcu^\eud_{y\leftarrow b} = \mcu^\eud_{y\leftarrow B}\circ \mcu^\eud_{B\leftarrow b},$$

one sees that $(\star)$ is in turn equivalent to

$$ \mcu^{\eue}_{B\leftarrow A} \circ \mcu^{\eud}_{A\leftarrow a}=  \mcu^{\eud}_{B\leftarrow b}\circ \mcu^{\eue}_{b\leftarrow a}.$$

Since $A\in \mcu^{\eud}_a=\mcs^{-\eud}_a$ and $b\in\mcu^\eue_a\subseteq \mcu^{-\eud}_a$, applying \autoref{118} with $\euc$ replaced with $-\eud$ guarantees this equality.

\done
\end{pf}

%% file: uzman_arithmeticity_affineextension.tex
\section{Measurable Covering Map and Diagonal Affine Extension\label{135}}

In this section we put together the affine structures as well as the affine holonomies we have established in the previous sections to obtain a diagonal affine extension of $(\mu,\alpha)$. We shall use invariant language to describe this affine extension; so that its phase space is $TM$ considered as an \ae-bundle over $M$.

\begin{rem}\label{128}

Let $\euc\in \Cham(\mu,\alpha)$. Note that for $x\in_\mu$ we have $T_xM \cong \RR^k \times \left(S^\euc_x \oplus U^\euc_x\right)$ as vector spaces. Let us denote for $v\in T_x M$ the corresponding components by $v=(v^o,v^s,v^u)$. Here and in what follows we use a direct product between the orbit and normal to the orbit directions instead direct sums (as opposed to \autoref{031}); we will indeed recover a certain solvable group structure on $T_xM$ where the twist involved in the semidirect product occurs exactly between the orbit directions $O_x\cong \RR^k$ and the directions $S^\euc_x\oplus U^\euc_x$ normal to the orbit  (see \autoref{150} below). More precisely we want a \emph{split} short exact sequence (in $\cat{Lie}$)

\begin{center}
\begin{tikzcd}
S^\euc_x\oplus U^\euc_x \arrow[r, tail] & T_xM \arrow[r, two heads]                  & \RR^k \\
                                        & \RR^k \arrow[u] \arrow[ru, "\id_{\RR^k}"'] &      
\end{tikzcd}

\end{center}

where the action $\RR^k\lact S^\euc_x\oplus U^\euc_x$ is \emph{not} trivial.

\done
\end{rem}

\subsection{Measurable Covering Map}

\begin{prp}\label{126}

Let $\euc\in\Cham(\mu,\alpha)$. Then for $x\in_\mu M$, the following formula defines a measurable map that is \ae-defined w/r/t the Haar measure class on $T_xM$:

\begin{align*}
\Phi_x=\Phi^\euc_x: &(T_xM,0)\cong (\RR^k\times (\mcs^\euc_x\oplus\mcu^\euc_x),(0,0,0)) \to (M,x),\\
 &v=(v^o,v^s,v^u)\mapsto \mcu^\euc_{\Upsilon^\euc_x(v^u)\leftarrow x}\circ \Gamma\Sigma^\euc_x(v^o,v^s).
\end{align*}

Moreover,

\begin{enumerate}

\item $\Phi:TM\to M$, $(x,v)\mapsto \Phi_x(v)$ is measurable.

\item The same map is also \ae-defined by the following formulas:

\begin{align*}
&\all x\in_\mu M, \all (v^o,v^s,v^u)\in_{\aem}T_xM:\\
&\phantom{\all x\in_\mu M, }\Phi^\euc_x(v^o,v^s,v^u)
=\alpha_{v^o}\circ \restr{\left[\mcu^\euc_{\Gamma\Upsilon^\euc_x(-v^o,v^u)\leftarrow \alpha_{-v^o}(x)}\right]}{\mcs^\euc_{\alpha_{-v^o}(x)}}\circ \Sigma^\euc_x(v^s)\\
&\phantom{\all x\in_\mu M, \Phi^\euc_x(v^o,v^s,v^u)}=\alpha_{v^o}\circ \mcu^\euc_{\Upsilon^\euc_x(v^u)\leftarrow x}\circ \Sigma^\euc_x(v^s)\\
&\phantom{\all x\in_\mu M, \Phi^\euc_x(v^o,v^s,v^u)}= \mcs^\euc_{\Sigma^\euc_x(v^s)\leftarrow x}\circ \Gamma\Upsilon^\euc_x(v^o,v^u)\\
&\phantom{\all x\in_\mu M, \Phi^\euc_x(v^o,v^s,v^u)}= \alpha_{v^o}\circ \restr{\left[\mcs^\euc_{\Gamma\Sigma^\euc_x(-v^o,v^s)\leftarrow \alpha_{-v^o}(x)}\right]}{U^\euc_{\alpha_{-v^o}(x)}}\circ \Upsilon^\euc_x(v^u)\\
&\phantom{\all x\in_\mu M, \Phi^\euc_x(v^o,v^s,v^u)}= \alpha_{v^o}\circ \mcs^\euc_{\Sigma^\euc_x(v^s)\leftarrow x}\circ \Upsilon^\euc_x(v^u).
\end{align*}

\item For any $x\in_\mu M$, $\Phi_x$ has a diagonal property: each coordinate line according to the splitting $T_xM = O_x\times \left(\bigoplus_{\chi\in\LSpec^\ast(\mu,\alpha)}L^\chi_x\right)$ is either mapped onto the leaf of a $1$-dimensional subfoliation of $\mco_x$ given by one of the infinitesimal generators of the action $\alpha$ xor to a $1$-dimensional Lyapunov manifold $\mcl_x^\chi$.  

\item For any $x\in_\mu M$,  $\Phi_x$ is independent of the chamber $\euc$. That is, if $\eud\in\Cham(\mu,\alpha)$ is an anonymous chamber, then

$$\all x\in_\mu M, \all v\in_{\aem} T_xM: \Phi^\euc_x(v)=\Phi^\eud_x(v).$$

\item We have for any $x\in_\mu M$, for any $v\in_{\aem} T_xM$ and for any chamber $\eud\in\Cham(\mu,\alpha)$:

\begin{align*}
v\in O_x &\implies \Phi_x(v)=\alpha_v(x),\\
v\in S^\eud_x &\implies \Phi_x(v)=\Sigma^\eud_x(v)\\
v\in U^\eud_x &\implies \Phi_x(v)=\Upsilon^\eud_x(v).\\
\end{align*}

\end{enumerate}

\done
\end{prp}

\begin{pf}

We have all the ingredients established already. That $\Phi_x=\Phi^\euc_x$ is well-defined  follows from \autoref{111}, \autoref{117} and \autoref{118}. That both the basepointed version $\Phi_x$ and the global version $\Phi$ are measurable follows from item (iii) of \autoref{109}, item (iv) of \autoref{111} and item (i) of \autoref{127}. Item (ii) is syntactic (see also \autoref{143}). Items (iii) and (iv) follow from the assembly properties \autoref{099} together with the observation that all diagonal affine holonomies can be ultimately decomposed as affine holonomies along Lyapunov \ae-foliations. Finally item (v) follows from the previous three items.

\done
\end{pf}

Next we establish further properties of $\Phi_x$ that will justify considering it as a measurable covering map. Let $\euc\in\Cham(\mu,\alpha)$ and define $\all x\in_\mu M:$ $O^\perp_x=S^\euc_x\oplus U^\euc_x$; $O^\perp_x$ is independent of the choice of the chamber $\euc$ (indeed it is canonically isomorphic to $\faktor{T_xM}{O_x}$ as a vector space). We also write $v=(v^o,v^\perp)\in O_x\times O^\perp_x$.

\begin{prp}\label{131}

For $x\in_\mu M$, let $\Phi_x:T_xM\to M$ be as in \autoref{126}. Denote for any $\delta\in\RRP$ by $\eta_{x,\delta}$ the Haar measure on $T_xM$ with $\eta_{x,\delta}(O_x[0|\leq\delta]\times O_x^\perp[0|\leq\delta])=1$. Then for any $x\in_\mu M$ and for any $\eps\in]0,1[$, there is a $\delta\in\RRP$, $c\in\RRP$ and a measurable $B_{x,\eps}\subseteq O_x[0|\leq\delta]\times O_x^\perp[0|\leq\delta]$ such that

\begin{enumerate}

\item $0\in B_{x,\eps}$,

\item $1-\eps<\eta_{x,\delta}(B_{x,\eps})$ and $1-\eps<(\eta_{x,\delta})^{O_x}_0(B_{x,\eps}\cap \{0\}\times O^\perp_x)$, where we disintegrate $\eta_{x,\delta}$ along $T_xM\to O_x$\fn{See \autoref{090}.}.

\item $\restr{\Phi_x}{B_{x,\eps}}: (B_{x,\eps},0)\to (M,x)$ is injective,

\item $\fwd{\Phi_x}(\restr{\eta_{x,\delta}}{B_{x,\eps}})=c\, \restr{\mu}{\fwd{\Phi_x}(B_{x,\eps})}$.

\item Furthermore we have $\mu\left(\fwd{\Phi_x}(T_xM)\right)=1$.

\end{enumerate}

\done
\end{prp}

\begin{pf}

Let $\eps\in]0,1[$. By \autoref{108} and the Lebesgue Density Theorem, for $x\in_\mu M$, $\exi \Lambda\in\RRP$ such that there is a $\delta'\in\RRP$ with

$$1-\eps<\mu(M[x|\leq \delta']\, |\, \Pi(\eps,\Lambda));$$

here we condition $\mu$ on the Lusin-Pesin set $\Pi(\eps,\Lambda)$. By items (iv) and (ii) in \autoref{109} and \autoref{127}, respectively, and the absolute continuity of holonomies on Pesin sets, fixing a chamber $\euc$ and defining

$$B_{x,\eps}=\bck{\Phi_x}(\mco_{x,\loc}\times \mcs_{x,\loc}^\euc\times \mcu_{x,\loc}^\euc\cap \Pi(\eps,\Lambda))$$

with the sizes of local leaves not larger than some $\delta\in\RRP$ gives the first three items. Note that $\fwd{\Phi_x}(B_{x,\eps})$ is measurable by the continuity of $\Phi$ on Pesin sets, item (iii) and the Lusin-Souslin Theorem\fn{See \cite[p.89, Thm.15.1]{MR1321597}.}; hence item (iv) follows from \autoref{129}. For item (v) note first that $\fwd{\Phi_x}(T_xM)$ is $\alpha$-invariant; thus it suffices by ergodicity to verify that it's measurable. That it's indeed measurable follows  from items (iii) and (iv), together with \autoref{146}.

\done
\end{pf}

\begin{dfn}\label{130}

In light of \autoref{131}, we call the measurable map $\Phi_x: (T_xM,0)\to (M,x)$ defined in \autoref{126} the \textbf{measurable covering} of the system $(\mu,\alpha)$ at $x\in_\mu M$.


\done
\end{dfn}

\subsection{Diagonal Affine Extension}

Next we discuss the effect of changing the basepoints on $\Phi_x$. One can change the basepoint either horizontally or vertically (or both simultaneously). The horizontal basepoint change is the change of $x$ along the manifold $M$. The vertical basepoint change is the change of the marked origin $0\in T_xM$ to some other vector; note that we consider $T_xM$ as an affine space. Of course under the measurable covering map $\Phi_x$ vertical basepoint changes transform into horizontal basepoint changes. We start with special cases.

\begin{prp}\label{132}

Let $\euc\in\Cham(\mu,\alpha)$ and $x\in_\mu M$. Then

\begin{enumerate}

\item For $v^s\in_{\aem} S^\euc_x$ put $y=\Phi_x(v^s)\in\mcs^\euc_x$. Then there is a unique affine isomorphism $\Phi_{(y,0)\ot (x,v^s)}\in \DAff(T_xM; T_yM)$ such that $\Phi_{(y,0)\ot (x,v^s)}(v^s)=0$ and $\Phi_{y}\circ\Phi_{(y,0)\ot (x,v^s)}=_{\aem} \Phi_x$ on $T_xM$. 

\item Similarly for $v^u\in_{\aem}U^\euc_x$, putting $z=\Phi_x(v^u)\in\mcu^\euc_x$, there is a unique affine isomorphism $\Phi_{(z,0)\ot (x,v^u)}\in \DAff(T_xM;T_zM)$ such that $\Phi_{(z,0)\ot (x,v^u)}(v^u)=0$ and $\Phi_{z}\circ\Phi_{(z,0)\ot (x,v^u)}=_{\aem} \Phi_x$ on $T_xM$.

\end{enumerate}

\done
\end{prp}

\begin{pf}

The proof is a matter of unfolding and transforming the definitions. We only verify the first item and merely give the analogous formula for the second item for future use. Let $w=(w^o,w^s,w^u)\in_{\aem} T_xM$. Then

\begin{align*}
\Phi_x(w)
&=\alpha_{w^o}\circ\mcs^\euc_{\Sigma^\euc_x(w^s)\ot x}\circ \Upsilon^\euc_x(w^u)\\
&=\alpha_{w^o}\circ\mcs^\euc_{\Sigma^\euc_y\left(\left(\Sigma^\euc_y\right)^{-1}\circ \Sigma^\euc_x(w^s)\right)\ot y}\circ \mcs^\euc_{y\ot x}\circ \Upsilon^\euc_x(w^u)\\
&=\alpha_{w^o}\circ\mcs^\euc_{\Sigma^\euc_y\left(\left(\Sigma^\euc_y\right)^{-1}\circ \Sigma^\euc_x(w^s)\right)\ot y}\circ \Upsilon^\euc_y \left(\left(\Upsilon^\euc_y\right)^{-1}\circ\mcs^\euc_{y\ot x}\circ \Upsilon^\euc_x(w^u)\right)\\
&=\Phi_y\left(w^o,\left(\Sigma^\euc_y\right)^{-1}\circ \Sigma^\euc_x(w^s), \left(\Upsilon^\euc_y\right)^{-1}\circ\mcs^\euc_{y\ot x}\circ \Upsilon^\euc_x(w^u)\right)\\
&=\Phi_y\circ \left[I_k\times\left(\left(\Sigma^\euc_y\right)^{-1}\circ \Sigma^\euc_x\right)\times \left(\left(\Upsilon^\euc_y\right)^{-1}\circ\mcs^\euc_{y\ot x}\circ \Upsilon^\euc_x\right)\right](w),
\end{align*}

so that the following formula defines the affine isomorphism $\Phi_{(y,0)\ot (x,v^s)}$ with the required properties \ae-uniquely:

\begin{align*}
&\Phi_{(y,0)\ot (x,v^s)}=I_k\times\left(\left(\Sigma^\euc_y\right)^{-1}\circ \Sigma^\euc_x\right)\times \left(\left(\Upsilon^\euc_y\right)^{-1}\circ\mcs^\euc_{y\ot x}\circ \Upsilon^\euc_x\right) \quad\quad (\star^s)\\
&\phantom{\Phi_{(y,0)\ot (x,v^s)}=I_k\times\left(\left(\Sigma^\euc_y\right)^{-1}\circ \Sigma^\euc_x\right)\times} \text{ for } y=\Phi_x(v^s), v^s\in_{\aem} S^\euc_x.
\end{align*}

Similarly we have

\begin{align*}
&\Phi_{(z,0)\ot (x,v^u)}=I_k\times\left(\left(\Sigma^\euc_z\right)^{-1}\circ\mcu^\euc_{z\ot x}\circ \Sigma^\euc_x\right) \times \left(\left(\Upsilon^\euc_z\right)^{-1}\circ \Upsilon^\euc_x\right) \quad\quad (\star^u)\\
&\phantom{\Phi_{(y,0)\ot (x,v^s)}=I_k\times\left(\left(\Sigma^\euc_y\right)^{-1}\circ \Sigma^\euc_x\right)\times} \text{ for } z=\Phi_x(v^u), v^u\in_{\aem} U^\euc_x.
\end{align*}

\done
\end{pf}

In \autoref{132} we have discussed the effect on $\Phi$ of changing the basepoint vertically along the stable and unstable directions of a chamber; next we discuss the effect of changing the basepoint along the orbit directions. Note that since $\all x\in_\mu M,$ $\RR^k\cong O_x\cong \mco_x$ canonically, we may interpret this change of basepoint either vertically or horizontally, without reference to the measurable covering map $\Phi_x$.

\begin{prp}\label{134}

Let $x\in_\mu M$, $v^o\in O_x$ and put $y=\Phi_x(v^o)\in\mco_x$. Then there is a unique affine isomorphism $\Phi_{(y,0)\ot (x,v^o)}\in \DAff(T_xM; T_yM)$ such that $\Phi_{(y,0)\ot (x,v^o)}(v^o)=0$ and $\Phi_{y}\circ\Phi_{(y,0)\ot (x,v^o)}=_{\aem} \Phi_x$ on $T_xM$.

\done
\end{prp}

\begin{pf}

Let us fix a chamber $\euc\in\Cham(\mu,\alpha)$ and consider the cocycles $\phi^{\Sigma^\euc}$ and $\phi^{\Upsilon^\euc}$ attached to the stable and unstable affine parameters of $\euc$ in item (ii) of \autoref{109}. Then arguing as in \autoref{132} we have for $w=(w^o,w^s,w^u)\in_{\aem} T_xM$:

$$\Phi_x(w)=\Phi_{y}\left(w^o-v^o,\phi^{\Sigma^\euc}(v^o,x)(w^s),\phi^{\Upsilon^\euc}(v^o,x)(w^u)\right).$$

Thus the following formula defines the affine isomorphism $\Phi_{(y,0)\ot (x,v^o)}$ with the desired properties:

\begin{align*}
&\Phi_{(y,0)\ot (x,v^o)} = \left(I_k\times \phi^{\Sigma^\euc}(v^o,x)\times \phi^{\Upsilon^\euc}(v^o,x)\right)-v^o \quad\quad (\star^o)\\
&\phantom{\Phi_{(y,0)\ot (x,v^o)} = \left(I_k\times \phi^{\Sigma^\euc}(v^o,x)\times\right.}\text{ for } y=\Phi_x(v^o), v^o\in O_x.
\end{align*}

Note that by virtue of the assembly properties \autoref{099}, this formula defines an affine isomorphism independently of the choice of the chamber $\euc$.

\done
\end{pf}

\begin{cor}\label{147}

Let $x\in_\mu M$ and $t\in\RR^k$. Then

$$\alpha_t\circ \Phi_x =_{\aem} \pm\, \Phi_{\alpha_t(x)}\circ T_x\alpha_t.$$

Consequently $\Phi: TM\to M$ displays $(\mu,\alpha)$ as the factor map of some system that acts via diagonal isomorphisms fiberwise.

\done
\end{cor}

\begin{pf}

The canonical identification $\all x\in_\mu M,$ $\RR^k\cong O_x\cong \mco_x$ explicitly refers to the fact that $\left(\alpha_\blt(\alpha_t(x))\right)^{-1}\circ \alpha_t\circ \alpha_\blt(x)=\id_{\RR^k}$, that is

\begin{center}
\begin{tikzcd}
\mco_x \arrow[r, "\alpha_t"]                                & \mco_{\alpha_t(x)}                           \\
\RR^k \arrow[u, "\alpha_\blt(x)"] \arrow[r, "\id_{\RR^k}"'] & \RR^k \arrow[u, "\alpha_\blt(\alpha_t(x))"']
\end{tikzcd}

\end{center}

Thus $T_x\alpha_t$ fixes any vector in $O_x$, when the Euclidean coordinates are used. By \autoref{115} we also have that the cocycles $\phi^{\Sigma^\euc}$ and $\phi^{\Upsilon^\euc}$ coincide with the derivative cocycles of $\alpha$ restricted to the stable and unstable \ae-subbundles of the chamber $\euc$ up to a sign, respectively; thus $(\star^o)$ in the proof \autoref{134} immediately gives the statement with $v^o$ replaced by $t$. The extension of $(\mu,\alpha)$ that $\Phi:TM\to M$ is a factor map that coincides with the derivative cocycle up to a sign that depends on $x\in_\mu M$ and $t\in\RR^k$.

\done
\end{pf}

\begin{prp}\label{133}

Let $x,y\in_\mu M$, $v\in_{\aem} T_xM$ and $w\in_{\aem} T_y M$. If $p$ $=\Phi_x(v)$ $=\Phi_y(w)$, then there is a unique affine isomorphism $\Phi_{(y,w)\ot (x,v)}$ $\in \DAff(T_xM;T_yM)$ with $\Phi_{(y,w)\ot (x,v)}(v)$ $=w$ such that $\Phi_y\circ \Phi_{(y,w)\ot (x,v)} =_{\aem} \Phi_x$, that is,

\begin{center}

\begin{tikzcd}
{(T_xM,v)} \arrow[rdd, "\Phi_x"'] \arrow[rr, "{\Phi_{(y,w)\ot(x,v)}}"] &         & {(T_yM,w)} \arrow[ldd, "\Phi_y"] \\
                                                                                 & \aem    &                                       \\
                                                                                 & {(M,p)} &                                      
\end{tikzcd}

\end{center}

\done
\end{prp}

\begin{pf}

The proof is based on the previous two propositions \autoref{132} and \autoref{134}; together with good choices of points that are to be abbreviated. Let us first establish the affine isomorphism $\Phi_{(p,0)\ot(x,v)}: (T_xM,v)\to (T_pM,0)$. Let $v=(v^o,v^s,v^u)$ and let us put 

$$q=\alpha_{-v^o}(p),r=\Phi_x(v^s)\in\mcs^\euc_x,z^u=\left(\Phi_r\right)^{-1}(q)\in U^\euc_r.$$

We have by \autoref{134} that $\Phi_q =_{\aem} \Phi_p\circ \Phi_{(p,0)\ot(q,v^o)}$. Since we also have

$$
q=\Phi_x(0,v^s,v^u)
= \mcu^\euc_{\Upsilon^\euc_x(v^u)\leftarrow x}\circ \Sigma^\euc_x(v^s)
=\mcu^\euc_{\Upsilon^\euc_x(v^u)\leftarrow x}(r)
\in\mcu^\euc_r,$$

by \autoref{132} we have $\Phi_r =_{\aem} \Phi_q\circ \Phi_{(q,0)\ot (r,z^u)}$. Again by \autoref{132} we have $\Phi_x=_{\aem} \Phi_r\circ \Phi_{(r,0)\ot (x,v^s)}$; thus we may \ae-uniquely define the affine isomorphism $\Phi_{(p,0)\ot (x,v)}:T_xM\to T_pM$ by

$$\Phi_{(p,0)\ot (x,v)} =_{\aem} \Phi_{(p,0)\ot(q,v^o)}\circ \Phi_{(q,0)\ot (r,z^u)}\circ \Phi_{(r,0)\ot (x,v^s)}.$$

Composing the formulas $(\star^s)$, $(\star^u)$ and $(\star^o)$ in the proofs of \autoref{132} and \autoref{134} with appropriate point substitutions we obtain:

\begin{align*}
\Phi_{(p,0)\ot (x,v)} 
&=_{\aem} I_k \\
&\phantom{=_{\aem}}\times \left(\phi^{\Sigma^\euc}(v^o,q)\circ \left(\Sigma^\euc_q\right)^{-1} \circ \mcu^\euc_{q\ot r}\circ \Sigma^\euc_x\right)\\
&\phantom{=_{\aem}}\times \left( \phi^{\Upsilon^\euc}(v^o,q)\circ \left(\Upsilon^\euc_q\right)^{-1}\circ \mcs^\euc_{r\ot x}\circ \Upsilon^\euc_x\right)-v^o.\\
&\phantom{=_{\aem}\times}\text{ for }p=\Phi_x(v), q=\alpha_{-v^o}(p),r=\Phi_x(v^s), v=(v^o,v^s,v^u)\in_{\aem} T_xM.
\end{align*}

With this formula it's also straightforward that indeed $\Phi_{(p,0)\ot (x,v)}(v)=0$. Finally we may define the affine isomorphism

$$\Phi_{(y,w)\ot (x,v)}=\left(\Phi_{(p,0)\ot (y,w)}\right)^{-1}\circ \Phi_{(p,0)\ot (x,v)}:(T_xM,v)\mapsto (T_yM,w).$$

\done
\end{pf}

%% file: uzman_arithmeticity_homoclinicgroup.tex
\section{Homoclinic Group and Conclusion of Proof of the Main Theorem\label{069}}

In this section we introduce the homoclinic groupoid and the homoclinic group and conclude the proof of \autoref{001}. In the literature homoclinic groups are either defined already with reference to algebraic structure (i.e. for actions on groups by automorphisms)\fn{See e.g. \cite{MR1359974,MR1738686,MR2013356,MR2275445,MR3082539,MR3314515,MR4366230}.}, or else with reference to a Smale structure (i.e. an abstract local product structure)\fn{See e.g. \cite{MR415675,MR943921,MR2667385,MR2586354}.}. Our approach is closer to the latter perspective; the measurable covering map $\Phi_x$ we constructed is a measurable analog of a global product structure and the recovery of the algebraic structure is essentially based on the fact that the two approaches coincide. Let us also note that often homoclinic groups are considered for discrete group systems. In fact the author was unable to find a case study for the homoclinic group of a non-discrete group system in the literature.

Note that by \autoref{133}, for $x\in_\mu M$ and for $v,w\in_{\aem} T_xM$, if $\Phi_x(v)=\Phi_x(w)$, then there is a unique affine isomorphism $\Phi_{(x,w)\ot(x,v)}:(T_xM,v)\to (T_xM,w)$ such that $\Phi_x\circ \Phi_{(x,w)\ot(x,v)}=_{\aem} \Phi_x$, thus in light of \autoref{131} we may consider $\Phi_{(x,w)\ot(x,v)}$ as a measurable deck transformation of the measurable covering map $\Phi_x$. 

\begin{dfn}\label{145}

We define the \textbf{homoclinic groupoid} $\mfH$ of $(\mu,\alpha)$ by 

$$\mfH_{y\ot x}=\{A\in \DAff(T_xM;T_yM)\,|\, \Phi_y\circ A=_{\aem} \Phi_x\}$$

for $x,y\in_\mu M$. Similarly for $x\in_\mu M$, we define the \textbf{ homoclinic group} $\mfH_x$ of $(\mu,\alpha)$ at $x$ by

$$\mfH_x=\mfH_{x\ot x}=\{A\in\DAff(T_xM)\,|\, \Phi_x\circ A=_{\aem }\Phi_x\}.$$

\done
\end{dfn}

One can see the reasoning behind calling $\mfH_x$ the homoclinic group as follows: fix a chamber $\euc\in\Cham(\mu,\alpha)$, and let for $x\in_\mu M$,  $y\in\mcs^\euc_x\cap \mcu^\euc_x$, so that $y$ is homoclinic to $x$ w/r/t $\alpha_t$ for any $t\in\euc\cup-\euc$. But if $y\neq x$, then we have two distinct vectors $v^s,v^u\in T_xM$ that are both sent to $y$ via the measurable covering map $\Phi_x$, namely, $v^s=\left(\Sigma^\euc_x\right)^{-1}(y)\in S^\euc_x$, and $v^u=\left(\Upsilon^\euc_x\right)^{-1}(y)\in U^\euc_x$. Then by \autoref{133}, we have a unique element $A=A_y\in \mfH_x$ such that $A:v^s\mapsto v^u$.

Let us denote by $\LL: \DAff(T_xM)\to \DGL(T_xM)$ the homomorphism that takes an affine automorphism to its linear part, and $\II:T_xM \to \DAff(T_xM)$ be the inclusion of the kernel of $\LL$. Denote by $\mfT\mfH_x=\bck{\II}(\mfH_x)\cong\ker(\LL)\cap \mfH_x$ the normal subgroup of $\mfH_x$ that is the translation part. We shall write $w$ for an element of $\mfT\mfH_x$ when we want to emphasize that $w$ is a vector in $T_xM$ and $\II(w)$ when we want to emphasize that $w$ is a symmetry of $\Phi_x$.

\begin{obs}\label{148}

Let $x\in_\mu M$, $t\in \RR^k$. Immediately by \autoref{147} we have that $T_x\alpha_t:T_xM\to T_{\alpha_t(x)}M$ conjugates the homoclinic group at $x$ and $\alpha_t(x)$. More precisely;

$$\all A\in \mfH_x: T_x\alpha_t\circ A\circ \left(T_x\alpha_t\right)^{-1}\in \mfH_{\alpha_t(x)}.$$

It's also straightforward that $\all w\in \mfT\mfH_x:$ $T_x\alpha_t\circ \II(w)\circ \left(T_x\alpha_t\right)^{-1}$ $=\II\circ  T_x\alpha_t(w)$, thus

$$\fwd{T_x\alpha_t}(\mfT\mfH_x)=\mfT\mfH_{\alpha_t(x)}.$$

\done
\end{obs}

By \autoref{148}, the linear span $\Span(\mfT\mfH_\blt)\leq TM$ of the translation part $\mfT\mfH_\blt$ of the family of the homoclinic groups is $\Ad^\alpha$-invariant; it's also straightforward that $\Span(\mfT\mfH_\blt)$ is a measurable\fn{See the proof of \autoref{153} for an argument for measurability.} polarization. Thus by ergodicity $\dim(\Span(\mfT\mfH_\blt)):M\to \UL{2k+2}$ is \ae-constant. Note that by essential freeness of $(\mu,\alpha)$ and \autoref{134}, for $x\in_\mu M, $ $\dim(\Span(\mfT\mfH_x))\cap O_x=0$, that is, $\Span(\mfT\mfH_x)\leq O_x^\perp$; we'll show that indeed $\dim(\Span(\mfT\mfH_\blt))=_\mu k+1$. We split the proof into two; first we'll show that $\dim(\Span(\mfT\mfH_\blt))$ is not degenerate and then we'll show that its rank is maximal, subject to the complementarity to the orbit directions $O_\blt$.

\begin{lem}\label{151}

For $x\in_\mu M$, $\mfT\mfH_x$ is a discrete subgroup of $O_x^\perp$. Here we consider $O_x^\perp$ with its standard abelian group structure.

\done
\end{lem}

\begin{pf}

Suppose otherwise. Then there is a sequence $w_\blt\subseteq \mfT\mfH_x$ of vectors such that $\lim_{n
\to \infty} w_n$ $= 0$. Fix $\eps\in]0,1[$, and let $B_{x,\eps}$ be as in \autoref{131}. If $v\in B_{x,\eps}\cap (B_{x,\eps}+w_n)$, then $v,v-w_n=\II(-w_n)(v)\in B_{x,\eps}$. But then $\Phi_x(v)=\Phi_x\circ \II(w_n)(v),$ so that $(\eta_{x,\delta})^{O_x}_0(B_{x,\eps}\cap (B_{x,\eps}+w_n))=0$ (recall that $(\eta_{x,\delta})^{O_x}_0$ is the fiber measure at $0$ of a Haar measure on $T_xM$ along the projection $T_xM\to O_x$). We also have $\lim_{n\to \infty}(\eta_{x,\delta})^{O_x}_0(B_{x,\eps}\cap (B_{x,\eps}+w_n))=(\eta_{x,\delta})^{O_x}_0(B_{x,\eps})>1-\eps>0$, a contradiction.

\done
\end{pf}

\begin{lem}\label{152}

For $x\in_\mu M$: $\dim(\Span(\mfT\mfH_x))>0$.

\done
\end{lem}

\begin{pf}

Suppose to the contrary $\dim(\Span(\mfT\mfH_\blt))=_\mu 0$, so that for $x\in_\mu M$, $\mfH_x$ has no translation part. We'll see that this implies that $\mfH_x$ is trivial, which in turn implies by \autoref{133} that $\Phi_x: T_xM \to M$ is \ae-injective. Conjugating $(\mu,\alpha)$ via $\Phi_x$; we have by \autoref{147},

$$\left(\Phi_x\right)^{-1}\circ \alpha_t\circ \Phi_x=_{\aem} \pm \left(\Phi_x\right)^{-1}\circ \Phi_{\alpha_t(x)}\circ T_x\alpha_t, $$

so that the conjugate system on $T_xM$ is by affine automorphisms. By item (iv) of \autoref{131} the transformed probability measure is of Haar class; but no affine automorphism of a vector space admits a \emph{probability} measure of positive entropy\fn{Indeed, by Poincar\'{e} Recurrence the support of any probability measure invariant under an affine automorphism is in the closure of set of recurrent points; the linear part of the affine automorphism restricted to this set must have all eigenvalues of modulus $1$, hence on the support of the invariant measure trajectories can separate at most polynomially.}, a contradiction.

Thus it suffices to show that $\mfH_x=\{\id_{T_xM}\}$. Suppose not. For the sake of readability we'll fix a basis $T_xM\cong \RR^{2k+1}$ and write e.g. $(A,a)$ for the affine automorphism $v\mapsto Av+a$. Since $\mfT\mfH_x$ is trivial, $\LL: \mfH_x\hookrightarrow \DGL(T_xM)$ is an embedding; in particular $\mfH_x$ is abelian. Then for $(A,a),(B,b)\in\mfH_x$, the abelian nature of $\mfH_x$ implies that $(A-I)(b)=(B-I)(a)$. If either of $A$ or $B$ is $I$ this equation is vacuous; for $A\neq I \neq B$ however we have $(B-I)^{-1}(b)=(A-I)^{-1}(a)$, so that the set

$$\{-(A-I)^{-1}(a)\,|\, (A,a)\in \mfH_x, A\neq I\}$$

is a singleton; let us denote by $v^\ast_x$ its unique element. A straightforward computation shows that $v^\ast_x\in \Fix(\mfH_x)= \{v\in T_xM \,|\, \all A\in\mfH_x: A(v)=v\}$, or alternatively that the affine automorphism $(I,v^\ast_x)\in \DAff(T_xM)$ conjugates the affine action $\mfH_x\lact T_xM$ to the linear action $\fwd{\LL}(\mfH_x)\lact T_xM$. Let us put $F_x=\Fix(\mfH_x)$ and $F_x^0=\Fix(\fwd{\LL}(\mfH_x))$, so that we have

$$F_x= F_x^0+v^\ast_x.$$

Since elements in $\DGL(T_xM)$ are diagonal w/r/t the Oseledets splitting $T_xM = O_x \oplus \bigoplus_{\chi} L^\chi_x$, $F_x^0$ is a direct sum of (possibly some) orbit directions and (possibly) some Lyapunov subspaces. In particular both $F_\blt$ and $F_\blt^0$ are measurable and $\Ad^\alpha$ invariant by \autoref{152}. Further, by item (v) of \autoref{126} and \autoref{114}, $O_x\leq F_x^0$. Hence there is a unique vector $v^\dagger_x\in T_xM$ that is an element of the sum of the Lyapunov subspaces that are not included in $F_x^0$ such that

$$F_x=F_x^0+v^\dagger_x.$$

$v^\dagger_\blt: M\to TM$ is an \ae-defined measurable vector field that is $\Ad^\alpha$-invariant. If it didn't vanish $\mu$-\ae, we could find a $t^\dagger\in\RR^k$ such that $n\mapsto |v_{\alpha_{nt^\dagger}(x)}|$  grows exponentially; which contradicts Poincar\'{e} Recurrence; whence $v^\dagger_\blt=_\mu 0$, so $F_\blt=_\mu F_\blt^0$, and $0\in F_x$ for $x\in_\mu M$. Thus if $(A,a)\in \mfH_x$, then $a=0$; i.e. $\mfH_x$ is \emph{equal} to $\LL(\mfH_x)$ under the standard embedding $\DGL(T_xM)\hookrightarrow \DAff(T_xM)$, $A\mapsto (A,0)$.

Let us fix a chamber $\euc\in\Cham(\mu,\alpha)$. Then any $A\in\mfH_x$ preserves $S^\euc_x$. Let $v^s\in S^\euc_x$, fix $\eps\in ]0,1[$ and let $B_{\blt,\eps}$ be as in \autoref{131}. Then there is a $t^\star\in \euc$ such that $T_x\alpha_{t^\star}(v^s), T_x\alpha_{t^\star}\circ A(v^s)\in B_{\alpha_{t^\star}(x),\eps}$; by \autoref{148}, the injectivity of $\Phi_{\alpha_{t^\star}(x)}$ on $B_{\alpha_{t^\star}(x),\eps}$ forces $A(v^s)=v^s$; that is to say $S^\euc_x\subseteq F_x$. Similarly $U^\euc_x\subseteq F_x$, so that $F_x= T_xM$. Thus $\mfH_x$ is trivial, a contradiction.

\done
\end{pf}

\begin{lem}\label{153}

For $x\in_\mu M:$  $\Span(\mfT\mfH_x)=O_x^\perp$ and consequently also $\mfT\mfH_x\cong \ZZ^{k+1}$.

\done
\end{lem}

\begin{pf}

Suppose otherwise, so that, after \autoref{152}, for $x\in_\mu M$, $0<\dim(\Span(\mfT\mfH_x))<k+1$. Let us denote by $d$ the $\mu$-\ae constant value of $\dim(\Span(\mfT\mfH_\blt))$. For $x\in_\mu M$, let us denote by $\omega_x$ the volume of the $d$-dimensional torus $\faktor{\Span(\mfT\mfH_x)}{\mfT\mfH_x}$; recall that at the beginning of \autoref{112} we had fixed a $C^\infty$ Riemannian metric $\mfg$ and via this $\omega_x$ is well defined for $x\in_\mu M$. Note that $\omega_\blt: M\to \RRP$ is measurable; indeed we may fix a basepoint $x^\ast\in_\mu M$ and by item (v) of \autoref{131} and by the proof of \autoref{133} $\mfH_\blt$ is the image of $\mfH_{x^\ast}$ under a measurable family of diagonal affine maps; similarly $\Span(\mfT\mfH_\blt))$ and consequently the restriction of $\mfg$ to $\Span(\mfT\mfH_\blt))$ is also measurable. We may bunch the Lyapunov exponents of $(\mu,\alpha)$ via an argument similar to the one we used in the proof of \autoref{109} to obtain a time $t^\star$ such that $n\mapsto \omega_{\alpha_{nt^\star}(x)}$ grows exponentially by the Oseledets Theorem\fn{More specifically here we use the intermediate-dimensional version of \textbf{ASYM2} of \autoref{031}.}; this contradicts Poincar\'{e} Recurrence. That $\mfT\mfH_x$ is isomorphic to $\ZZ^{k+1}$ follows from \autoref{151}.

\done
\end{pf}

Next observe more generally that $O_x$ is $\mfH_x$ invariant. Indeed, let $A\in\mfH_x$ and $v^o\in O_x$. We may write $A(v^o,0) = (A^o(v^o)+a^o,a^\perp)$. If $a^\perp\neq 0$, then as in the final paragraph of the proof of \autoref{152}, choosing an appropriate $t^\star\in\RR^k$ we get a contradiction to the local injectivity of $\Phi_{\alpha_{t^\star}(x)}$.  There is a similar invariance property in the $O^\perp_x$ direction. Fix a chamber $\euc\in\Cham(\mu,\alpha)$ and let $A\in\mfH_x$. If $v^\perp=(v^s,v^u)\in_{\aem} O_x^\perp$, then by  \autoref{134} we have

\begin{align*}
\Phi_x(0,v^s,v^u) 
&= \Phi_x(a^o,A^s(v^s)+a^s, A^u(v^u)+a^u)\\
&=\Phi_x(0, \phi^{\Sigma^\euc}\circ A^s(v^s)+\phi^{\Sigma^\euc}(a^o,x)(a^s), \phi^{\Upsilon^\euc}\circ A^u(v^u)+\phi^{\Upsilon^\euc}(a^o,x)(a^u)).
\end{align*}

Factoring $\mfT\mfH_x$ out from $T_xM$ we have a $\TT^{k+1}$ factor coming from $O^\perp_x$ and further that $\RR^k\cong O_x$ acts on the factor $\faktor{T_xM}{\mfT\mfH_x}$. The return times of the $\RR^k$ action on the $\TT^{k+1}$ factor is a cocompact subgroup; note that identifying this subgroup as $\ZZ^k$ produces the constant time change $\kappa$ in the statement of \autoref{001}. Thus we get that $\mfH_x$ has a subgroup isomorphic to $\ZZ^k$ acts on the $\TT^{k+1}$ factor by affine automorphisms and by the above calculation the $\RR^k$ action is a suspension of the $\ZZ^k$ action on the $\TT^{k+1}$ factor. Furthermore since the Lyapunov hyperplanes are in general position this $\ZZ^k$ action on $\TT^{k+1}$ is maximal Cartan. Again using \autoref{131} we have that the implied subgroup $\ZZ^{k+1}\rtimes \ZZ^k$ is of finite index in $\mfH_x$ and indeed the factor $F=\faktor{\mfH_x}{\ZZ^{k+1}\rtimes \ZZ^k}$ is a finite group of automorphisms of $\TT^{k+1}$. Since the $\ZZ^k$ action on $\TT^{k+1}$ is maximal Cartan $F$ is either $\{I_{k+1}\}$ xor $\{\pm I_{k+1}\}$. We have just proved:

\begin{prp}\label{150}

Let $x\in_\mu M$. Then the homoclinic group $\mfH_x$ is isomorphic to $(\ZZ^{k+1}\rtimes F)\rtimes \ZZ^k$ for $F=\{I_{k+1}\}$ xor $F=\{\pm I_{k+1}\}$. and consequently $\faktor{T_xM}{\mfH_x}$ is a $T^{k+1}$ bundle over $\TT^{k}$, where $T^{k+1}$ is either $\TT^{k+1}$ xor the $\pm$-infratorus $\faktor{\TT^{k+1}}{\pm I_{k+1}}$.

\done
\end{prp}

By definition $\Phi_x$ descends to a measure theoretical isomorphism $\faktor{\Phi_x}{\mfH_x} :\left(\faktor{T_xM}{\mfH_x},\mfH_x\right)\to (M,x)$ whose inverse $\Phi_{(\mu,\alpha)}=\Phi_{(\mu,\alpha),x}=\left(\faktor{\Phi_x}{\mfH_x}\right)^{-1}$ transforms $\alpha$ to the suspension of an affine Cartan action and $\mu$ to the suspension measure induced by Haar measure on $T^{k+1}$. Applying a Journ\'{e} lemma by de la Llave\fn{\cite[pp.304-305, Thm.5.7; pp.312-313, Prop.5.13]{MR1194019}} gives the smoothness properties of $\Phi_{(\mu,\alpha)}$ and concludes the proof.